
\input amstex
\documentstyle{amsppt}

\magnification=1250
\hsize=6.5truein
\vsize=9truein

\font \smallrm=cmr10 at 10truept
\font \smallbf=cmbx10 at 10truept
\font \smallit=cmti10 at 10truept
 at 10truept
 at 10truept

\baselineskip = 12pt

\def \loongrightarrow {\relbar\joinrel\relbar\joinrel\rightarrow}
\def \llongrightarrow {\relbar\joinrel\relbar\joinrel\relbar\joinrel
\rightarrow}

\def \longhookrightarrow {\lhook\joinrel\relbar\joinrel\rightarrow}
\def \llonghookrightarrow
{\lhook\joinrel\relbar\joinrel\relbar\joinrel\relbar\joinrel\rightarrow}

\def \longtwoheadrightarrow {\relbar\joinrel\twoheadrightarrow}
\def \llongtwoheadrightarrow
{\relbar\joinrel\relbar\joinrel\relbar\joinrel\twoheadrightarrow}

\def \smallcirc {\, {\scriptstyle \circ} \,}
\def \aij {a_{ij}}

\def \N {\Bbb N}
\def \Z {\Bbb Z}
\def \Q {\Bbb Q}
\def \C {\Bbb C}

\def \tm {t^{-1}}
\def \qm {q^{-1}}
\def \Cq {\C(q)}

\def \ehat {\widehat E}
\def \fhat {\widehat F}
\def \ebar {\overline E}
\def \fbar {\overline F}

\def \otimeshat {\,\widehat\otimes\,}

\def \gerg {{\frak g}}
\def \gerh {{\frak h}}
\def \gerk {{\frak k}}
\def \gerK {{\frak K}}
\def \gert {{\frak t}}
\def \gern {{\frak n}}
\def \gerb {{\frak b}}
\def \gerE {{\frak E}}
\def \gerU {{\frak U}}
\def \gerF {{\frak F}}

\def \ghat {{\hat \gerg}}
\def \that {{\hat \gert}}
\def \hhat {{\hat \gerh}}

\def \calU {\hbox{$ {\Cal U} $}}
\def \calF {\hbox{$ {\Cal F} $}}
\def \calE {\hbox{$ {\Cal E} $}}

\def \ehat {\widehat E}
\def \fhat {\widehat F}

\def \phire {\Phi^{\text{re}}}

\def \phipre {\Phi_+^{\text{re}}}
\def \phipim {\Phi_+^{\text{im}}}
\def \phitildep {\widetilde{\Phi}_+}
\def \phitildepim {\widetilde{\Phi}_+^{\text{im}}}

\def \ug {U(\ghat)}
\def \uh {U \big( \hhat \big)}

\def \uqg {U_q(\ghat)}
\def \uqh {U_q \big( \hhat \big)}

\def \uqMg#1 {U_q^{\scriptscriptstyle#1}(\ghat)}
\def \gerUMg#1 {\gerU^{\scriptscriptstyle#1}(\ghat)}
\def \gerUunoMg#1 {\gerU_1^{\scriptscriptstyle #1}(\ghat)}
\def \gerUepsilonMg#1
   {\gerU_\varepsilon^{\scriptscriptstyle #1}(\ghat)}
\def \calUMg#1 {{\Cal U}^{\scriptscriptstyle #1}(\ghat)}
\def \calUunoMg#1 {{\Cal U}_1^{\scriptscriptstyle #1}(\ghat)}
\def \calUepsilonMg#1
   {{\Cal U}_\varepsilon^{\scriptscriptstyle #1}(\ghat)}

\def \uzM#1 {U_0^{\scriptscriptstyle#1}}
\def \calUzM#1 {{\Cal U}_0^{\scriptscriptstyle#1}}
\def \gerUzM#1 {\gerU_0^{\scriptscriptstyle#1}}

\def \um {U_{\scriptscriptstyle -}}
\def \up {U_{\scriptscriptstyle +}}
\def \calUm {{\Cal U}_{\scriptscriptstyle -}}
\def \gerUm {\gerU_{\scriptscriptstyle -}}
\def \calUp {{\Cal U}_{\scriptscriptstyle +}}
\def \gerUp {\gerU_{\scriptscriptstyle +}}

\def \uMbm#1 {U^{\scriptscriptstyle#1}_{\scriptscriptstyle \leq}}
\def \calUMbm#1
   {{\Cal U}^{\scriptscriptstyle#1}_{\scriptscriptstyle \leq}}
\def \gerUMbm#1
   {\gerU^{\scriptscriptstyle#1}_{\scriptscriptstyle \leq}}

\def \uMbp#1 {U^{\scriptscriptstyle#1}_{\scriptscriptstyle \geq}}
\def \calUMbp#1
   {{\Cal U}^{\scriptscriptstyle#1}_{\scriptscriptstyle \geq}}
\def \gerUMbp#1
   {\gerU^{\scriptscriptstyle#1}_{\scriptscriptstyle \geq}}

\def \uqMh#1 {U_q^{\scriptscriptstyle#1} \big( \hhat \big)}
\def \calUMh#1 {{\Cal U}^{\scriptscriptstyle #1} \big( \hhat \big)}
\def \calUunoMh#1
   {{\Cal U}_1^{\scriptscriptstyle #1} \big( \hhat \big)}
\def \calUepsilonMh#1
   {{\Cal U}_\varepsilon^{\scriptscriptstyle #1} \big(
\hhat \big)}
\def \gerUMh#1 {\gerU^{\scriptscriptstyle#1} \big( \hhat \big)}
\def \gerUunoMh#1 {\gerU_1^{\scriptscriptstyle #1} \big( \hhat \big)}
\def \gerUepsilonMh#1
   {\gerU_\varepsilon^{\scriptscriptstyle #1} \big( \hhat \big)}

\def \uMbs #1#2 {U^{\scriptscriptstyle#1}_{\scriptscriptstyle #2}}
\def \gerUMbs #1#2
   {\gerU^{\scriptscriptstyle#1}_{\scriptscriptstyle #2}}
\def \calUMbs #1#2
   {{\Cal U}^{\scriptscriptstyle#1}_{\scriptscriptstyle #2}}

\def \uMbm#1 {U^{\scriptscriptstyle#1}_{\scriptscriptstyle \leq}}
\def \calUMbm#1
   {{\Cal U}^{\scriptscriptstyle#1}_{\scriptscriptstyle \leq}}
\def \gerUMbm#1
   {\gerU^{\scriptscriptstyle#1}_{\scriptscriptstyle \leq}}

\def \uMbp#1 {U^{\scriptscriptstyle#1}_{\scriptscriptstyle \geq}}
\def \calUMbp#1
   {{\Cal U}^{\scriptscriptstyle#1}_{\scriptscriptstyle \geq}}
\def \gerUMbp#1
   {\gerU^{\scriptscriptstyle#1}_{\scriptscriptstyle \geq}}

\def \gerFrg {{\frak F}{\frak r}_\ghat}
\def \gerFrh {{\frak F}{\frak r}_\hhat}
\def \calFrg {{{\calF}r}_\ghat}
\def \calFrh {{{\calF}r}_\hhat}


\document

\topmatter

{\ }

\vskip-33pt

\hfill   {\sl Mathematische Zeitschrift  {\bf 234}  (2000), 9--52  
\hskip11pt   ---   \hskip11pt   DOI: 10.1007/s002090050502 }

\vskip41pt

\title
   Dual affine quantum groups
\endtitle

\author
   Fabio Gavarini
\endauthor

\affil
   Universit\`a degli Studi di Roma ``Tor Vergata'' ---
Dipartimento di Matematica  \\
   \hbox{ Via della Ricerca Scientifica 1,
I-00133 Roma --- ITALY }  \\
\endaffil

\address\hskip-\parindent
   Universit\`a degli Studi di Roma ``Tor Vergata'' \newline
   Dipartimento di Matematica  \newline
   Via della Ricerca Scientifica 1  \newline
   I-00133 Roma --- ITALY  \newline
   e-mail: \ gavarini\@mat.uniroma2.it  \newline
\endaddress

 \abstract
   Let  $ \ghat $  be an untwisted affine Kac-Moody algebra, with its
Sklyanin-Drinfel'd structure of Lie bialgebra, and let  $ \hhat $  be the dual
Lie bialgebra.  By dualizing the quantum double construction   --- via
formal Hopf algebras ---   we construct a new quantum group  $ \uqh $,  dual of
$ \uqg $.  Studying its specializations at roots of 1 (in particular, its
classical limits), we prove that it yields quantizations of  $ \hhat $  and
$ \widehat{G}^\infty $  (the formal group attached to  $ \ghat $),  and we
construct new quantum Frobenius morphisms. The whole picture extends to the
untwisted affine case the results known for quantum groups of finite type.
 \endabstract

\endtopmatter


\footnote""{ AMS 1991 {\it Mathematics Subject Classification,}
Primary 17B37, Secondary 81R50 }

\footnote""{ Partially supported by a post-doc fellowship of the  {\it
Consiglio Nazionale delle Ricerche\/}  \, (Italy) }

\vskip0,7truecm

 \centerline{ \bf  Introduction }

\vskip13pt

\hfill  \hbox{\vbox{ \hbox{\it  \   "Dualitas dualitatum }
                     \hbox{\it \ \ \;\, et omnia dualitas" }
                     \vskip4pt
                     \hbox{\sl    N.~Barbecue, "Scholia" } } \hskip1truecm }

\vskip11pt

   Let  $ \ghat $  be an untwisted affine complex Kac-Moody algebra, with the
Sklyanin-Drin-fel'd structure of Lie bialgebra; let  $ \hhat $  be its dual Lie
bialgebra.  Let  $ R $  be the subring of complex rational functions having no
poles at roots of 1.  Let  $ \uqg $  be the quantum group   --- over the field
$ \Cq $  ---   associated to  $ \ghat $:  then there exists an integer form  $ \gerU(\ghat) $  of
$ \uqg $  over  $ R $  which for  $ \, q \rightarrow 1 \, $  specializes to
$ \ug $  as a Poisson Hopf coalgebra (cf.~[Lu2]).  On the other hand, another integer form  $ \calU(\ghat) $
exists which for  $ \, q \rightarrow 1 \, $  specializes (as a Poisson Hopf
algebra) to  $ F \big[ \widehat{H} \big] $,  the function algebra of an
infinite dimensional proalgebraic Poisson group
$ \widehat{H} $  whose tangent
Lie bialgebra is  $ \hhat $  (cf.~[BK]).
   All this can be seen as an application of (a "global version" of) the
{\it quantum duality principle\/}:  this claims (cf.~[Dr], \S 7, or
[CP], \S 6; see also [Ga3] for a proof) that the quantization of
a Lie bialgebra   --- via a quantum universal enveloping algebra
(QUEA) ---   provides also a quantization of the dual Lie
bialgebra (through its associated formal Poisson group)   --- via a
quantum formal series Hopf algebra (QFSHA) ---   and, conversely, a QFSHA
which quantizes a Lie bialgebra (via its associated formal Poisson group)
yields a QUEA for the dual Lie bialgebra as well.
%
%
   In addition, both  $ \gerU(\ghat) $  and  $ \calU
(\ghat) $  can be specialized at roots of 1, and special
{\it quantum Frobenius morphisms}  $ \, \gerU_\varepsilon(\ghat) \longtwoheadrightarrow \gerU_1(\ghat) \, $  and  $ \, \calU_1(\ghat) \longhookrightarrow \calU_\varepsilon(\ghat) \, $  exist which are quantum analogues (in characteristic zero!) of the
Frobenius morphisms  $ \, U(\ghat_{\Z_p}) \longtwoheadrightarrow
U(\ghat_{\Z_p}) \, $  and  $ \, F \big[ \widehat{H}_{\Z_p} \big] \longhookrightarrow F \big[ \widehat{H}_{\Z_p} \big] \, $  which exist in characteristic  $ p $.  Such results are not predicted by the quantum duality principle, and are typical of the Jimbo-Lusztig's approach to quantum groups.
                                                   \par
   Our aim is to find an analogue of  $ \uqg $  for the algebra  $ \hhat $
instead of  $ \ghat $;  inspired by the quantum duality principle, and encouraged by the finite-type case (cf.~[Ga1]), we choose as a reasonable
candidate the linear dual  $ {\uqg}^* $,  which has a natural structure of
formal Hopf algebra.  This dual can be studied by dualizing
Drinfel'd's construction of the quantum double and using Tanisaki's pairings between quantum Borel (sub)algebras.
So we find a description of  $ {\uqg}^* $,  as a topological algebra with
formal Hopf algebra structure, in terms of generators and relations: we
call this
algebra  $ \uqh $,  for in fact we prove that it is for  $ \hhat $  what
$ \uqg $  is for  $ \ghat $.  In particular,  $ \uqh $  has an integer form  $ \gerU \big( \hhat \big) $
(over  $ R $)  which is a quantization of  $ \uh $;  moreover,  $ \uqh $  has
also a second integer form  $ \calU \big( \hhat \big) $  which is a quantization of  $ \, F^\infty \big[
\widehat{G} \big] \, $,  where  $ \widehat{G} $  of course is a Kac-Moody Poisson group with  $ \ghat $  as tangent Lie bialgebra.  More in general, both  $ \gerU \big( \hhat
\big) $  and  $ \calU \big( \hhat \big) $  can be specialized at roots
of 1, and quantum Frobenius morphisms exist (for both kind of forms), which are
dual of those of  $ \uqg $  and have a similar description.
                                                      \par
   Finally, a brief sketch of the main ideas of the paper.  First, since  $ \uqg $  is a quotient of a quantum double
$ \, D_q(\ghat) := D \left( U_q \big( \hat{\gerb}_- \big), U_q \big( \hat{\gerb}_+ \big), \pi \right) \, $,  its
linear dual  $ {\uqg}^* $  embeds into  $ {D_q(\ghat)}^* $.  Second, since  $ \, D_q(\ghat) \cong U_q \big( \hat{\gerb}_+ \big) \otimes U_q \big( \hat{\gerb}_-
\big)  \, $  (as coalgebras) we have
$ \, {D_q(\ghat)}^* \cong {U_q \big(
\hat{\gerb}_+ \big) }^* \otimeshat
{U_q \big( \hat{\gerb}_- \big) }^* \, $
(as algebras), where  $ \, \otimeshat \, $
denotes topological tensor product.  Third, since quantum Borel algebras of opposite sign are perfectly paired, their
linear duals are suitable completions of quantum Borel algebras of opposite sign: thus we
find a presentation of  $ {\uqg}^* $  (as a topological algebra) by generators and relations which leads us
to  {\sl define}  $ \, \uqh := {\uqg}^* \, $
(actually, one has to keep track
of some choice of lattices too, involved in the toral parts).  From this, all
claimed results follow.  In particular, the form  $ \calU \big( \hhat \big) $  is the subset (of  $ \, \uqh := {\uqg}^* \, $)  of linear functions on  $ \uqg $  which are $ R $--valued  on  $ \gerU(\ghat) $,  so  $ \, \calU \big( \hhat \big) \cong {Hom}_R \big( \gerU(\ghat), R \big) \, $  whence all results about specialisations of  $ \calU \big( \hhat \big) $  and its quantum Frobenius morphisms follow from those about  $ \gerU(\ghat) $.  On the other hand, the form  $ \gerU \big( \hhat \big) $  is a  {\sl proper}  subset of  $ \, {Hom}_R \big( \calU(\ghat), R \big) \, $,  for sort of a (non-trivial) "locality condition" is required for elements of  $ \, {Hom}_R \big( \calU(\ghat), R \big) \, $  to belong to  $ \gerU  \big( \hhat \big) $.

\vskip7pt
%
%

\centerline { ACKNOWLEDGEMENTS }

\vskip4pt

  The author wishes to thank M.~Rosso and E.~Strickland for
many helpful discussions, and I.~Damiani and J.~Beck for
explanations about their papers.  This article was partly
prepared during a one-year stay of the author at the Institut
de Recherche Math\'ematique Avanc\'ee of Strasbourg (France):
the author takes this opportunity to thank all the colleagues
and the staff of I.R.M.A.~for the warm hospitality.

\vskip1,7truecm

\centerline{ \bf  \S \; 1 \,  The classical objects }

\vskip10pt

   {\bf 1.1  Cartan data.}  \  Let  $ \gerg $  be a simple
finite dimensional complex Lie algebra, and consider the
folllowing data.
                                                          \par
   We take  $ I_0 = \{1, \dots, n \} \, $  to be the set of vertices of the
Dynkin diagram of  $ \gerg $  (see [Bo] for the identification between  $ I_0 $
and  $ \{1, \dots, n \} $);  $ \, A_0 = {\big( a_{ij} \big)}_{i,j \in I_0} \, $
the Cartan matrix of  $ \gerg $;  $ \, D_0 = \text{diag}(d_1, \dots, d_n) \, $
the (unique) diagonal matrix with relatively prime positive integral entries
such that  $ \, D_0 A_0 \, $  is symmetric;  $ \gert $,  a Cartan subalgebra of
$ \gerg $,  with a fixed basis  $ \, \{h_1, \dots, h_n\} = \{\, h_i \mid i \in
I_0 \,\} \, $;  $ \Phi_0 = \Phi_{0,+} \cup \left( -\Phi_{0,+} \right) \subseteq
\gert^* \, $  the root system of  $ \gerg $,  with  $ \, \Phi_{0,+} \, $  the
set of positive roots, and  $ \, \Pi := \{\alpha_1, \dots, \alpha_n\} = \{\,
\alpha_i \mid i \in I_0 \,\} \, $  the set of simple roots;  $ Q_0 :=
\sum_{\alpha \in \Phi_0} \Z \alpha = \oplus_{i \in I_0} \Z \alpha_i \, $  the
root lattice of  $ \gerg $,  and  $ Q_0^\vee := \sum_{\alpha \in \Phi_0} \Z
\alpha^\vee = \oplus_{i \in I_0} \Z \alpha_i^\vee \, $  the coroot lattice;
$ W_0 $  the Weyl group of  $ \gerg $.  Finally, we fix a function  $ \, o : I_0
\longrightarrow \{ \pm 1 \} \, $  such that  $ \, a_{ij} < 0 \Longrightarrow
o(i) \, o(j) = -1 \, $.
                                                          \par
  We denote  $ \ghat $  the untwisted affine Kac-Moody algebra
associated to  $ \gerg $  and we consider its loop-algebra like
realization as  $ \; \ghat = \gerg \otimes_\C \C \left[t,\tm\right]
\oplus \C \cdot c \oplus \C \cdot \partial \; $  with the Lie bracket
given by: \  $ \, [c,z] = 0 \, $,  $ \; \left[ \partial, x \otimes t^m
\right] = m x \otimes t^m \, $,  $ \; \left[ x \otimes t^r, y \otimes
t^s \right] = [x,y] \otimes t^{r+s} + \delta_{r,-s} r \, (x,y) \,
c \, $  for all  $ \, z \in \ghat $,  $ \, x, y \in \gerg \, $,
$ \, m, r, s \in \Z \, $  where  $ \, (\,\cdot\, , \,\cdot\, ) \, $
is the Killing form of  $ \gerg $,  normalized in such a way that
$ \, \big( h_i, h_j \big) = {\, a_{ij} \, \over \, d_j \,} \, $.
                                                          \par
  For  $ \ghat $  we define:  $ \, I := \{0,1,\dots,n\} \supset I_0 \, $  to be
the set of vertices of the Dynkin diagram, and  $ \, I_\infty := I \, \cup
\{\infty\} \, $;  $ \, A = {\big( a_{ij} \big)}_{i,j \in I} \, $  the
(generalized) Cartan matrix and  $ \, D = \text{diag}(d_0, d_1, \dots, d_n)
\, $  with  $ \, d_0 = 1 \, $  (so that  $ \, D A \, $  is symmetric);  $ \that
:= \gert \oplus \C \cdot c \oplus \C \cdot \partial \, (\, \subseteq \ghat \,)
\, $;  $ \, \Phi = \Phi_+ \cup (-\Phi_+) \, \left(\, \subset {(\gert \oplus \C
\cdot c)}^* \subset \that^* \, \right) \, $  the root system,  $ \, \Phi_+ =
\phipre \cup \phipim \, $  the set of positive roots,  $ \, \{ \alpha_0,
\alpha_1, \dots, \alpha_n \} = \{\, \alpha_i \mid i \in I \,\} \, $  the set of
simple roots,  $ \, \phipim = \, \{\, m \delta \mid m \in \N_+ \, \} \, $  the
set of imaginary positive roots (where  $ \, \delta = \sum_{i \in I} d_i
\alpha_i = \theta + \alpha_0 \, $  and  $ \theta $  is the longest positive root
of  $ \gerg $),  $ \, \phipre = \Phi_{0,+} \cup \{\, \alpha + m \delta \mid
\alpha \in \Phi_0 \, ,  m > 0 \,\} \, $  the set of real positive roots.  Then
$ \ghat $  has a decomposition into direct sum of  $ \that $  and root spaces
$ \, \ghat = \that \oplus \big( \oplus_{\alpha \in \Phi} \ghat_\alpha \big)
\, $,  and  $ \, \text{dim}_{\Bbb C} \big( \ghat_\alpha \big) = 1 \, $  if  $ \,
\alpha \in \phire \, $,  $ \, \text{dim}_{\Bbb C} \big( \ghat_\alpha \big) = \#
\big(I_0\big) = n \, $  if  $ \, \alpha \in \phipim \, $;  therefore we define
the set  $ \phitildep $  of "positive roots with multiplicity" as  $ \,
\phitildep := \phipre \cup \phitildepim \, $,  where  $ \, \phitildepim :=
\phipim \times I_0 \, $;  then we denote  $ \, p : \phitildep \rightarrow
\Phi_+ \, $  the natural projection map.
                                             \par
   Furthermore, we have:  the root lattice (of  $ \ghat $)  $ \, Q =
\sum_{\alpha \in \Phi} \Z \cdot \alpha = \oplus_{i \in I} \Z \cdot \alpha_i = \Z
\cdot \alpha_0 \oplus Q_0 = Q_0 \oplus \Z \cdot \delta \, $,  $ \, Q_\infty := Q
\oplus \Z \cdot \alpha_\infty \, $  (where  $ \, \alpha_\infty \in \that^* \, $  is defined by  $ \, \langle \alpha_\infty, \gert \rangle = 0 $,  $ \, \langle \alpha_\infty, \partial \rangle = 0 $,  $ \, \langle \alpha_\infty, c \rangle = 1 $ ),  and the order
relation  $ \, \leq \, $  on  $ Q_\infty $  given by  $ \, \alpha \leq \beta
\iff \beta - \alpha \in Q_+ \, $,  with  $ \, Q_+ := \sum_{i \in I} \N \cdot
\alpha_i \, $;  the extended Cartan matrix  $ \, A_\infty = {\big( a_{ij}
\big)}_{i,j \in I_\infty} \, $,  defined by setting  $ \, a_{i \infty}
:= \delta_{i,0} $,  $ \, a_{\infty j}:= \delta_{0,j} $;  the diagonal
matrix  $ \, D_\infty = \text{diag}(d_0, d_1, \dots, d_n, d_\infty)
\, $  with  $ \, d_\infty = 1 \, $;  the non-degenerate symmetric
bilinear form on  $ \, Q_\infty \otimes_\Z {\Bbb R} \, $
given by  $ \, (\alpha_i, \alpha_j) = d_i a_{ij} \,
(\, \forall\, i, j \in I_\infty \,) $;  the group  $ \, W = W_0 \ltimes Q_0^\vee
\, $,  the subset of simple reflections  $ \, \{ s_0, s_1, \dots, s_n \} = \{\,
s_i \mid i \in I \,\} \, (\, \subseteq W \,) \, $,  and the length function
$ \, l \colon \, W \rightarrow \N \, $;  the braid group  $ {\Cal B} $
(associated to  $ W \, $),  generated by  $ \, \{ T_0, T_1, \dots, T_n \} = \{\,
T_i \mid i \in I \,\} \, $,  and the section  $ \, T \colon \, W \rightarrow
{\Cal B} \, $  such that  $ \, T_w = T_{i_1} \cdots T_{i_r} \, $  for
all  $ \, w = s_{i_1} \cdots s_{i_r} \in W \, $  with  $ l(w) = r $.
Notice that the form  $ (\, \cdot \, , \, \cdot \,) $  is
$ W $--invariant.  Finally, set  $ \, Q_- :=  -Q_+ =
\sum_{i \in I} (-\N) \cdot \alpha_i \, $.
                                             \par
   We define the  {\it weight lattice\/}  $ \, P_\infty := {Hom}_\Z \big(
Q_\infty, \Z \big) $  to be the dual lattice of  $ Q_\infty $,  and we fix the
$ \Z $--basis  $ \{\, \omega_i \mid i \in I_\infty \,\} $  such that  $ \,
\langle \alpha_i \vert \omega_j \rangle = \delta_{ij} \, $,  where  $ \, \langle
\ \vert \ \rangle \colon Q_\infty \times P_\infty \rightarrow \Z \, $  is the
natural pairing; we set  $ \, P_\infty^+ := \sum^n_{i=1} \N \omega_i \, $  (the
subset of  {\it dominant integral weights\/}).  Remark that  $ W $  acts on
$ Q_\infty $  too.  Via the form  $ \, (\, \cdot \, , \, \cdot \,) $  we can
embed  $ Q_\infty $  into  $ P_\infty $,  so that  $ \, \alpha_i = \sum_{j \in
I_\infty} d_i a_{ij} \omega_j \, $  for all  $ \, i \in I_\infty \, $.  We
extend the form  $ \, (\, \cdot \, , \, \cdot \,) \colon Q_\infty \times
Q_\infty \rightarrow \Z \, $  to a
(non-degenerate symmetric) pairing  $ \, (\, \cdot \, , \, \cdot \,) \colon \Q
Q_\infty \times \Q Q_\infty \rightarrow \Q \, $  of  $ \Q $--vector  spaces by
scalar extension (hereafter  $ \, \Q T_\infty := \Q \otimes_\Z T_\infty (T = Q, P)
\, $):  then restriction gives a pairing  $ \, (\, \cdot \,
, \,\cdot \,) \colon P_\infty \times P_\infty \rightarrow \Q \, $
(looking at  $ P_\infty $  as a sublattice of  $ \, \Q P_\infty =
\Q Q_\infty \, $),  which extends  $ \, (\, \cdot \, , \,\cdot \,)
\colon Q_\infty \times P_\infty \rightarrow \Z \, $  and takes
values in $ \Z \left[ \Delta_\infty^{\,-1} \right] $,  where
$ \, \Delta_\infty := det \left( {(\aij)}_{i,j \in I_\infty}
\right) \, $.  Finally we define  $ \, d_\alpha := {(\alpha,\alpha)
\over 2} \, $  for all  $ \, \alpha \in \phipre \, $.
                                                           \par
  Given any pair of lattices  $ (M,M') $,  with  $ \, Q_\infty \leq M, M' \leq
\Q P_\infty \, $,  we say that they are  {\it dual of each  other\/}  if  $ \,
M' = \big\{\, y \in \Q P_\infty \,\big\vert\, \langle M, y \rangle \subseteq \Z
\,\big\} \, $,  $ \, M = \big\{\, x \in \Q P_\infty \,\big\vert\, \langle x, M'
\rangle \subseteq \Z \,\big\} \, $,  the two conditions being equivalent; then
for any lattice  $ M $  with  $ \, Q_\infty \leq M \leq \Q P_\infty \, $  there
exists a unique dual lattice  $ M' $  such that  $ \, Q_\infty \leq M' \leq \Q
P_\infty \, $  and  $ \, (\, \cdot \, , \, \cdot \,) \colon \Q P_\infty \times
\Q P_\infty \rightarrow \Q \, $  restricts to a perfect pairing  $ \, (\, \cdot
\, , \, \cdot \,) \colon M \times M' \rightarrow \Z \, $;  in particular  $ \,
P_\infty' = Q_\infty \, $  and  $ \, Q_\infty' = P_\infty \, $.  In the sequel
we denote by  $ \, \{\, \mu_i \mid i \in I_\infty \,\} \, $  and  $ \, \{\, \nu_i \mid i \in I_\infty \,\} \, $  fixed  $ \Z $--bases  of  $ M $  and  $ M'
$  dual of each other, i.~e.~such that  $ \, (\mu_i \vert \nu_j) = \delta_{i j}
\, $  for all  $ \, i, j \in I_\infty \, $,  and we set  $ \, M_+ := M \cap
P_\infty^+ \, $.  In the following our constructions will work in general for the pairs of
dual lattices  $ \big( P_\infty, Q_\infty \big) $  and  $ \big( Q_\infty,
P_\infty ) $;  but in the simply laced case (in which  $ \, \langle \ ,\ \rangle
= (\ ,\ ) \, $)  $ \left( M, M' \right) $  will be  {\sl any}  pair of dual
lattices.

\vskip7pt

   {\bf 1.2  The classical Manin triple.}  \  Let be given for  $ \ghat $  the
usual presentation by Chevalley-type generators  $ f_i $,  $ h_j $,  $ e_i $
($ i \in I $,  $ j \in I_\infty $)  and relations; then let  $ \hat\gern_+ $,
resp.~$ \hat\gern_- $,  be the Lie subalgebra of  $ \ghat $  generated by the
$ e_i $'s,  resp.~the  $ f_i $'s  ($ i \in I \, $).  Now let  $ \hat\gerk $  be
the Lie algebra  $ \, \hat\gerk := \ghat \oplus \ghat \, $;  inside it we find a
diagonal copy of  $ \ghat $  and a second Lie algebra
  $$  \hhat := \Big\{\, \big( n_- t_-, t_+ n_+ \big) \;\Big\vert\; n_\pm \in
\hat\gern_\pm \, , \, t_\pm \in \that , t_- + t_+ = 0 \,\Big\}  \quad  \big(\,
\leq \hat\gerb_- \times \hat\gerb_+ \leq \hat\gerk \big) \, ;  $$
in particular  $ \, \hhat = \hat\gern_- \oplus \that \oplus \hat\gern_+ \, $  as
vector spaces.  Define a bilinear form on  $ \hat\gerk $  by the formula  $ \;
\big\langle x_1 \oplus y_1 , x_2 \oplus y_2 \big\rangle := {\,1\, \over
\,2\,} \, (y_1,y_2) - {\,1\, \over \,2\,} \, (x_1,x_2) \, $,  \;  where  $ \,
(\, \cdot \, , \, \cdot \,) \, $  is the Killing form; this form makes  $ \left(
\hat\gerk, \ghat, \hhat \right) $  into a Manin triple; in particular
$ \ghat $  and  $ \hhat $  are Lie bialgebras, and the bilinear form on
$ \hat\gerk $  gives by restriction a non-degenerate pairing  $ \, \langle \ ,
\ \rangle \colon \hhat \otimes \ghat \rightarrow \C \, $  which respect the Lie
bialgebra structure on both sides, that is  $ \; \big\langle x, [y_1,y_2]
\big\rangle = \big\langle \delta_{\hhat} (x) \, , \, y_1 \otimes y_2 \big\rangle
\, $,  $ \, \big\langle [x_1,x_2] \, , \, y \big\rangle = \big\langle x_1
\otimes x_2, \delta_{\ghat} (y) \big\rangle \, $,  \; where  $ \delta $  is the
Lie cobracket.

\vskip7pt

   {\bf 1.3  The Poisson Hopf coalgebra  $ \uh $.}  \  The presentation of
$ \ghat $  by generators and relations gives a similar one for  $ \ug $,  with
the same generators.  From the latter we take for  $ \uh $  the following
presentation, where  $ \, \text{f}_i = f_i \oplus 0 \, $,  $ \, \text{h}_i = h_i
\oplus h_i \, $,  $ \, \text{e}_i = 0 \oplus e_i \, $:  $ \uh $  is the
associative  \hbox{$ \C $--alge}bra  with 1 generated by  $ \text{f}_i $,
$ \text{h}_i $,  $ \text{e}_i $  ($ i = \in I_\infty $)  with relations
  $$  \eqalign {
       \text{h}_r \text{h}_s - \text{h}_s \text{h}_r = 0 \, ,  \quad  \text{e}_i
\text{f}_j - \text{f}_j \text{e}_i = 0  &  \, ,  \quad  \text{h}_r \text{f}_j -
\text{f}_j \text{h}_r = a_{rj} \, \text{f}_j \, ,  \quad  \text{h}_r \text{e}_j
- \text{e}_j \text{h}_r = a_{rj} \, \text{e}_j  \cr
       \sum_{k=0}^{1-\aij} \! {(-1)}^k {1 - \aij \choose k}
\text{f}_i^{\, 1-\aij-k} \text{f}_j \text{f}_i^{\, k}  &  = 0 \, ,
\;  \sum_{k=0}^{1-\aij} \! {(-1)}^k {1 - \aij \choose k}
\text{e}_i^{\, 1-\aij-k} \text{e}_j \text{e}_i^{\, k} = 0  \cr }
   \hfill \hskip3pt  (1.1)  $$
for  $ \, r, s \in I_\infty \, $,  $ \, i, j \in I \, $  ($ \, i \neq j \, $  in
the bottom row); its Hopf structure is given by
  $$  \Delta (x) = x \otimes 1 + 1 \otimes x \, ,  \;\,
S(x) = - x \, ,  \;\,  \epsilon(x) = 0 \phantom{\; .}  \quad
\forall \;\; x \in \big\{ \text{f}_i, \text{h}_j, \text{e}_i
\big\vert i \in I, j \in I_\infty \big\}   \hfill \hskip9pt  (1.2)  $$
(the natural one) and the co-Poisson structure  $ \delta \! = \!
\delta_{\hhat} \colon \uh \! \rightarrow \! \uh \! \otimes \! \uh $  by
  $$  \eqalign {
      \delta(\text{f}_i)  &  = d_i \cdot \big( \text{h}_i \otimes \text{f}_i -
\text{f}_i \otimes \text{h}_i \big) + 2 \, d_i^{-1} \cdot \!
         \!\!\!\!\!\! \sum_{\Sb  \alpha, \beta \in \phitildep  \\
                                 p(\alpha) - p(\beta) = -\alpha_i  \\  \endSb}
\!\!\!\!\! c^{i,+}_{\alpha,\beta} \, d_\alpha d_\beta \cdot \big(
\text{e}_\alpha \otimes \text{f}_\beta -  \text{f}_\beta \otimes
\text{e}_\alpha
\big)  \cr
      \hskip9pt  \delta(\text{h}_i)  &  = 4 \, d_i^{-1} \cdot
\!\!\!\!\! \sum_{\Sb  \alpha, \beta \in \phitildep  \\
                     p(\alpha) - p(\beta) = 0  \\   \endSb}
\!\!\!\!\! c_{\alpha,\beta} \, d_\alpha \, \big( \alpha_i \big\vert p(\alpha)
\big) \cdot \big( \text{e}_\alpha \otimes \text{f}_\beta - \text{f}_\beta
\otimes \text{e}_\alpha \big)   \hskip82pt (1.3)  \cr
      \delta(\text{e}_i)  &  = d_i \cdot \big( \text{e}_i \otimes \text{h}_i -
\text{h}_i \otimes \text{e}_i \big) + 2 \, d_i^{-1} \cdot \!
         \!\!\!\!\!\! \sum_{\Sb  \alpha, \beta \in \phitildep  \\
                                 p(\alpha) - p(\beta) = +\alpha_i  \\  \endSb}
\!\!\!\!\! c^{i,-}_{\alpha,\beta} \, d_\alpha d_\beta \cdot \big(
\text{e}_\alpha \otimes \text{f}_\beta - \text{f}_\beta \otimes \text{e}_\alpha
\big)  \cr }  $$
for all  $ \, i \in I $,  $ j \in I_\infty \, $;  here the
$ \text{e}_\gamma $'s  and the  $ \text{f}_\gamma $'s  are suitable
"root-with-multiplicity vectors" (respectively of "weight"  $ \, p(\gamma) \, $
and  $ \, - p(\gamma) \, $)  such that  $ \, \big\langle \text{e}_\gamma, f_\eta
\big\rangle = + \delta_{\gamma,\eta} d_\gamma \big/ 2 \, $,  $ \, \big\langle
\text{f}_\gamma, f_\eta \big\rangle = - \delta_{\gamma,\eta} d_\gamma \big/ 2
\, $  ($ f_\eta $  and  $ e_\eta $  being root-with-multiplicity vectors of
$ \ghat $),  the  $ \, c_{\alpha,\beta} \, $'s  are given by the equations  $ \,
\big[ e_\beta, f_\alpha \big] = c_{\alpha,\beta} \, h_{p(\alpha)} \, $  and the
$ c^{i,\pm}_{\alpha,\beta} $'s  by  $ \, \big[ f_\alpha, e_\beta \big] =
c^{i,-}_{\alpha,\beta} \cdot f_i \, $,  $ \, \big[ f_\alpha, e_\beta \big] =
c^{i,+}_{\alpha,\beta} \cdot e_i \, $.  Note that the formulas above contains
infinite sums, so  $ \delta $  in fact takes values in a certain completion of
$ \, \uh \otimes \uh \, $;  hence to be precise  $ \uh $  is a Hopf algebra
which is co-Poisson only in a larger sense; and similarly for  $ \hhat $  as a
Lie bialgebra).

\vskip1,7truecm

\centerline{ \bf  \S \; 2 \,  Quantum Borel algebras and DRT
pairings }

\vskip10pt

   {\bf 2.1  Quantum Borel algebras.}  \  From now on  $ M $  will be any
lattice such that  $ \, Q_\infty \leq M \leq P_\infty \, $;  then  $ M' $  will
be the dual lattice, as in \S 1.1 (upon conditions therein).
                                                   \par
  For all  $ \, s, n \in \N \, $,  let  $ \, {(n)}_q := {q^n - 1 \over q - 1}
\; (\in \Z[q]) \, $,  $ \, {(n)}_q! := \prod_{r=1}^n {(r)}_q $,
       $ {({n \atop s})}_q := {{(n)}_q! \over {(s)}_q! {(n-s)}_q! } $\break
\noindent   $ (\in \Z[q]) \, $,  and  $ \, {[n]}_q := {q^n - q^{-n} \over q -
\qm} \; (\in \Z \left[ q, \qm \right]) \, $,  $ \, {[n]}_q! := \prod_{r=1}^n
{[r]}_q $,  $ {[{n \atop s}]}_q := {{[n]}_q!
\over {[s]}_q! {[n-s]}_q! } \; (\in \Z \left[ q, \qm \right]) \, $.  Let  $ \,
q_\alpha := q^{d_\alpha} \, $  for all  $ \, \alpha \in \phipre \, $,  and  $ \,
q_i := q_{\alpha_i} \, $  for all  $ \, i \in I \, $;  then set  $ \, q_{(r
\delta,i)} := q_i \, $  for every positive imaginary root  $ \, (r \delta,i) \in
\phipim \, $.
                                                   \par
  We define  $ \, U_q^{\scriptscriptstyle M} \! \big( \hat{\gerb}_- \big) \, $,
resp.~$ \, U_q^{\scriptscriptstyle M} \! \big( \hat{\gerb}_+ \big)  \, $,  to be
the associative  $ \Cq $--algebra  with 1 generated by  $ \, L_\mu $  ($ \mu \in
M \, $),  $ F_1, \dots, F_n $,  resp.~$ \, L_\mu $  ($ \mu \in M \, $),  $ E_1,
\dots, E_n $,  with  relations
  $$  \eqalign{
      L_0 = 1 \, , \quad  &  \qquad \quad  L_\mu L_\nu = L_{\mu+\nu} \, ,  \cr
      L_\mu F_j = q^{-( \alpha_j | \mu )} F_j L_\mu \, ,  &
\qquad  \sum_{p+s=1-\aij} {(-1)}^s {\left[ {1 - \aij \atop s} \right]}_{\! q_i}
F_i^p F_j F_i^s = 0  \cr
      \hbox{resp. \ \ }  L_\mu E_j = q^{( \alpha_j | \mu )}
E_j L_\mu \, ,  &  \qquad  \sum_{p+s=1-\aij} {(-1)}^s {\left[ {1 - \aij \atop
s} \right]}_{\! q_i} E_i^p E_j E_i^s = 0  \,  \cr }   \eqno (2.1)  $$
for all  $ \, i, j \in I $,  $ i \neq j $,  and  $ \, \mu, \nu \in M
\, $;  these are both Hopf algebras, with
  $$  \matrix
   \Delta (F_i) = F_i \otimes L_{-\alpha_i} + 1 \otimes F_i \, ,  &  \qquad
\epsilon (F_i) = 0 \, ,  &  \qquad  S (F_i) = - F_i L_{\alpha_i}  \\
   \Delta (L_\mu) = L_\mu \otimes L_\mu \, ,  &  \qquad  \epsilon (L_\mu) = 1
\, ,  &  \qquad  S (L_\mu) = L_{-\mu}  \\
   \Delta (E_i) = E_i \otimes 1 + L_{\alpha_i} \otimes E_i \, ,  &  \qquad
\epsilon (E_i) = 0 \, ,  &  \qquad  S (E_i) = - L_{-\alpha_i} E_i  \\
      \endmatrix  $$
for all  $ \, i \in I \, $,  $ \, \mu \in M \, $.  We also consider the
subalgebras  $ U_q^{\scriptscriptstyle M} \! (\gert) $  (generated by
the  $ L_\mu $'s),  $ U_q (\gern_-) $  (generated by the  $ F_i $'s),
$ U_q (\gern_+) $  (generated by the  $ E_i $'s).  In the sequel we
shall use the notation  $ \, K_\alpha := L_\alpha \, $,  $ \, M_\mu
:= L_\mu \, $,  $ \, \Lambda_\nu := L_\nu \, (\forall \, \alpha
\in Q_\infty, \mu \in M, \nu \in M' \,) \, $  (and in particular
$ \, K_i := K_{\alpha_i} \, $,  $ \, M_i := M_{\mu_i} \, $  $ \,
\Lambda_i := \Lambda_{\nu_i} \, $),  and  $ \, \uMbm{M} :=
U_q^{\scriptscriptstyle M} \! \big( \hat{\gerb}_- \big) \, $,
$ \, \uMbp{M} := U_q^{\scriptscriptstyle M} \! \big( \hat{\gerb}_+
\big)  \, $,  $ \, \uzM{M} := U_q^{\scriptscriptstyle M} \! (\gert)
\, $,  $ \, \um := U_q (\gern_-) \, $,  $ \, \up := U_q (\gern_+)
\, $.  Multiplication yields various linear isomorphisms, which we
shall refer to as  {\it triangular decompositions},  namely
  $$  \uMbm{M} \cong \um \otimes \uzM{M} \cong \uzM{M} \otimes \um \, ,  \qquad
\uMbm{M} \cong \up \otimes \uzM{M} \cong \uzM{M} \otimes \up  $$
                                                     \par
   A natural  $ Q_+ $--grading,  resp.~$ Q_- $--grading,  (of algebras)
is defined on  $ \uMbp{M} $,  resp.~$ \uMbp{M} $,  by setting
$ \, deg(E_i) := \alpha_i \, $,  $ \, deg(F_i) := -\alpha_i \, $
($ i \in I \, $);  then these are also gradings  {\sl of Hopf
algebras},  inherited by the various subalgebras defined above.

\vskip7pt

   {\bf 2.2  DRT pairings.}  \;  If  $ H $ is any Hopf algebra, we let
$ H^{op} $  be the same coalgebra with opposite multiplication, and
$ H_{op} $  the same algebra with opposite comultiplication.
                                       \par
   There exists perfect (i.~e.~non-degenerate) pairings of graded Hopf algebras
(cf.~[Ta])
  $$  \eqalign{
   \pi \colon \, {\left( {\uMbm M } \right)}_{op} \otimes {\uMbp {M'} }
\longrightarrow \Cq  \, ,  &  \quad \qquad  \pi \colon \, {\uMbm M } \otimes
{\left( \uMbp{M'} \right)}^{op} \longrightarrow \Cq  \cr
   \overline{\pi} \colon \, {\left( {\uMbp M } \right)}_{op} \otimes {\uMbm {M'}
} \longrightarrow \Cq  \, ,  &  \quad \qquad  \overline{\pi} \colon \, {\uMbp
M } \otimes {\left( {\uMbm {M'} } \right)}^{op} \longrightarrow \Cq  \cr }  $$
  $$  \eqalign{
   \pi (L_\mu,L_\nu) = q^{-(\mu,\nu)} \, ,  \;\ \pi (L_\mu,E_j) = 0
\, ,  &  \;\ \pi (F_i,L_\nu) = 0 \, ,  \;\  \pi (F_i,E_j) =
-\delta_{ij} {\left( q_i - q_i^{-1} \right)}^{-1}  \cr
   \overline{\pi} (L_\mu,L_\nu) = q^{+(\mu,\nu)} \, , \;\
\overline{\pi} (E_i,L_\nu) = 0 \, ,  &  \;\  \overline{\pi}
(L_\mu,F_j) = 0 \, ,  \;\ \overline{\pi} (E_i,F_j) = +\delta_{ij}
{\left( q_i - q_i^{-1}\right)}^{-1}  \cr }  $$
   \indent   These pairings were introduced by Drinfel'd, Rosso,
Tanisaki, and others, so we shall call them   {\it DRT pairings\/}.
If  $ \pi $  is any DRT pairing we shall also set  $ \, {\langle x,
y \rangle}_\pi \, $  for  $ \, \pi(x,y) \, $.

\vskip7pt

   {\bf 2.3  Quantum root vectors.}  \  We define quantum root vectors along the
lines of Beck's work (cf.~[Be1], [Be2]), but fixing conventions as in [Da],
[Ga2].  It is possible to define a total order  $ \, \preceq \, $  on the set
$ \phitildep $  such that  $ \; \beta_1 \preceq \beta_2 \preceq \beta_3 \preceq
\cdots \preceq \beta_{k-1} \preceq \beta_k \preceq \beta_{k+1} \preceq \cdots
\preceq \big( (r+1) \delta,n) \preceq (r\delta,1) \preceq (r\delta,2) \cdots
\preceq (2\delta,n) \preceq (\delta,1) \preceq (\delta,2) \preceq \cdots \preceq
(\delta,n) \preceq \cdots \preceq \beta_{-(k+1)} \preceq \beta_{-k} \preceq
\beta_{-(k-1)} \cdots \preceq \beta_{-2} \preceq \beta_{-1} \preceq \beta_0
\, $  \; and moreover  $ \, \big\{\, \beta_k \,\big\vert\, k \geq 1
\,\big\} = \big\{\, r \delta - \alpha \,\big\vert\, r > 0, \alpha \in
\Phi_{0,+} \,\big\} \, $  and  $ \, \big\{\, \beta_k \,\big\vert\, k \leq 0
\,\big\} = \big\{\, r \delta + \alpha \,\big\vert\, r \geq 0, \alpha \in
\Phi_{0,+} \,\big\} \, $.  Let such an order be fixed: then one defines   --- as
in [Ga2], \S 2.2 ---   (quantum root) vectors  $ \, E_\gamma \, $,  for all  $ \,
\gamma \in \phitildep \, $,  and also root vectors  $ F_\gamma $  (for all  $ \,
\gamma \in \phitildep \, $)  associated to negative roots.  Definitions give
$ \; E_\gamma \in {\left( U_+ \right)}_{p(\gamma)} \, $,  $ \; F_\gamma \in
{\left( U_- \right)}_{-p(\gamma)} \; $  for all  $ \, \gamma \in \phitildep
\, $.  For later use we recall an important property of imaginary root vectors:

\vskip5pt

   $ \underline{\hbox{\sl Claim}} $:  \, {\it  All root vectors
attached to imaginary roots commute with each other.}

\vskip7pt

   {\bf 2.4  PBW bases and orthogonality.}  \  It is proved in [Be2] that the
set  $ B_+ $  of ordered monomials in the root vectors  $ E_\alpha $'s
(according to the order  $ \, \preceq \, $  on  $ \phitildep \, $),  namely the
$ \, \prod_{\alpha \in \Phi_+} E_\alpha^{n_\alpha} $'s,  is a  $ \Cq $--basis
of  $ \up \, $  ({\sl Remark:\/}  {\it hereafter, when dealing with such a kind
of monomials, or the like, the function  $ \, \alpha \mapsto n_\alpha \in \N \, $
will always be zero
             for almost all\footnote{Hereafter by  {\it "almost all"}  we shall
always mean  {\it "all but a finite number of"}.}
$ \alpha \, $}).  Similarly, the set  $ B_- $  of ordered monomials in the root
vectors  $ F_\alpha $'s  is a  $ \Cq $--basis  of  $ \um \, $;  in addition, the
set  $ \, B_0^{\scriptscriptstyle M} := \{\, K_\alpha \mid \alpha \in Q_\infty
\,\} \equiv \left\{\, \prod_{i \in I_\infty} K_i^{l_i} \,\Big\vert\, l_i \in \Z
\,\; \forall \, i \in I_\infty \,\right\} \, $  is a  $ \Cq $--basis  of
$ \uzM{M} $.  Then from triangular decompositions one concludes that the sets of
ordered monomials  $ \, B_+ \cdot B_0^{\scriptscriptstyle M} \, $  and  $ \,
B_0^{\scriptscriptstyle M} \cdot B_+ \, $  are  $ \Cq $--bases  of  $ \up \, $,
and similarly  $ \, B_- \cdot B_0^{\scriptscriptstyle M} \, $  and  $ \,
B_0^{\scriptscriptstyle M} \cdot B_- \, $  are  $ \Cq $--bases  of  $ \um \, $.
From [Da] we know that
  $$  \displaylines{
   \pi \left( F_\alpha , E_\beta \right) = {\, \delta_{\alpha,\beta} \, \over
\, \left( q_\alpha^{-1} - q_\alpha \right) \,} \, ,  \;\;  \pi \left( F_\alpha
, E_\gamma \right) = 0 \, ,  \;\;  \pi \left( F_\gamma , E_\alpha \right) = 0
\qquad  \forall \, \alpha, \beta \in \phipre \, , \, \gamma \in \phitildepim
\cr
   \pi \left( F_{(r \delta, i)} , E_{(s \delta, j)} \right) = \delta_{r,s}
{\big( o(i) o(j) \big)}^r {\, {[r a_{ij}]}_{q_i} \, \over \, r \left( q_j^{-1}
- q_j \right) \,}  \qquad \forall \, (r, \delta, i), (s \delta, j) \in
\phitildepim  \cr }  $$
and similarly for  $ \overline{\pi} $,  for we have  $ \, \overline{\pi} \left(
E_\alpha , F_\beta \right) = - \pi \left( F_\beta , E_\alpha \right) \, $;
these formulas are the starting point to build up orthogonal bases of
$ \um $  and  $ \up \, $:  the key result is the following.

\vskip7pt

\proclaim{Lemma 2.5 (cf.~[Da], [Ga2])}  For all  $ r \in \N_+ \, $,  let
$ V_r $,  resp.~$ W_r $,  be the  $ \Cq $--vector  space with basis  $ \{\,
E_{(r \delta, i)} \,\vert\, i \in I_0 \,\} $, resp.~$ \{\, F_{(r \delta, i)}
\,\vert\, i \in I_0 \,\} $,  and let  $ \, \{\, x_{r,i} \mid i \in I_0 \,\}
\, $  and  $ \, \{\, y_{r,j} \mid j \in I_0 \,\} \, $  be bases of  $ V_r $
and  $ W_r $  orthogonal of each other with respect to  $ \pi $,  namely  $ \,
\pi \big( y_{r,i}, x_{r,i} \big) \not= 0 \, $  and  $ \,
\pi \big( y_{r,j}, x_{r,i} \big) = 0 \, $  for all  $ i, j \in I_0 \, $.  Then
$ \, \big\{\, \prod_{k \leq 0} E_{\beta_k}^{n_k} \cdot \! \prod_{r \in \N, i \in
I_0} x_{r,i}^{n_{r,i}} \cdot \prod_{k > 0} E_{\beta_k}^{n_k} \,\big\vert\, n_k,
n_{r,i} \in \N \; \forall \, k, i \,\big\} \, $  and  $ \, \big\{\, \prod_{k
\leq 0} F_{\beta_k}^{m_k} \cdot \! \prod_{r \in \N, j \in I_0} y_{s,j}^{m_{s,j}}
\cdot \prod_{k > 0} F_{\beta_k}^{m_k} \,\big\vert\, m_k, m_{s,j} \in \N \;
\forall \, k, j \,\big\} \, $  are bases of  $ \up $  and  $ \um $  which are
orthogonal of each other.  And similarly for the pairing  $ \overline{\pi} $.
More precisely, we have
  $$  \displaylines{
   {\ }  \pi \left( \prod_{k \leq 0} F_{\beta_k}^{m_k} \cdot \!
\prod_{r \in \N_+, j \in I_0} y_{r,j}^{m_{r,j}} \cdot \prod_{k > 0}
F_{\beta_k}^{m_k} \, , \, \prod_{h \leq 0} E_{\beta_h}^{n_h} \cdot \! \prod_{s
\in \N_+, i \in I_0} x_{s,i}^{n_{s,i}} \cdot \prod_{h > 0} E_{\beta_h}^{n_h}
\right) =   \hfill {\ }  \cr
  {\ } \hfill   = \prod_{\alpha \in \phipre} \delta_{n_\alpha, m_\alpha}
q_\alpha^{\left( n_\alpha \atop 2 \right)} {\, {[n_\alpha]}_{q_\alpha} ! \,
\over \, {\left( q_\alpha^{-1} - q_\alpha \right)}^{n_\alpha} \,} \, \cdot
\prod_{r \in \N, i \in I_0} \delta_{n_{r,i}, m_{r,i}} n_{r,i}! \, \pi {\left(
y_{r,i} \, , x_{r,i} \right)}^{n_{r,i}}  \, ,  \; \phantom{\square}   \cr
   {\ }  \overline{\pi} \left( \prod_{k \leq 0} E_{\beta_k}^{n_k} \cdot \!
\prod_{r \in \N_+, i \in I_0} x_{r,i}^{n_{r,i}} \cdot \prod_{k > 0}
E_{\beta_k}^{n_k} \, , \, \prod_{h \leq 0} F_{\beta_h}^{m_h} \cdot \! \prod_{s
\in \N_+, j \in I_0} y_{s,j}^{m_{s,j}} \cdot \prod_{h > 0} F_{\beta_h}^{m_h}
\right) =    \hfill {\ }  \cr
  {\ } \hfill   = \prod_{\alpha \in \phipre} \delta_{n_\alpha, m_\alpha}
q_\alpha^{\left( n_\alpha \atop 2 \right)} {\, {[n_\alpha]}_{q_\alpha}! \,
\over \, {\left( q_\alpha - q_\alpha^{-1} \right)}^{n_\alpha} \,} \, \cdot
\prod_{r \in \N, i \in I_0} \delta_{n_{r,i}, m_{r,i}} n_{r,i}! \, \pi {\left(
x_{r,i} \, , y_{r,i} \right)}^{n_{r,i}}  \, .  \; \square  \cr }  $$
\endproclaim

\vskip7pt

   {\bf 2.6 Integer forms.}   \  Let  $ R $  be the subring of
$ \C(q) $  of all rational functions having no poles at roots of
unity of odd order.  Let  $ \gerUMbm{M} $  be the  $ R $--subalgebra
of  $ \uMbm{M} $  generated by the elements  $ \, F_i^{(m)} := F_i^m
\Big/ {[m]}_{q_i}! \, $,  $ \, \left( M_j ; \, c \atop t \right) :=
\prod_{s=1}^t {{q_j^{c-s+1} M_j - 1} \over {q_j^s - 1}} \, $  (the
so-called  $ q $--{\it divided powers\/})  and  $ \, M_j^{-1} \, $
for all  $ \, m, c, t \in \N $,  $ i \in I $,  $ j \in I_\infty \, $.
It is known (cf.~[Lu2]) that  $ \gerUMbm{M} $  is a Hopf subalgebra
of  $ \uMbm{M} $.  It is proved in [Ga2] that  $ \gerUMbm{M} $  has
a PBW basis (as an  $ R $--module)  of increasing ordered monomials
$ \; \Big\{\, \prod_{j \in I_\infty} \!\! \left( M_j ; \, 0 \atop t_i
\right) M_j^{-Ent({t_j/2})} \cdot \prod_{\alpha \in \phitildep} \!\!
F_\alpha^{(n_\alpha)}
                         \, \,\Big\vert $\break
\noindent   $ \, n_\alpha \! \in \! \N, t_j \! \in \! \N \,\Big\} $
(recall that the  $ n_\alpha $'s  are almost all zero), where we
use notation
  $$  F_\gamma^{(n)} := {\, F_\gamma^n \, \over \, {[n]}_{q_\gamma}! \,}  \quad
\forall \; \gamma \in \phipre \, ,  \qquad  F_{(r \delta,i)}^{(n)} := {\,
{\left( r \big/ {[r]}_{q_i} \cdot F_{(r \delta,i)} \right)}^n \, \over \, n!
\,}  \quad  \forall \; (r \delta, i) \in \phitildepim \, ;  $$
a similar PBW basis of decreasing ordered monomials also exists; thus
$ \gerUMbm{M} $  is an  $ R $--form  of  $ \uMbm{M} $.  Similarly we
define the Hopf subalgebra  $ \gerUMbp{M} $  and find PBW bases for it.
                                               \par
  Let  $ \, \ebar_\alpha := \left( q_\alpha - q_\alpha^{-1} \right) E_\alpha
\, $  for all $ \, \alpha \in \phitildep \, $,  and let  $ \, \calUMbp{M} \, $
be the  $ R $--subalgebra  of  $ \uMbp{M} $  generated by all the
$ \ebar_\alpha $'s  and all the  $ M_i $'s;  then (cf.~[BK])  $ \calUMbp{M} $
is a Hopf subalgebra of  $ \uMbp{M} $,  with a PBW basis (as an
$ R $--module)  $ \; \left\{\, \prod_{j \in I_\infty} M_j^{t_j} \cdot
\prod_{\alpha \in \phitildep} \ebar_\alpha^{\, n_\alpha} \, \Big\vert \,
n_\alpha \! \in \! \N, \, t_j \! \in \! \N \,\right\} \; $  of increasing ordered monomials and
a similar PBW basis of decreasing ordered monomials; in particular
$ \calUMbp{M} $  is an  $ R $--form  of  $ \uMbp{M} \, $.  The same procedure
yields the definition of the Hopf subalgebra  $ \calUMbm{M} \, $  and provides
PBW bases for it.
                                               \par
   Similar constructions and results hold for the algebras $ \um $,
$ \uzM{M} $,  $ \up $,  providing integer forms  $ \gerUm $,  $ \calUp $,  and
so on, all of them endowed with PBW bases of ordered "monomials".  Finally, we
have triangular decompositions
  $$  \displaylines{
   \gerUMbm{M} \cong \gerUm \otimes \gerUzM{M} \cong
\gerUzM{M} \otimes \gerUm \, ,  \; \gerUMbp{M} \cong \gerUp
\otimes \gerUzM{M} \cong \gerUzM{M} \otimes \gerUp  \cr
   \calUMbm{M} \cong \calUm \otimes \calUzM{M} \cong
\calUzM{M} \otimes \calUm \, ,  \; \calUMbp{M} \cong \calUp
\otimes \calUzM{M} \cong \calUzM{M} \otimes \calUp
\cr }  $$

\vskip7pt

   {\bf 2.7  $ R $--duality among integer forms.}   \  For all  $ r \in \N_+
\, $,  let  $ V_r $  and  $ W_r $  be the  $ \Cq $--vector  spaces defined in
Lemma 2.5.  Consider the elements  $ \, \ehat_{(r \delta, i)} := \big( r \big/
{[r]}_{q_i} \big) \cdot E_{(r \delta,i)} \in V_r \, $  ($ \, i \in I_0 \, $);
then  $ \, \Big\{\, \ehat_{(r \delta, i)} \,\Big\vert\, i \in I_0 \,\Big\}
\, $ is a basis of  $ V_r $.  The basis of  $ W_r $  dual of  $ \,
\Big\{\, \ehat_{(r \delta, i)} \,\Big\vert\, i \in I_0 \,\Big\} \, $  with
respect to  $ \pi $,  is the subset  $ \, \Big\{\,\dot{\fbar}_{(r
\delta,j)} \,\Big\vert\, j \in I_0 \,\Big\} \subset W_r \, $  such that
$ \, \pi \Big( \dot{\fbar}_{(r \delta,j)}, \ehat_{(r \delta,i)} \Big) =
\delta_{j{}i} \, $,  for all  $ \, i, j \in I_0 \, $.  Similarly, we
define the basis  $ \, \Big\{\, \,\dot{\!\fhat}_{(r \delta,j)}
\,\Big\vert\, j \in I_0 \,\Big\} \, $  of  $ W_r $  dual of  $ \, \Big\{\,
\ebar_{(r \delta, i)} \,\Big\vert\, i \in I_0 \,\Big\} \, $  with respect
to $ \pi $.  Similar definitions hold with  $ \overline{\pi} $  instead of
$ \pi $,  and reversing the roles of  $ E $  and  $ F $.
                                              \par
   Definitions and formulas in \S 2.4 give the following (see [Ga2],
Proposition 4.6):

\vskip5pt

   $ \underline{\hbox{\sl Claim {\it (a)}}} $:  \, {\it for
all  $ \, r \in \N_+ \, $,  the basis  $ \, \Big\{\,\dot{\fbar}_{(r
\delta,j)} \,\Big\vert\, j \in I_0 \,\Big\} \subset W_r \, $  of
$ \, W_r \, $  dual of  $ \, \Big\{\, \ehat_{(r \delta, i)} \,\Big\vert\,
i \in I_0 \,\Big\} \, $  (with respect to  $ \pi $  or to
$ \overline{\pi} \, $)  lies in the  $ R $--linear  span of  $ \, \big\{\,
\fbar_{(r \delta,j)} \,\big\vert\, j \in I_0 \,\big\} \subset W_r \, $:
in particular, it lies in  $ \calUm \, $.  Similarly, the basis  $ \,
\Big\{\, \,\dot{\!\fhat}_{(r \delta,j)} \,\Big\vert\, j \in I_0 \,\Big\}
\subset W_r \, $  of  $ \, W_r \, $  dual of  $ \, \big\{\, \ebar_{(r
\delta, i)} \,\big\vert\, i \in I_0 \,\big\} \, $  (with respect to
$ \pi $  or to  $ \overline{\pi} \, $)  lies in the  $ R $--linear  span
of  $ \, \Big\{\, \fhat_{(r \delta,j)} \,\Big\vert\, j \in I_0 \,\Big\}
\subset W_r \, $:  in particular, it lies in  $ \gerUm \, $.  Similar
statements hold when reversing the roles of  $ E $  and  $ F $.}

\vskip5pt

  Then Lemma 2.5, and results in [DL], \S 3 (see also Lemma 4.4 later on) give:

\vskip5pt

   $ \underline{\hbox{\sl Claim {\it (b)}}} $:  \, {\it If in one integer form
(of a subalgebra of a quantum Borel algebra) we fix a PBW basis, then on the
other hand the form (of the opposite algebra) of opposite typographic "font"
does contain a PBW basis   --- made with "dotted" root vectors (which are equal
to the old ones for real roots) ---   which is dual to the initial one}.

\vskip5pt

  Therefore  {\sl integer forms of opposite "fonts" are  {\it  $ R $--dual\/}
of each other},  in the following sense: for every DRT pairing, if we take
$ \gerU $  on one side, then the form  $ \calU $  on the other side is equal to
the subset of all elements which paired with  $ \gerU $  give a value in  $ R \,
$;  and similarly exchanging  $ \gerU $  and  $ \calU $.  For instance
  $$  \displaylines{
    \gerUzM{M} = \Big\{\, y \in \uzM{M} \,\Big\vert\, \pi \Big( \calUzM{{M'}} ,
\, y \Big) \subseteq R \,\Big\} = \Big\{\, x \in \uzM{M} \,\Big\vert\,
\overline{\pi} \Big( x \, , \, \calUzM{{M'}} \Big) \subseteq R \,\Big\}  \cr
    \calUzM{M} = \Big\{\, y \in \uzM{M} \,\Big\vert\, \pi \Big( \gerUzM{{M'}} ,
\, y \Big) \subseteq R \,\Big\} = \Big\{\, x \in \uzM{M} \,\Big\vert\,
\overline{\pi} \Big( x \, , \, \gerUzM{{M'}} \Big) \subseteq R \,\Big\} \cr
    \gerUm = \Big\{\, x \in \um \,\Big\vert\, \pi \Big( x \, , \, \calUp \Big)
\subseteq R \,\Big\} = \Big\{\, y \in \um \,\Big\vert\, \overline{\pi} \Big(
\calUp \, , \, y \Big) \subseteq R \,\Big\}  \cr
    \calUMbp{M} = \Big\{\, x \in \uMbp{M} \,\Big\vert\, \overline{\pi} \Big( x
\, , \, \gerUMbm{{M'}} \Big) \subseteq R \,\Big\} = \Big\{\, y \in \uMbp{M}
\,\Big\vert\, \pi \Big( \gerUMbm{{M'}} , \, y \Big) \subseteq R \,\Big\}
\cr }  $$

\eject

\centerline{ \bf  \S \; 3 \,  The quantum group  $ \uqMg{M} $ }

\vskip10pt

   {\bf 3.1  Quantum double.}  \  Let  $ H_- $,  $ H_+ $
be two arbitrary Hopf algebras on a ground field (or ring)
$ F $,  and let  $ \, \pi \colon {\big( H_- \big)}_{op} \otimes H_+
\rightarrow F \, $  be any arbitrary Hopf pairing.  The  {\it quantum
double \/}  $ \, D = D \big( H_-, H_+, \pi \big) \, $  is the algebra
$ \, T \big( H_- \oplus H_+ \big) \big/ {\Cal R} \, $,  where  $ T \big( H_-
\oplus H_+ \big) $  is the tensor algebra on  $ \, H_- \oplus H_+ \, $  and
$ {\Cal R} $  is the ideal of relations
  $$  \eqalign{
    1_{H_-} = 1 = 1_{H_+}  \,  ,  &  \qquad  x \otimes y
= x y  \qquad \qquad \;  \hbox{for} \; \, x, y \in H_+  \; \hbox{ or }
\; x, y \in H_-  \cr
    \sum_{(x),(y)} \pi \left( y_{(2)}, x_{(2)} \right) \, x_{(1)} \otimes
y_{(1)} =  &  \sum_{(x),(y)} \pi \left( y_{(1)},x_{(1)} \right) \, y_{(2)}
\otimes x_{(2)}  \quad \;  \hbox{for}  \; \, x \in H_+ , \,  y \in H_- \, .
\cr }  $$
   \indent   Then (cf.~[CP], \S \ 4.2.D)  $ D $  has a canonical
structure of Hopf algebra such that  $ \, H_- $,  $ H_+ \, $
are Hopf subalgebras of it and multiplication yields isomorphisms
of coalgebras
  $$  H_+ \otimes H_- \longhookrightarrow D \otimes D {\buildrel {m} \over
\llongrightarrow} D \, ,  \qquad  H_- \otimes H_+ \longhookrightarrow D
\otimes D {\buildrel {m} \over \llongrightarrow} D \, .   \eqno (3.1)  $$
   \indent   Now take  $ \, D_q^{\scriptscriptstyle M} \! (\ghat) := D
\left( \uMbm{Q_\infty} , \uMbp{M} , \pi \right) \, $:  by definition,
$ D_q^{\scriptscriptstyle M} (\ghat) $  is generated by  $ \, K_\alpha $,
$ L_\mu $,  $ F_i $,  $ E_i $   --- identified with  $ 1 \otimes K_\alpha $,
$ L_\mu \otimes 1 $,  $ 1 \otimes F_i $,  $ E_i \otimes 1 $  when thinking at
$ \, D_q^{\scriptscriptstyle M} (\ghat) \cong \uMbp{M} \otimes \uMbm{Q_\infty}
\, $ ---   ($ \alpha \in Q_\infty $,  $ \mu \in M $,  $ i \in I \, $),  while
the relations defining  $ {\Cal R} $  reduce to
  $$  \displaylines{
   K_\alpha L_\mu = L_\mu K_\alpha \; ,  \qquad  K_\alpha E_j = q^{+(\alpha_j
\vert \alpha)} E_j K_\alpha \; ,  \qquad  L_\mu F_j = q^{-(\alpha_j \vert
\mu)} F_j L_\mu  \cr
   E_i F_j - F_j E_i = \delta_{ij} {{L_{\alpha_i} - K_{-\alpha_i}} \over {q_i -
q_i^{-1}}}  \cr }  $$
   \indent   Finally, PBW bases of quantum Borel algebras provide PBW bases
of  $ D_q^{\scriptscriptstyle M} \! (\ghat) $.  In the sequel we shall also use
the notation  $ \, D_{\scriptscriptstyle M} := D_q^{\scriptscriptstyle M} \!
(\ghat) \, $.

\vskip7pt

   {\bf 3.2 The quantum algebra  $ \, \uqMg{M} \, $.}  \  Let
$ \gerK^{\scriptscriptstyle M} $  be the ideal of  $ D_q^{\scriptscriptstyle M}
\! (\ghat) $  generated by the elements  $ \, L \otimes 1 - 1 \otimes L \, $,
$ \, L \in \uzM{M} \, $;  then  $ \gerK^{\scriptscriptstyle M} $  is in fact a
Hopf ideal, whence  $ \, D_q^{\scriptscriptstyle M} \! (\ghat) \big/
\gerK^{\scriptscriptstyle M} \, $  is a Hopf algebra.  The presentation above
yields one of  $ \, \uqMg{M} := D_q^{\scriptscriptstyle M} \! (\ghat) \Big/
\gerK^{\scriptscriptstyle M} \, $:  it is the associative  $ \Cq $--algebra
with 1 given by generators  $ \, F_i $,  $ L_\mu $,  $ E_i \, $  and relations
  $$  \displaylines {
   \hfill   L_0 = 1 \, ,  \quad  L_\mu L_\nu = L_{\mu + \nu} = L_\nu L_\mu \, ,
\quad  E_i F_h - F_h E_i = \, \delta_{ih} {\, L_{\alpha_i} -
L_{-\alpha_i} \, \over \, q_i - q_i^{-1} \,}   \hfill \phantom{(3.2)}  \cr
   \hfill   L_\mu F_i = q^{-(\alpha_i | \mu)} F_i L_\mu \, ,  \qquad
\sum_{k=0}^{1-a_{ij}} (-1)^k \left[ \! {1-a_{ij} \atop k} \!
\right]_{\!q_i} \! F_i^{1-\aij-k} F_j F_i^k = 0   \hfill (3.2)  \cr
   \hfill   L_\mu E_i = q^{+(\alpha_i | \mu)} E_i L_\mu \, ,  \qquad
\sum_{k = 0}^{1-a_{ij}} (-1)^k \left[ \! {1-a_{ij} \atop k} \!
\right]_{\!q_i} \! E_i^{1-\aij-k} E_j E_i^k = 0   \hfill \phantom{(3.2)}  \cr }  $$
(for all  $ \, \mu \in M $,  $ i, j, h \in I \, $ with  $ \, i \neq j \, $)
with the Hopf structure given by
  $$  \matrix
       \Delta (F_i) = F_i \otimes L_{-\alpha_i} + 1 \otimes F_i \, ,  \quad  &
\quad  \epsilon (F_i) = 0  \, ,  \quad  &  \quad  S (F_i) = - F_i L_{\alpha_i}
\\
       \Delta (L_\mu) = L_\mu \otimes L_\mu \, ,  \quad  &  \quad  \epsilon
(L_\mu) = 1 \, ,  \quad  &  \quad  S (L_\mu) = L_{-\mu}  \\
       \Delta (E_i) = E_i \otimes 1 + L_{\alpha_i} \otimes F_i \, ,  \quad  &
\quad  \epsilon (E_i) = 0 \, ,  \quad  &  \quad  S (F_i) = - L_{-\alpha_i}
E_i    \\
      \endmatrix   \eqno (3.3)  $$
   \indent   Finally, let  $ \, pr_{\scriptscriptstyle M}
\colon \, D_q^{\scriptscriptstyle M} \! (\ghat) \longtwoheadrightarrow
D_q^{\scriptscriptstyle M} \! (\ghat) \Big/ \gerK^{\scriptscriptstyle M} =:
\uqMg{M} \, $  be the canonical Hopf algebra epimorphism; we shall also use
notation  $ \, K_\alpha := L_\alpha $,  $ M_\mu := L_\mu $,  for all  $ \,
\alpha \in Q_\infty, \mu \in M \, $.

\vskip7pt

   {\bf 3.3 Integer forms of  $ \uqMg{M} $.} \  Let  $ \gerUMg{M} $  be
the  $ R $--subalgebra  of $ {\uqMg M } $  generated by  $ \, \left\{\,
F_i^{(\ell)}, \left( M_j ; \, c \atop t \right), M_j^{-1}, E_i^{(m)} \,
\bigg\vert \, \ell,c,t,m \in \N; i \in I, j \in I_\infty \,\right\} \, $;  this
is a  {\sl Hopf}  subalgebra of  $ {\uqMg M } $  (cf.~[Lu1], [DL]), with PBW basis over
$ R $  (cf.~[Ga2])
  $$  \left\{\, \prod_{\alpha \in \phitildep} E_\alpha^{(n_\alpha)} \cdot
\prod_{j \in I_\infty} \left( M_j \, ; \, 0 \atop t_j \right)
M_j^{-Ent({t_j/2})} \cdot \prod_{\alpha \in \phitildep} F_\alpha^{(m_\alpha)} \,
\Bigg\vert \, n_\alpha, t_j, m_\alpha \in \N, \; \forall\, \alpha, j \,\right\}
\; ;  $$
this is also a  $ \Cq $--basis  of  $ \uqMg{M} $,  hence  $ \gerUMg{M} $  is an
$ R $--form  of $ \uqMg{M} \, $.
                                                 \par
   Let  $ \, \calUMg{M} \, $  be the  $ R $--subalgebra of
$ \uqMg{M} $  generated by  $ \; \big\{\, \fbar_\alpha \, \big\vert
\, \alpha \in \phitildep \,\big\} \cup \big\{ M_j^{\pm 1} \, \big\vert
\, j \in I_\infty \,\big\} \cup \big\{\, \ebar_\alpha \, \big\vert \,
\alpha \in \phitildep \,\big\} \, $  (cf.~[BK], \S 3); this also is
a  {\sl Hopf}  subalgebra, with PBW basis
  $$  \left\{\, \prod_{\alpha \in \phitildep} \ebar_\alpha^{\, n_\alpha} \cdot
\prod_{j \in I_\infty} M_j^{t_j} \cdot \prod_{\alpha \in \phitildep}
\fbar_\alpha^{\, m_\alpha} \, \Bigg\vert \, t_i \in \Z, \, n_\alpha, m_\alpha
\in \N, \; \forall\, j, \alpha \,\right\}  $$
(over  $ R $);  the latter is also a  $ \Cq $--basis  of  $ \uqMg{M} $,  hence
$ \calUMg{M} $  is an  $ R $--form  of  $ \uqMg{M} $.
                                                 \par
   Like for quantum Borel algebras, the same forms can also be defined using
modified root vectors, hence they have also PBW bases of ordered monomials in
the  $ M_i $'s  and the modified root vectors.

\vskip7pt

   {\bf 3.4  Specialization at roots of 1 and quantum Frobenius
morphisms.}  \  In sight of specializations,  $ \C $  will be thought
of as an  $ R $--algebra  via  $ \, \C \cong R \big/ (q-c) \, $,  for
all  $ \, c \in \C \setminus \{0\} \, $.  Let  $ \varepsilon $  be
a primitive  $ \ell $--th  root of 1, for  $ \ell $  {\it odd\/},
$ \, \ell > d:= \max_i {\{d_i\}}_i \, $,  or  $ \, \ell = 1 \, $.
Then we set  $ \, \gerUepsilonMg{M} := \, \gerUMg{M} \Big/ (q -
\varepsilon) \, \gerUMg{M} \, \cong \, \gerUMg{M} \otimes_R \C \, $.
When  $ \, \ell =1 \, $  (i.~e.~$ \, \varepsilon = 1 \, $)  it is
well-known (cf.~e.~g.~[DL])\footnote{[{\it loc.~cit.}]  deals
with the finite type case and quantum groups of adjoint type,
but its arguments apply to the present situation too.}
that  $ \, \gerUunoMg{M} \, $  is a Poisson Hopf coalgebra, and
we have a Poisson Hopf coalgebra isomorphism
  $$  \gerUunoMg{M} \cong \ug \; ;   \eqno (3.4)  $$
in a word,  $ \gerUMg{M} $  specializes to  $ \ug $  for
$ \, q \rightarrow 1 \, $:  in symbols,  $ \, \gerUMg{M}
@>{q \rightarrow 1}>> \ug \, $.

\vskip5pt

   $ \underline{\hbox{\sl Remark}} $:  \, this point deserves some care.
Following Lusztig (cf.~for instance [Lu2]), people usually specialize
$ \uqg $  at  $ \, q = 1 \, $  by taking an integer form which is slightly
different from our  $ \, \gerUMg{M} \, $:  to be precise, it is defined like
$ \, \gerUMg{Q} \, $  but with  $ \, \left[ K_i ; \, c \atop t \right] :=
\prod_{s=1}^t {{q_i^{c-s+1} K_i - q_i^{-c+s-1} K_i^{-1}} \over  {q_i^s -
q_i^{-s}}} \, $  instead of  $ \, \left( K_i ; \, c \atop t \right) :=
\prod_{s=1}^t {{q_i^{c-s+1} K_i - 1} \over  {q_i^s - 1}} \, $  (hence
its toral part is smaller).  In particular one has  $ \, K_i - K_i^{-1}
= \left( q_i - q_i^{-1} \right) \cdot \left[ K_i ; \, 0 \atop 1 \right]
\, $,  which at  $ \, q = 1 \, $  gives  $ \, K_i^2 = 1 \, $:  thus
for  $ \, q \rightarrow 1 \, $  what one gets is only a central
extension of  $ \ug $,  namely  $ \; \ug \otimes_{\C} \C\,[K_0, K_1,
\dots, K_n] \, $,  \; with the  $ K_i $'s  central, group-like, and
idempotents.
                                             \par
   On the other hand, in  {\sl our}    $ \, \gerUMg{Q} \, $  (and similarly
for  $ \, \gerUMg{M} \, $,  in general)  we have  $ \, K_i - 1 = (q_i - 1)
\cdot \left( K_i ; \, 0 \atop 1 \right) \, $,  which at  $ \, q = 1 \, $
gives  $ \, K_i = 1 \, $,  so that (3.4) holds.

\vskip5pt

   When  $ \, \ell > 1 \, $  we have an epimorphism (cf.~[Lu1], [DL])
  $$  \gerFrg : \, \gerUepsilonMg{M} \llongtwoheadrightarrow \gerUunoMg{M} \cong
\ug   \eqno (3.5)  $$
of Hopf algebras defined by (for all  $ \, i \in I $,  $ j \in I_\infty \, $,
with  $ \, M_j := L_{\mu_j} \, $)
  $$  \gerFrg \colon
 \cases
   \! F_i^{(s)} \Big\vert_{q=\varepsilon} \!\!\!\! \mapsto F_i^{(s / \ell)}
\Big\vert_{q=1},  \left( M_j ; \, 0 \atop s \right) \!
\bigg\vert_{q=\varepsilon} \!\!\!\! \mapsto \left( M_j ; \, 0 \atop s / \ell
\right) \! \bigg\vert_{q=1}, E_i^{(s)} \Big\vert_{q=\varepsilon} \!\!\!\!
\mapsto E_i^{(s / \ell)} \Big\vert_{q=1}  \, \hbox{ if } \; \ell \Big\vert s  \\
   \! F_i^{(s)} \Big\vert_{q=\varepsilon} \!\!\! \mapsto 0,  \quad  \left(
M_j ; \, 0 \atop s \right) \! \bigg\vert_{q=\varepsilon} \!\!\! \mapsto 0,
\quad  E_i^{(s)} \Big\vert_{q=\varepsilon} \!\!\! \mapsto 0  \; \hbox{ \,
otherwise \, }   \hfill \hskip70pt  (3.6)  \\
   \! M_i^{-1} \Big\vert_{q=1} \!\!\! \mapsto 1  \\
 \endcases  $$
   \indent   As it is usual in literature,  we call
$ \gerFrg $  a  {\it
quantum Frobenius morphism\/}  (cf.~[Lu1]).
                                                  \par
   Similarly, we set  $ \, \calUepsilonMg{M} := \, \calUMg{M} \Big/
(q - \varepsilon) \, \calUMg{M} \, \cong \, \calUMg{M}
\otimes_R \C \, $;  when  $ \, \ell = 1 \, $  it is
  $$  \calUunoMg{M} \cong F \big[ \widehat{H}_{\scriptscriptstyle M} \big]
\eqno (3.7)  $$
as Poisson Hopf algebras over  $ \C \, $  (see [BK], \S 5): here
$ \widehat{H}_{\scriptscriptstyle M} $  is a Poisson proalgebraic group with
tangent Lie bialgebra  $ \hhat \, $:  it is simply an extension   --- by adding
two copies of  $ \C^* $  to the maximal torus ---   of the Poisson proalgebraic
group  $ \Omega $  which is defined in [BK], \S 4.  In a word,  $ \calUMg{M} $
specializes to  $ F \big[ \widehat{H}_{\scriptscriptstyle M} \big] $  as
$ \, q \rightarrow 1 \, $,  or  $ \, \calUMg{M} @>{\ q \rightarrow 1 \ }>> F
\big[ \widehat{H}_{\scriptscriptstyle M} \big] \, $.
                                                  \par
  When  $ \, \ell > 1 \, $,  assume in addition that the following technical
condition is satisfied:
  $$  \hbox{ $ \matrix
     g.c.d.(\ell,n+1) = 1  &  \qquad  \text{if  $ \gerg $  is of type
$ A_n $  or  $ C_n $}  \\
     g.c.d.(\ell,2n-1) = 1  &  \qquad  \text{if  $ \gerg $  is of type
$ B_n $}  \\
     g.c.d.(\ell,n-1) = 1  &  \qquad  \text{if  $ \gerg $  is of type
$ D_n $}  \\
     g.c.d.(\ell,3) = 1  &  \qquad  \text{if  $ \gerg $  is of type
$ E_6 $  or  $ G_2 $}  \\
               \endmatrix $ }   \eqno (3.8)  $$
Then from [BK], \S 3, we record the existence of a Hopf algebra monomorphism
  $$  \calFrg : \, F \big[ \widehat{H}_{\scriptscriptstyle M} \big]
\cong \calUunoMg{M} \llonghookrightarrow \calUepsilonMg{M}   \eqno (3.9)  $$
which we call again  {\it quantum Frobenius morphism\/};  it is uniquely
determined by
  $$  \calFrg : \quad \; \fbar_\alpha \Big\vert_{q=1} \mapsto
{\fbar_\alpha}^{\!\! \ell} \Big\vert_{q=\varepsilon} \, ,  \;\; L_\mu
\Big\vert_{q=1} \mapsto {L_\mu}^{\!\! \ell} \Big\vert_{q=\varepsilon} \, ,
\;\; \ebar_\alpha \Big\vert_{q=1} \mapsto {\ebar_\alpha}^{\!\! \ell}
\Big\vert_{q=\varepsilon}   \eqno (3.10)  $$
for all  $ \, \alpha \! \in \! \phitildep \, $,  $ \, \mu \! \in \!
M \, $  (this is a renormalized version of that in [BK]), and it enjoys
  $$  \calFrg : \; \fbar_{(r \delta,i)} \Big\vert_{q=1} \!\! \mapsto
\ell \, \fbar_{(r \ell \delta,i)} \Big\vert_{q=\varepsilon} \, ,
\; \ebar_{(r \delta,i)} \Big\vert_{q=1} \!\! \mapsto \ell \,
\ebar_{(r \ell \delta,i)} \Big\vert_{q=\varepsilon}  \quad  \forall
\; (r \delta,i) \in \phitildepim \, ;   \hskip11pt  (3.11)  $$
finally the image of  $ \calFrg $  is the center of  $ \calUepsilonMg{M} $,
i.e.
  $$  \calFrg \big( \calUunoMg{M} \big) = Z \big( \calUepsilonMg{M} \big) =:
Z_\varepsilon \; .   \eqno (3.12)  $$

\vskip1,7truecm

 \centerline{ \bf  \S \; 4 \,  Quantum formal groups }

\vskip10pt

   {\bf 4.1 Formal Hopf algebras and quantum formal groups.} \  In this
subsection we introduce the notion of  {\it quantum formal group\/}.
Recall (cf.~[Di], ch.~I) that formal groups can be defined in a category of a
special type of commutative topological algebras, whose underlying vector space
(or module) is linearly compact; following Drinfel'd's philosophy,we define
quantum formal groups by simply dropping out any commutativity assumption of the
classical notion of formal group; thus now we quickly outline how to modify the
latter (following [Di], ch.~I) in order to define our new quantum objects.
                                                 \par
   Let  $ E $  be any vector space over a field  $ K $  (one can then generalize
more or less wathever follows to the case of free modules over a ring), and let
$ E^* $  be its (linear) dual; we write  $ \, \langle x^*, x \rangle \, $  for
$ x^*(x) $  for  $ \, x \in E \, $,  $ \, x^* \in E^* \, $.  We consider on
$ E^* $  the  {\it weak  $ \ast $--topology},  i.~e.~the coarsest topology such
that for each  $ \, x \in E \, $  the linear map  $ \, x^* \mapsto \langle x^*,
x \rangle \, $  of  $ E^* $  into  $ K $  is continuous, when  $ K $  is given
the discrete topology.  We can describe this topology by choosing a basis
$ {\{e_i\}}_{i \in I} $  of  $ E \, $:  to each  $ i \in I $  we associate the
linear (coordinate) form  $ e_i^* $  on  $ E $  such that  $ \, \langle e_i^*,
e_j \rangle = \delta_{i j} \, $,  and we say that the family  $ {\{e^*_i\}}_{i
\in I} $  is the  {\it pseudobasis\/}  of  $ E^* $  dual to  $ {\{e_i\}}_{i \in
I} $;  then the subspace  $ E' $  of  $ E $  which is (algebraically)
generated by the  $ e_i^* $  is dense in  $ E^* $,  and  $ E^* $ is nothing
but the  {\it completion\/}  of  $ E' $,  when  $ E' $  is given the topology
for which a fundamental system of neighborhoods of  $ 0 $  consists of the
vector subspaces containing almost all the  $ e_i^* $;  thus elements of
$ E^* $  can be described by series in the  $ e_i^* $'s  which in the given
topology are in fact convergent.  Finally, the topological vector spaces
$ E^* $  are characterized by the property of linear compactness.
                                                 \par
   Let now  $ E $,  $ F $  be any two vector spaces over  $ K $,  and  $ \, u
\colon E \rightarrow F \, $  a linear map; then the dual map  $ \, u^* \colon
F^* \rightarrow E^* \, $ is  continuous, and conversely for any linear map  $ \,
v \colon F^* \rightarrow E^* \, $  which is continuous there exists a unique
linear map  $ \, u \colon E \rightarrow F \, $  such that  $ \, v = u^* \, $.
                                              \par
   The tensor product  $ E^* \otimes F^* \, $  is naturally identified to a
subspace of  $ {(E \otimes F)}^* $  by  $ \, \langle x^* \otimes y^*, x \otimes
y \rangle = \langle x^*, x \rangle \cdot \langle y^*, y \rangle \, $;  thus if
$ {\{e_i\}}_{i \in I} $  and  $ {\{f_j\}}_{j \in J} $  are bases of  $ E $  and
$ F $,  and  $ {\{e^*_i\}}_{i \in I} $  and  $ {\{f^*_j\}}_{j \in J} $  their
dual pseudobases in  $ E^* $  and  $ F^* $,  then  $ {\{e^*_i \otimes
f^*_j\}}_{i \in I, j \in J} $  is the dual pseudobasis of  $ {\{e_i \otimes
f_j\}}_{i \in I, j \in J} $  in  $ {(E \otimes F)}^* $.  Thus  $ {(E \otimes
F)}^* $  is the completion of  $ E^* \otimes F^* $  for the tensor product
topology, i.~e.~the topology of  $ E^* \otimes F^* $  for which a fundamental
system of neighborhoods of  $ 0 $  consists of the sets  $ \, E^* \otimes V +
W \otimes F^* \, $  where  $ V $,  resp.  $ W $,  ranges in a fundamental
system of neighborhoods of  $ 0 $  made of vector subspaces; we denote this
completion by  $ E^* \otimeshat F^* $,  and we call it the  {\it completed\/}
(or  {\it topological\/})  {\it tensor product\/}  of  $ E^* $  and  $ F^* $;
the embedding  $ \, E^* \otimes F^* \longhookrightarrow {(E \otimes F)}^* = E^*
\otimeshat F^* \, $  is then continuous.  Finally, if  $ \, u \colon E_1
\rightarrow  E_2 \, $,  $ \, v \colon F_1 \rightarrow F_2 \, $  are linear maps,
then  $ \, {(u \otimes v)}^* \colon {(E_2 \otimes F_2)}^* = {E_2}^* \otimeshat
{F_2}^* \longrightarrow {(E_1 \otimes F_1)}^* = {E_1}^* \otimeshat {F_1}^* \, $
coincides with the continuous extension to  $ {E_2}^* \otimeshat {F_2}^* $
of the continuous map  $ \, u^* \otimes v^* \colon {E_2}^* \otimes {F_2}^*
\rightarrow {E_1}^* \otimes  {F_1}^* \, $;  thus it is also denoted by  $ u^*
\otimeshat v^* \, $.
                                                  \par
   We define a  {\it linearly compact algebra}  to be a topological algebra
whose underlying vector space (or free module) is linearly compact: then
linearly compact algebras form a full subcategory of the category of topological
algebras; morever, for any two objects  $ A_1 $  and  $ A_2 $  in this category,
their topological tensor product  $ A_1 \otimeshat A_2 $  is defined.  Dually,
within the category of linearly compact vector spaces we define  {\it linearly
compact coalgebras}  as triplets  $ (C, \Delta, \epsilon) $  with  $ \, \Delta
\colon C \rightarrow C \otimeshat C \, $  and  $ \, \epsilon \colon C
\rightarrow K \, $  satisfying the usual coalgebra axioms.  The arguments in
[Di] (which never require commutativity nor cocommutativity) show that  $ \, {(\
)}^* \colon (A, m, 1) \mapsto (A^*, m^*, 1^*) \, $  defines a contravariant
functor from algebras to linearly compact coalgebras, while  $ {(\ )}^* \colon
(C, \Delta, \epsilon) \mapsto (C^*, \Delta^*, \epsilon^*) \, $  defines a
contravariant functor from coalgebras to linearly compact algebras.  Finally, we
define a  {\it formal Hopf algebra}  as a datum  $ (H, m, 1, \Delta, \epsilon,
S) $  such that $ (H, m, 1) $  is a linearly compact algebra,  $ (H, \Delta,
\epsilon) $  is a linearly compact coalgebra, and the usual compatibility axioms
of Hopf algebras are satisfied.  "Usual" Hopf algebras are particular cases of
formal Hopf algebras.
                                                 \par
   We define  {\bf quantum formal group}  the  {\it spectrum\/}  of a formal
Hopf algebra (whereas  {\it classical\/}  formal groups are spectra of  {\it
commutative\/}  formal Hopf algebras: cf.~[Di], ch.~I).
                                           \par
   Our goal is to study  $ \, {\uqMg M }^* \, $.  Since  $ {\uqMg M } $  is a
Hopf algebra, its linear dual  $ {\uqMg M }^* $  is a  {\sl formal}  Hopf
algebra.  The functor  $ {(\ \;)}^* $  turns the natural epimorphism  $ \,
{pr}_{\! \scriptscriptstyle M} \colon D_{\scriptscriptstyle M}
\longtwoheadrightarrow \uqMg{M} \, $ into a monomorphism  $ \,
j_{\scriptscriptstyle M'} \! := {\left( {pr}_{\! \scriptscriptstyle M}
\right)}^* \colon {\uqMg M }^* \longhookrightarrow {D_{\scriptscriptstyle M}}^*
\, $  of formal Hopf algebras: therefore we begin by studying
$ {D_{\scriptscriptstyle M}}^* $.  The following is straightforward:

\vskip7pt

\proclaim{Proposition 4.2}  Let  $ H_- $,  $ H_+ $  be Hopf
$ F $--algebras,  let \  $ \, \pi \colon {(H_-)}_{op} \otimes H_+
\longrightarrow F \, $  \  be an arbitrary Hopf pairing, and let
$ \, D:= D(H_-,H_+,\pi) \, $  be the corresponding quantum
double.  Then there exist  $ F $--algebra isomorphisms
  $$  D^* \cong {H_+}^* \otimeshat {H_-}^* \; ,  \quad  D^*
\cong {H_-}^* \otimeshat {H_+}^*  $$
dual of the  $F$--coalgebra  isomorphisms  $ \, D \cong H_+
\otimes H_- \, $,  $ \, D \cong H_- \otimes H_+ \, $  (cf.~\S 3.1).
$ \square $
\endproclaim

\vskip7pt

   {\bf 4.3  Quantum enveloping algebras as function algebras.}  \  The DRT
pairings induce several linear embeddings, namely
  $$  \matrix
   \um \longhookrightarrow {\up}^{\! *} \, ,  &
{im}_{\scriptscriptstyle M} \colon \,
{\uzM M } \longhookrightarrow {\uzM {M'} }^{\, *} \, ,  &
{\uMbm M } \longhookrightarrow {\uMbp {M'} }^{\, *}  &  \quad
\hbox{(induced by  $ \, \pi \, $)}  \\
   \up \longhookrightarrow {\um}^{\! *} \, ,  &
{\overline{im}}_{\scriptscriptstyle M} \colon \, {\uzM M }
\longhookrightarrow {\uzM {M'} }^{\, *} \, ,  &  {\uMbp M }
\longhookrightarrow {\uMbm {M'} }^{\, *}  &  \quad  \hbox{(induced by
$ \, \overline{\pi} \, $)}  \\
      \endmatrix   \eqno (4.1)  $$
the right-hand-side ones being also embeddings of formal Hopf algebras;
thus we identify the various quantum algebras with their images in the
corresponding dual spaces.

\vskip7pt

\proclaim{Lemma 4.4}
                                                \hfill\break
   \indent   (a) \  $ \calUm $  contains the pseudobasis of
$ {\up}^{\! *} $  dual of the PBW basis of  $ \gerUp $  of decreasing ordered
monomials, and  $ \gerUm $  contains the pseudobasis of  $ {\up}^{\! *} $  dual
of the PBW basis of  $ \calUp $  of decreasing ordered monomials.  Moreover, the
PBW bases of  $ \calUm $  and of  $ \gerUm $  are also pseudobases of
$ {\up}^{\! *} $.  A similar statement holds with the roles of  $ \um $  and
$ \up $  reversed.
                                                \hfill\break
   \indent   (b) \  $ {\calUzM M } $  (hence  $ {\uzM M } \, $)
contains the pseudobasis  $ {\Cal B}_{\scriptscriptstyle M} $  (relative to
$ {im}_{\scriptscriptstyle M} $),  resp.~$ {\overline{\Cal
B}}_{\scriptscriptstyle M} $  (relative to
$ {\overline{im}}_{\scriptscriptstyle M} $),  of  $ {\uzM {M'} }^{\, *} $
dual of the PBW basis of  $ \gerUzM{M'} $.
                                              \hfill\break
   \indent   (c) \  $ {\calUMbm M } \, $,  resp.~$ {\calUMbp M } \, $  (hence
$ {\uMbm M } \, $,  resp.~$ {\uMbp M } \, $)  contains the pseudobasis of
$ {\uMbp {M'} }^{\, *} $,  resp.~of  $ {\uMbp {M'} }^{\, *} $,  dual of the PBW
basis of  $ \gerUMbm{M'} \, $,  resp.~of  $ \gerUMbp{M'} \, $.  Moreover, the
PBW basis of  $ \calUMbm{M} $,  resp.~$ {\calUMbp M } \, $,  is itself a
pseudobasis of  $ {\uMbp {M'} }^{\, *} $,  resp.~of  $ {\uMbp {M'} }^{\, *} $.
\endproclaim

\demo{Proof}  Part  {\it (a)}  of the statement follows from Lemma 2.5 and the
{\sl Claim}  in \S 2.7.
                                                      \par
   As for  {\it (b)},  let  $ \, u_\tau := \prod_{i \in \infty}
\left( \Lambda_i; \, 0 \atop t_i \right) \cdot
\Lambda_i^{-Ent({t_i/2})} \, $  ($ \tau = (t_0, t_1,
\dots, t_n, t_\infty) \in \N^{n+2} \, $)  be any monomial
in the PBW basis of  $ \gerUzM{M'} $.  Then direct computation gives
  $$  {\Big\langle L_{-\mu}, u_\tau \Big\rangle}_{\! \pi} =
\prod_{i \in \infty} {m_i \choose t_i}_{\!\! q_i} \cdot
q^{ - d_i m_i \cdot Ent(t_i/2)}  \; \qquad  \forall\, \mu,
\tau \in \N^{n+2}  $$
where we identify  $ \, M_+ \cong \N^{n+2} \, $  so that
$ \, M_+ \ni \mu = m_0 \, \mu_0 + m_1 \, \mu_1 + \cdots + m_n \,
\mu_n + m_\infty \, \mu_\infty \cong (m_0, m_1, \dots, m_n, m_\infty)
\in \N^{n+2} \, $.  Then endowing  $ \N^{n+2} $  with the product
ordering (of the natural ordering of  $ \N $)  we have
  $$  \eqalign{
   {\Big\langle L_{-\mu}, u_\tau \Big\rangle}_{\! \pi}  &  \neq 0  \iff  \tau
\preceq \mu  \cr
   {\Big\langle L_{-\tau}, u_\tau \Big\rangle}_{\! \pi}  &  = q^{-T(\tau)}
\qquad \qquad \qquad \forall\, \tau \in \N^n  \cr }  $$
where  $ \, T(\tau) := \sum_{i \in I_\infty} d_i t_i Ent(t_i/2) \, $;  in
particular  $ \, q^{-T(\tau)} \, $  {\sl is invertible in}  $ R $.  Thus we have formulas  (for all  $ \, \tau \in \N^{n+2} \, $)
  $$  {} \;  L_{-\tau} = q^{-T(\tau)} \cdot u_\tau^{\,*} + \sum_{\tau' \prec
\tau} {\Big\langle L_{-\tau}, u_{\tau'} \Big\rangle}_{\! \pi} \cdot
u_{\tau'}^{\,*}  $$
which tell us that  $ \, \big\{\, L_{-\tau} \,\big\vert\, \tau \in \N^{n+2}
\,\big\} \, $  is obtained from  $ \, \big\{\, u_\tau^{\,*} \,\big\vert\, \tau
\in \N^{n+2} \,\big\} \, $  by means of the matrix  $ \, {\Bbb M} :=  {\Big(
{\big\langle L_{-\tau}, u_{\tau'} \big\rangle}_{\! \pi} \Big)}_{\tau, \tau' \in
\N^{n+2}} $  which has lower triangular shape, all entries in  $ R $,  and
diagonal entries  {\sl invertible}  in  $ R \, $;  then the inverse matrix
$ {\Bbb M}^{-1} $  has the same properties, whence  {\it (b)}  follows.
                                                      \par
   Finally, for  {\it (c)}  note that  $ \, \pi \big( y \cdot \ell, x \cdot m
\big) = \pi(y,x) \cdot \pi(\ell,m) \, $  for all  $ \, y \in \um \, $,  $ \,
\ell \in \uzM{M} \, $,  $ \, m \in \uzM{M'} \, $,  $ \, x \in \up \, $,  as one
sees at once from Lemma 2.5 (and similarly for  $ \overline{\pi} \, $);
therefore, since the PBW basis of, say,  $ \gerUMbm{M'} $  is the tensor
product
of the like bases of  $ \gerUp $  and  $ \gerUzM{M'} $  then the dual
pseudobasis of  $ {\left( \gerUMbm{M'} \right)}^* $  is the tensor product of
the pseudobases of  $ {\gerUp}^* $  and  $ {\left( \gerUzM{M'} \right)}^* $
given in  {\it (a)}  and  {\it (b)}.  With similar arguments all of  {\it (c)}
is proved.   $ \square $
\enddemo

\vskip7pt

   {\bf 4.5 Remark.}  \;  Since  $ \, D_{\scriptscriptstyle {M'} }
\cong \uMbp{M'} \otimes \uMbm{Q_\infty} \cong \up \otimes \uzM{M'}
\otimes \uzM{Q_\infty} \otimes \um \, $,  we have  $ \,
{D_{\scriptscriptstyle {M'} }}^* \cong {\uMbp {M'} }^* \otimeshat
{\uMbm{Q_\infty} }^* \cong {\up}^* \otimeshat {\uzM {M'} }^*
\otimeshat {\uzM{Q_\infty} }^* \otimeshat {\um}^* \, $;
hence from Lemma 4.4 we deduce
                                                       \par
   {\it  Every element  $ \, f \in {D_{\scriptscriptstyle {M'} }}^* \, $  has
a unique expression as formal series
  $$  f = \sum_{{\Cal F}, {\Cal M}, {\Cal L}, {\Cal E}} a_{{\Cal F}, {\Cal M},
{\Cal L}, {\Cal E}} \cdot \calF \cdot {\Cal M } \cdot {\Cal L} \cdot \calE  $$
in which  $ \, a_{{\Cal F}, {\Cal M}, {\Cal L}, {\Cal E}} \in \Cq \, $,
$ \, {\Cal M} \in {\Cal B}_{\scriptscriptstyle M} \, $,  $ \, {\Cal L} \in
{\overline{\Cal B}}_{\scriptscriptstyle P_\infty} \, $,  and the  $ \calF$'s,
resp.~the  $ \calE $'s,  are ordered monomials in the  $ F_\alpha $'s,
resp.~in
the  $ E_\alpha $'s;  then the natural evaluation pairing  $ \,
{D_{\scriptscriptstyle {M'} }}^* \otimes D_{\scriptscriptstyle {M'} }
\rightarrow \Cq \, $  is given by  $ \, \pi \otimes \pi \otimes \overline{\pi}
\otimes \overline{\pi} \, $.
                                                    \par
   In particular, every  $ \, f \in {D_{\scriptscriptstyle {M'} }}^* \, $  can
be uniquely expressed as a formal series in the  $ F_\alpha $'s  and the
$ E_\alpha $'s  ($ \, \alpha \in \phitildep \, $)  with coefficients in  $ \,
{\left( {\uzM {M'} } \otimes {\uzM{Q_\infty} } \right)}^* \cong {\uzM {M'} }^*
\otimeshat {\uzM{Q_\infty} }^* \, $.}
                                                    \par
   Similarly the triangular decompositions  $ \, \up \otimes {\uzM {M'} }
\otimes \um \cong {\uqMg {M'} } \cong \um \otimes {\uzM {M'} } \otimes \up \, $
give  $ \, {\up}^{\! *} \otimeshat {\uzM {M'} }^* \otimeshat {\um}^{\! *} \cong
{\uqMg {M'} }^{\! *} \cong {\um}^{\! *} \otimeshat {\uzM {M'} }^* \otimeshat
{\up}^{\! *} \, $,  whence Lemma 4.4 implies
                                                     \par
   {\it  Every  $ \, f \in {\uqMg {M'} }^* \, $  can be uniquely expressed as
a formal series in the  $ F_\alpha $'s  and the  $ E_\alpha $'s  ($ \, \alpha
\in \phitildep \, $)  with coefficients in  $ \, {\uzM {M'} }^* \, $;  then the
natural evaluation pairing  $ \, {\uqMg {M'} }^* \otimes {\uqMg {M'} }
\rightarrow \Cq \, $  can be described by  $ \, \pi \otimes \overline{\pi}
\otimes \overline{\pi} \, $  or by  $ \,  \pi \otimes \overline{\pi} \otimes
\overline{\pi} \, $.}
                                                    \par
  In particular, a pseudobasis of  $ {\uqMg {M'} }^* $  can be given
by taking the tensor product of pseudobases   --- for instance, those
provided by Lemma 4.4 ---   of  $ {\up}^{\! *} $,  $ {\uzM {M'} }^* $,
and  $ {\um}^{\! *} $;  moreover, the (tensor) PBW basis of  $ \,
\calUm \otimes \calUzM{M} \otimes \calUp \, $  is a pseudobasis
of  $ {\uqMg {M'} }^* $.
                                                      \par
   In the sequel when considering the composed embedding  $ \, {\uzM M }
\hookrightarrow {\uzM {M'} }^* \hookrightarrow {\uqMg {M'} }^* \, $  we
shall always mean that the first embedding is induced by
$ \overline{\pi} $  (cf.~(4.1)); accordingly, the evaluation pairing  $ \,
{\uqMg {M'} }^* \otimes {\uqMg {M'} } \rightarrow \Cq \, $  will be described
by  $ \,  \pi \otimes \overline{\pi} \otimes \overline{\pi} \, $,  hence the
pseudobasis of  $ {\uzM {M'} }^* $  to use is  $ {\overline{\Cal
B}}_{\scriptscriptstyle M} $.

\vskip7pt

\proclaim{Proposition 4.6}  The monomorphism  $ \; j_{\scriptscriptstyle M}
\colon \, {\uqMg {M'} }^* \longhookrightarrow {D_{\scriptscriptstyle
M'}}^{\! *} \; $  (cf.~\S 4.1) is given by
  $$   j_{\scriptscriptstyle M} \colon  \; \quad  F_i \mapsto F_i \otimes 1 \,
,
 \quad  L_\mu \mapsto L_{-\mu} \otimes L_\mu \, ,  \quad  E_i \mapsto 1 \otimes
E_i  \; \qquad  \forall \, i, \mu \, ;   \eqno (4.2)  $$
in particular the image of  $ j_{\scriptscriptstyle M} $  is the closure of
the subalgebra generated by the set
  $$   \Big\{\, F_i \otimes 1, \, L_{-\mu} \otimes L_\mu, \, 1 \otimes E_i
\,\Big\vert\, i\in I_\infty, \, \mu \in M \,\Big\} \; .  $$
\endproclaim

\demo{Proof}  For PBW monomials we have  $ \, {pr}_{\scriptscriptstyle M}
\Big( E \cdot L \otimes K \cdot F \Big) =  E \cdot L \cdot K \cdot F \, $;
therefore (4.2) comes out of the definition  $ \, j_{\scriptscriptstyle M} :=
{\left( {pr}_{\scriptscriptstyle M} \right)}^{\! *} \, $.  As an example
  $$  \displaylines{
   \Big\langle j_{\scriptscriptstyle M} \big( L_\mu \big), E \cdot L_\nu
\otimes K_\alpha \cdot F \Big\rangle = \Big\langle L_\mu, E \cdot L_\nu
\cdot K_\alpha \cdot F \Big\rangle = \delta_{E,1}
\cdot \delta_{F,1} \cdot q^{(\mu \vert \nu + \alpha)}  \cr
   {\Big\langle L_{-\mu} \otimes L_\mu, E \cdot
L_\nu \otimes K_\alpha \cdot F \Big\rangle}_{\! \pi \otimes
\overline{\pi}} = \delta_{E,1} \cdot \delta_{F,1} \cdot
q^{(\mu \vert \nu + \alpha)}  \cr }  $$
whence  $ \, j_{\scriptscriptstyle M} \big( L_\mu \big) =
L_{-\mu} \otimes L_\mu \, $.  Since  $ \, j_{\scriptscriptstyle M} := {\left(
{pr}_{\scriptscriptstyle M'} \right)}^{\! *} \, $  is continuous (cf.~\S 4.1),
by Lemma 4.4 and Remark 4.5 it is uniquely determined by (4.2).
$ \square $
\enddemo

\vskip7pt

{\bf  Remark 4.7.}  \  Thanks to the previous results, we can identify
$ j_{\scriptscriptstyle M} \left( {\uqMg {M'} }^* \right) $  with the
space of formal series in the  $ F_\alpha $'s  and the  $ E_\alpha $'s
($ \, \alpha \in \phitildep \, $)  with coefficients in  $ \,
{\uzM {M'} }^* \, $,  using Proposition 4.6 and the similar
identification for  $ {\uqMg {M'} }^* $  (cf.~Remark 4.5).
Then similar remarks to those in Remark 4.5 hold.

\vskip7pt

   {\bf 4.8  Integer forms.}  \;  We want to study the subspaces of linear
functions on  $ {\uqMg {M'} } $  which are "integer-valued" on its
integer forms.  Thus we define
  $$  \displaylines {
   {\gerU^{\scriptscriptstyle M'}\!(\ghat)}^* \!\! := \Big\{ f \in
{U_q^{\scriptscriptstyle M'}\!(\ghat) }^* \,\Big\vert
\left\langle f, \gerU^{\scriptscriptstyle M'}\!(\ghat)
\right\rangle \!\subseteq\! R \Big\} ,  \,  \hfill
{{\Cal U}^{\scriptscriptstyle M'}\!(\ghat)}^* \!\! := \Big\{ f \in
{U_q^{\scriptscriptstyle M'}\!(\ghat)}^* \,\Big\vert
\left\langle f, {\Cal U}^{\scriptscriptstyle M'}\!(\ghat)
\right\rangle \!\subseteq\! R \Big\}  \cr
   {\frak I}^{\scriptscriptstyle M} \! := \Big\{ f \in
j_{\scriptscriptstyle M} \! \Big( \hskip-0,4pt
{U_q^{\scriptscriptstyle M'}\!(\ghat)}^* \Big) \Big\vert
\left\langle f, \gerU^{\scriptscriptstyle M'}\!(\ghat)
\right\rangle \!\subseteq\! R \Big\} ,
\hskip1,8pt  \hfill  {\Cal I}^{\scriptscriptstyle M} \! :=
\Big\{ f \in j_{\scriptscriptstyle M} \! \Big( \hskip-0,4pt
{U_q^{\scriptscriptstyle M'}\!(\ghat)}^* \Big) \Big\vert
\left\langle f, {\Cal U}^{\scriptscriptstyle M'}\!(\ghat)
\right\rangle \!\subseteq\! R \Big\}  \cr
   {\gerUzM {M'} }^* := \Big\{\, f \in {\uzM {M'} }^* \,\Big\vert\,
\left\langle
f, {\gerUzM {M'} } \right\rangle \subseteq R \,\Big\} \, ,  \qquad  {\calUzM
{M'} }^* := \Big\{ f \in {\uzM{M'} }^* \,\Big\vert\, \left\langle f, {\calUzM
{M'} } \right\rangle \subseteq R \,\Big\} \, ;  \cr }  $$
notice that  $ j_{\scriptscriptstyle M} $  restricts
to isomorphisms  $ \; j_{\scriptscriptstyle M} \colon \, {\gerUMg
{M'} }^* \, {\buildrel \cong \over \loongrightarrow} \, {\frak
I}^{\scriptscriptstyle M} \, $,  $ \; j_{\scriptscriptstyle M}
\colon \, {\calUMg {M'} }^* \, {\buildrel \cong \over \loongrightarrow}
\, {\Cal I}^{\scriptscriptstyle M} \, $.

\vskip7pt

\proclaim{Proposition 4.9}
                                      \hfill\break
   \indent   (a)  $ {\gerUMg {M'} }^* $  is the  $ R $--submodule
(of  $ {\uqMg {M'} }^* \, $)
                of all formal series (cf.~\S 4.5)\break
$ \; \sum_{{\Cal F}, \psi, {\Cal E}} \calF \cdot \psi \cdot {\Cal E} \; $  in
which  $ \, \psi \in {\gerUzM {M'} }^* \, $  and the  $ \calF $'s,  resp.~the
$ \calE $'s,  are monomials of the PBW basis of  $ \calUm $,  resp.~of
$ \calUp $.  So  $ {\gerUMg {M'} }^* $  is a formal Hopf subalgebra
of  $ {\uqMg {M'} }^* $.
                                      \hfill\break
   \indent   (b)  $ {\calUMg {M'} }^* $  is the  $ R $--submodule  (of  $ {\uqMg
{M'} }^* \, $)
                        of all formal series (cf.~\S 4.5)\break
$ \; \sum_{{\frak F}, \phi, {\frak E}} {\frak F} \cdot \phi \cdot {\frak E}
\; $  in which  $ \, \phi \in {\calUzM {M'} }^* \, $  and the  $ {\frak F} $'s,
resp.~the  $ {\frak E} $'s,  are monomials of the PBW basis of  $ \gerUm $,
resp.~of  $ \gerUp $.  So  $ {\calUMg {M'} }^* $  is a formal Hopf
subalgebra of  $ {\uqMg {M'} }^* $.
\endproclaim

\demo{Proof}  Let us prove  {\it (b)},  the proof for  {\it (a)}  being
completely similar.
                                                          \par
   Let  $ \, f \in {\uqMg {M'} }^* \, $  be given, and expand it as a series
$ \; f = \sum_{\phi, \eta} \gerF_\phi \cdot \varPhi_{\phi,\eta} \cdot \gerE_\eta
\; $  in which the  $ \gerF_\phi  $'s,  resp.~the  $ \gerE_\eta $'s,  are PBW
monomials of  $ \gerUm $,  resp.~of  $ \gerUp $,  and  $ \, \varPhi_{\phi,\eta}
\in {\uzM {M'} }^* \, $.  Now fix  $ \bar\eta $,  and let  $ \calE_{\bar\eta}
\, $  be the PBW-like monomial of  $ \calUp $  which is (up to  $ \pm q^s $,
for some  $ s \in \Z $)  the dual of  $ \gerF_{\bar\eta} \, $:  according to
Lemma 2.5, this is an ordered product of the  $ \ebar_\gamma $'s  ($ \, \gamma
\in \phipre \, $)  and of the duals of the  $ {\, r \, \over \, {[r]}_{q_i} \,}
F_{(r \delta,i)} $'s  ($ \, (r \delta,i) \in \phitildepim \, $);  similarly, fix
$ \bar\phi $  and let  $ \calF_{\bar\phi} \, $  be the PBW-like monomial of
$ \calUm $  which is the dual (up to  $ \pm q^t $,  for some  $ t \in \Z $) of
$ \gerE_{\bar\phi} \, $:  thanks to the  {\sl Claim}  in \S 2.7, we
have  $ \, \calE_{\bar\eta}, \calF_{\bar\phi} \in \calUMg{M} \, $.
Fix also  $ \, \nu \in M' \, $:  then
  $$  \Big\langle f, \, \calE_{\bar\eta} \cdot L_\nu \cdot \calF_{\bar\phi}
\Big\rangle = \sum_{\phi,\eta} \left\langle \gerF_\phi \cdot
\varPhi_{\phi,\eta} \cdot \gerE_\eta, \, \calE_{\bar\eta} \cdot L_\nu \cdot \calF_{\bar\phi}
\right\rangle = \pm q^{s+t} \left\langle \varPhi_{\bar\phi,\bar\eta}, L_\nu
\right\rangle  $$
   \indent   Therefore, since the monomials  $ \, \calE_{\bar\eta} \cdot L_\nu
\cdot \calF_{\bar\phi} \, $  form a (PBW-like) basis of  $ {\calUMg {M'} }
\, $,  we have  $ \, f \in {\uqMg {M'} }^* \, $  if and only if  $ \,
\varPhi_{\bar\phi,\bar\eta} \in {\calUzM {M'} }^* \, $,  q.e.d.
                                                            \par
   Now consider the Hopf structure.  Let  $ \, f \in {\calUMg {M'} }^* \, $,
and expand  $ \Delta(f) $
                                         as a series\break
$ \; \Delta(f) = \sum_\sigma \left( \gerF_\sigma \cdot \varPhi_\sigma
\cdot \gerE_\sigma \right) \otimes \left( \gerF'_\sigma \cdot \varPhi'_\sigma
\cdot \gerE'_\sigma \right) \; $  so that  $ \, \varPhi_\sigma \otimes
\varPhi'_\sigma \neq \varPhi_\tau \otimes \varPhi'_\tau \, $
($ \, \in {\uzM {M'} }^* \otimeshat {\uzM {M'} }^* \! = {\left( \uzM{M'} \otimes \uzM{M'}
\right)}^* \, $)  for all  $ \sigma $,  $ \tau $,  such that  $ \, \left(
\gerF_\sigma, \gerE_\sigma, \gerF'_\sigma, \gerE'_\sigma \right) \neq \left(
\gerF_\tau, \gerE_\tau, \gerF'_\tau, \gerE'_\tau \right) \, $  (this is always
possible).  As  $ \, f \in {\calUMg {M'} } \, $,  then  $ \Delta(f) $  is
integer-valued on  $ \, {\calUMg {M'} } \otimes {\calUMg {M'} } $.  Fix any
$ \bar\sigma $:  as above there exist suitable PBW-like monomials  $ \calE_{\bar
\sigma} $,  $ \calF_{\bar \sigma} $,  $ \calE'_{\bar \sigma} $,
$ \calF'_{\bar\sigma} $  such that
  $$  \left\langle \Delta(f), \big( {\Cal E}_{\bar\sigma} \otimes
{\Cal E}'_{\bar\sigma} \big) \cdot \big( L_\nu \otimes L_{\nu'} \big) \cdot
\big( \calF_{\bar\sigma} \otimes \calF'_{\bar\sigma} \big) \right\rangle =
\pm q^z \cdot \big\langle \phi_\sigma \otimes \phi'_\sigma, L_\nu \otimes
L_{\nu'} \big\rangle  $$
for all  $ \, \nu, \nu' \in M' \, $  (for some  $ \, z \in \Z \, $);  since
$ \Delta(f) $  is integer-valued,  $ \phi_{\bar \sigma} \otimes \phi'_{\bar
\sigma} $  is integer-valued on  $ {\calUzM {M'} } \otimes {\calUzM {M'} } $,
that is  $ \, \phi_{\bar \sigma} \otimes \phi'_{\bar \sigma} \in {\left(
{\calUzM {M'} } \otimes {\calUzM {M'} } \right)}^* = {\calUzM {M'} }^*
\otimeshat {\calUzM {M'} }^* $;  but  $ \, \phi_{\bar \sigma} \otimes
\phi'_{\bar\sigma} \in {\uzM {M'} }^* \otimes {\uzM {M'} }^* $,  thus  $ \,
\phi_{\bar \sigma}, \phi'_{\bar\sigma} \in {\calUzM {M'} }^* $,  \, q.e.d.
                                               \par
   Finally, we have  $ \, 1 \in {\calUMg{M'} }^* \, $,  because  $ \, 1
:= \epsilon \, $,  $ \, \epsilon \big( {\calUMg{M'} }^*
\big) \subseteq R \, $  because  $ \, \epsilon := 1^* \, $  and  $ \, 1
\in {\calUMg{M'} } \, $,  and  $ \, S \big( {\calUMg{M'} }^* \big) =
{\calUMg{M'} }^* \, $  because  $ \, S := S^* \, $  and  $ \, S \big(
{\calUMg{M'} } \big) = {\calUMg {M'} } \, $.  Thus  $ {\calUMg{M'} }^* $  is a
formal Hopf subalgebra of  $ {\uqMg{M'} }^* $,
q.e.d.   $ \square $
\enddemo

\vskip7pt

\proclaim{Definition 4.10} \  We call  $ A^{\scriptscriptstyle M} $  the
subalgebra of  $ \, \uMbm{M} \otimes \uMbp{P_\infty} \; \big( \, \subset
{D_{\scriptscriptstyle M'}}^* \, \big) \, $  generated by  $ \; \big\{ \, F_i
\otimes 1, \, L_{-\mu} \otimes L_\mu, \, 1 \otimes E_i \,\big\vert\, i \in I, \,
\mu \in M \,\big\} \, $.  Then we set
  $$  \eqalign {
     {\frak A}^{\scriptscriptstyle M}  &  := \big\{\, f \in
A^{\scriptscriptstyle M} \,\big\vert\, \big\langle f, {\gerUMg {M'} }
\big\rangle \subseteq R \,\big\} = A^{\scriptscriptstyle M} \cap {\frak
I}^{\scriptscriptstyle M}  \cr
     {\Cal A}^{\scriptscriptstyle M}  &  := \big\{\, f \in
A^{\scriptscriptstyle
M} \,\big\vert\, \big\langle f, {\calUMg {M'} } \big\rangle \subseteq R
\,\big\}
= A^{\scriptscriptstyle M} \cap {\Cal I}^{\scriptscriptstyle M} \, .  \cr }  $$
\endproclaim

\vskip7pt

\proclaim{Lemma 4.11}
                                               \hfill\break
   \indent   (a) \  $ {\frak A}^{\scriptscriptstyle M} $  is an  $ R $--integer
form of  $ A^{\scriptscriptstyle M} $,  generated as an  $ R $--subalgebra by
  $$  \left\{\, \fbar_\alpha \otimes 1, \, L_{-\mu} \otimes L_\mu, \, 1 \otimes
\ebar_\alpha \,\Big\vert\, \alpha \in \phitildep, \, \mu \in M \,\right\}
\, .  $$
   \indent   (b) \  $ {\Cal A}^{\scriptscriptstyle M} $  is an  $ R $--integer
form of  $ A^{\scriptscriptstyle M} $,  generated as an  $ R $--subalgebra by
  $$  \bigg\{ F_i^{(a)} \otimes 1, \, \left( L_{-\mu_j} \otimes L_{\mu_j}; \,
c \atop t \right), \, L_{\mu_j} \otimes L_{-\mu_j}, \, 1 \otimes E_i^{(d)}
\;\bigg\vert\; i \in I, \, j \in I_\infty \, ; \; a, t, d \in \N \, ; \; c \in
\Z \,\bigg\} .  $$
\endproclaim

\demo{Proof}  Definitions yield a linear isomorphism  $ \,
\Phi_{\scriptscriptstyle M} \colon A^{\scriptscriptstyle M} \,
{\buildrel \cong \over \longrightarrow} \, \um \otimes {\uzM M }
\otimes \up \, $
                                       given by\break
$ \Phi_{\scriptscriptstyle M} \colon  \; F_i \! \otimes \! 1 \mapsto F_i \!
\otimes \! 1 \! \otimes \! 1, \, L_{-\mu} \! \otimes \! L_\mu \mapsto 1 \!
\otimes \! L_\mu \! \otimes \! 1, \, 1 \! \otimes \! E_i \mapsto 1 \! \otimes \!
1 \! \otimes \! E_i \, $;  but this restricts to  $ \, \Phi_{\scriptscriptstyle
M} \colon \, {\Cal A}^{\scriptscriptstyle M} \, {\buildrel \cong \over
\longrightarrow} \, \calUm \otimes {\calUzM M } \otimes \calUp \, $,  $ \,
\Phi_{\scriptscriptstyle M} \colon \, {\frak A}^{\scriptscriptstyle M} \,
{\buildrel \cong \over \longrightarrow} \, \gerUm \otimes {\gerUzM M } \otimes
\gerUp \, $,  so \S 3.3 yields the claim.   $ \square $
\enddemo
%
%
%
%

\vskip7pt

   {\bf 4.12  Gradings.}  \  Recall (cf.~\S 2.1) that  $ \uMbp{M} $
has a  $ Q_+ $--grading  $ \, \uMbp{M} = \oplus_{\alpha \in Q_+}
{\left( \uMbp{M} \right)}_{\! \alpha} \, $  given by decomposition
in direct sum of weight spaces for the adjoint action of
$ {\uzM M } $;  also  $ {\uMbm M } $  has an analogous
$ Q_- $--grading.  These are gradings of Hopf algebras (in the usual
obvious sense), inherited by the integer forms, and DRT pairings respect them,
that is e.g.  $ \, \pi \left({\left( {\uMbm M } \right)}_{\! \beta}, {\left(
{\uMbp {M'} } \right)}_{\! \gamma} \right) = 0  \, $  for all  $ \, \beta \in
Q_- $,  $ \gamma \in Q_+ \, $  such that  $ \, \beta + \gamma \neq 0 \, $.
                                                    \par
   The gradings of quantum Borel subalgebras induce a  $ Q $--grading of the
Hopf algebra  $ \, D_{\scriptscriptstyle M} := \uMbp{M} \otimes \uMbm{Q_\infty}
\, $  (inherited by its quotient Hopf algebra  $ {\uqMg M } $),  where the
subspace  $ \, {\left( {\uMbp M } \right)}_{\! \beta} \otimes {\left(
\uMbm{Q_\infty} \right)}_{\! \gamma} \, $  has degree  $ \beta + \gamma \, $,
and also a  $ Q $--grading  of the subalgebra  $ \, \uMbp{M} \otimes
\uMbm{Q_\infty} \, $  of  $ {D_{\scriptscriptstyle M'}}^{\! *} \, $;  since
$ {D_{\scriptscriptstyle M'}}^{\! *} $  is a completion (via formal series) of
this subalgebra, it inherits on its own sort of a "pseudograding", in the sense
that every element of  $ {D_{\scriptscriptstyle M'}}^{\! *} $  is a (possibly
infinite) sum of terms each of whom has a well-defined degree: namely, given
$ \, f \in {D_{\scriptscriptstyle M'}}^{\! *} \, $  with formal series expansion
(cf.~Remark 4.5)  $ \, f = \sum_{F,\phi,E} F \cdot \phi \cdot E \, $  (where
$ \, \phi \in {\left( {\uzM {M'} } \otimes \uzM{Q_\infty} \right)}^* \, $  and
the  $ F $'s  and the  $ E $'s  are PBW monomials), we define the degree of its various summands by

  $$  deg \big( F \cdot \phi \cdot E \big) := deg \big( F \big) + deg \big( E
\big)  $$
where  $ \, deg \left( \prod_{\alpha \in \phitildep} F_\alpha^{\, f_\alpha}
\right) := - \sum_{\alpha \in \phitildep} f_\alpha \alpha \, $,  $ \, deg \left(
\prod_{\alpha \in \phitildep} E_\alpha^{\, e_\alpha} \right) := \sum_{\alpha \in
\phitildep} e_\alpha \alpha \, $  (this degree is again a weight for a suitable
action of  $ {\uzM M } $  on  $ \, {\uMbm M } \otimes {\uMbp {P_\infty} }
\, $).  Now  $ \, {\uMbm M } \otimes {\uMbp {P_\infty} } \, $  is dense in
$ {D_{\scriptscriptstyle M'}}^{\! *} $,  and the restriction of the pairing
$ \, {D_{\scriptscriptstyle M'}}^{\! *} \otimes D_{\scriptscriptstyle M'}
\rightarrow \Cq \, $  to  $ \, \left( \uMbm{M} \otimes {\uMbp {P_\infty} }
\right) \otimes \left( {\uMbp {M'} } \otimes \uMbm{Q_\infty} \right) \, $  is
nothing but  $ \, \big( \pi \otimes \overline{\pi} \big) \smallcirc \tau_{2,3}
\, $  (with  $ \, \tau_{2,3} \colon \, x \otimes y \otimes z \otimes w \mapsto x
\otimes z \otimes y \otimes w \, $;)  therefore, since  $ \pi $  and
$ \overline{\pi} $  respect the gradings, also the pairing   $ \,
{D_{\scriptscriptstyle M'}}^{\! *} \otimes D_{\scriptscriptstyle M'} \rightarrow
\Cq \, $  respects the pseudogradings we are dealing with.
                                                        \par
   Finally, the pseudograding of  $ {D_{\scriptscriptstyle M'}}^{\! *} $  is
compatible with the formal Hopf structure.  For example, look at  $ S(x) $,  for
homogeneous  $ \, x \in {D_{\scriptscriptstyle M'}}^{\! *} \, $:  given
homogeneous  $ \, y \in D_{\scriptscriptstyle M'} \, $,  we have  $ \,
\big\langle S(x), y \big\rangle = \big\langle x, S(y) \big\rangle = \big\langle
x, y' \big\rangle \, $  where  $ \, y' := S(y) \, $  is homogeneous on its own
of degree  $ \, deg(y') = deg(y) \, $  (for the grading of
$ D_{\scriptscriptstyle M'} $  is compatible with the Hopf structure);
therefore we have  $ \, deg\big(S(x)\big) = deg(x) \, $,  because
  $$  \big\langle S(x), y \big\rangle \neq 0 \Longrightarrow deg(y) =
deg(y') = deg(x) \Longrightarrow S(x) \in {\left( {D_{\scriptscriptstyle
M'}}^{\! *} \right)}_{deg(x)} \, .  $$

\vskip7pt

   {\bf 4.13  Some umbral calculus.}  \;  In this section we provide concrete
information about the Hopf structure of our quantum formal groups.  This
will be especially important to define integer forms and speciale them
at roots of 1.  To be short, we set
  $$  F^\otimes_i := F_i \otimes 1 \, ,  \quad  1^\otimes := 1 \otimes 1 \, ,
\quad  E^\otimes_i := 1 \otimes E_i \, ,  \quad  L^\otimes_\mu := L_{-\mu}
\otimes L_\mu  \qquad \forall\; i \in I, \, \mu \in M  $$
                                                 \par
   The counit  $ \, \epsilon \colon \, {D_{\scriptscriptstyle
M'}}^{\! *} \rightarrow \Cq \, $  is  $ \, \epsilon := 1^* \, $,
hence  $ \, \epsilon (x^*) := \big\langle x^*, 1 \big\rangle \;
\forall \, x^* \in {D_{\scriptscriptstyle M'}}^{\! *} \, $;  thus
  $$  \epsilon \big( F_i^\otimes \big) = 0 \, ,  \qquad  \epsilon
\big( L^\otimes_\mu \big) = 1 \, ,  \qquad  \epsilon \big( E^\otimes_i \big) =
0 \; ;   \eqno (4.4)  $$
the elements above generate the algebra  $ \, j_{\scriptscriptstyle M} \big(
{\uqMg {M'\!} }^* \big) \, $  (in topological sense, cf.~Proposition 4.6), hence
(4.4) uniquely determines  $ \, \epsilon \colon \, j_{\scriptscriptstyle M}
\big( {\uqMg {M'\!} }^* \big)  \longrightarrow \Cq \, $.
                                                    \par
   The antipode of  $ {D_{\scriptscriptstyle M'}}^{\! *} \, $  is by definition
the dual of the antipode of  $ D_{\scriptscriptstyle M'} $,  hence it is
characterized by  $ \, \big\langle S(x^*), x \big\rangle = \big\langle x^*,
S(x) \big\rangle \, $,  for all  $ \, x^* \in {D_{\scriptscriptstyle M'}}^{\!
*} $,  $ x \in D_{\scriptscriptstyle M'} $.  Now consider  $ \, {F_i^\otimes}^f
= F_i^f \otimes 1 \in \uMbm{M} \otimes {\uMbp {P_\infty} } \leq {D_{\scriptscriptstyle M'}}^{\! *} \, $,  $ \, f \in \N \, $:  it is
homogeneous of degree  $ \, - f \alpha_i \, $,  whence  $ S \left( F_i^f \otimes
1 \right) $  has the same degree.  Thus writing  $ S \left( {F_i^\otimes}^f
\otimes 1 \right) $  as a series
  $$  S \left( {F_i^\otimes}^f \right) = \sum_\sigma F_\sigma \cdot
\varPhi_\sigma \cdot E_\sigma  $$
we have  $ \, deg \left( F_\sigma \cdot \varPhi_\sigma \cdot E_\sigma \right) :=
deg(F_\sigma) + deg(E_\sigma) = - f \alpha_i \, $.  Now, the pseudograding of
$ {D_{\scriptscriptstyle M'}}^{\! *} $  induces a pseudograding of
$ {\frak I}^{\scriptscriptstyle M} $  too; hence, since
$ {\frak I}^{\scriptscriptstyle M} $  is a formal Hopf subalgebra of
$ {D_{\scriptscriptstyle M'}}^{\! *} $  (Proposition 4.9), we can apply the same
procedure and get
  $$  S \left( {\fbar_i^\otimes}^f \right) = \sum_\sigma \calF_\sigma \cdot
\varphi_\sigma \cdot \calE_\sigma   \eqno (4.5)  $$
where  $ \, \varphi_\sigma \in {\gerUzM {{M'} } }^{\! *} \, $  and the
$ \calF_\sigma $'s,  resp.~$ \calE_\sigma $'s, are PBW monomials of
$ \calUm $,  resp.~$ \calUp $,  such that  $ \, deg \! \left( \calF_\sigma
\right) + deg \! \left( \calE_\sigma \right) = - f \alpha_i \, $.  An entirely
similar argument yields
  $$  S \left( {F_i^\otimes}^{(f)} \otimes 1 \right) = \sum_\sigma \gerF_\sigma
\cdot \phi_\sigma \cdot \gerE_\sigma   \eqno (4.6)  $$
where  $ \, \phi_\sigma \in {\calUzM {M'} }^* \, $  and the
$ \gerF_\sigma $'s,  resp.~$ \gerE_\sigma $,  are PBW monomials of  $ \gerUm $,
resp.~$ \gerUp $,  such that  $ \, deg \! \left( \gerF_\sigma \right) + deg \!
\left( \gerE_\sigma \right) = - f \alpha_i \, $.
                                      \par
   Now,  $ {\frak I}^{\scriptscriptstyle M} $  and
$ {\Cal I}^{\scriptscriptstyle M} $  can be compared through the natural
embedding  $ \, {\Cal I}^{\scriptscriptstyle M} \cong {\gerUMg {M'} }^*
\hookrightarrow {\calUMg {M'} }^* \cong {\frak I}^{\scriptscriptstyle M} \, $
(dual of  $ \, {\calUMg {M'} } \hookrightarrow {\gerUMg {M'} } \, $);  directly
from definitions we get
  $$  \eqalignno{
   \fbar_\gamma^{\,f} = \prod_{s=1}^f \left( q_\gamma^{\,s} - q_\gamma^{\,-s}
\right) \cdot F_\gamma^{\,(f)} \, ,  &  \quad \;  \ebar_\gamma^{\,e} =
\prod_{s=1}^e \left( q_\gamma^{\,s} - q_\gamma^{\,-s} \right) \cdot
E_\gamma^{\,(e)}  \hskip83pt   &   \forall \; \gamma \in \phipre  \;  \cr
   \fbar_{(r \delta,i)}^{\,f} = {\left( q_i^r - q_i^{-r} \right)}^f {\, f! \,
\over \, r^f \,} \cdot F_{(r \delta,i)}^{\,(f)} \, ,  &  \quad \;  \ebar_{(r
\delta,i)}^{\,e} = {\left( q_i^r - q_i^{-r} \right)}^e {\, e! \, \over \, r^e
\,} \cdot E_{(r \delta,i)}^{\,(e)}  \hskip83pt   &   \forall \; (r \delta,i) \in
\phitildepim  \;  \cr }  $$
thus comparing (4.5) and (4.6) we find
  $$  \eqalign{
   {\ }  \gerF_\sigma \cdot \phi_\sigma \cdot \gerE_\sigma \in \Bigg(
\prod_{u=1}^f  &  {\big( q_i^u - q_i^{-u} \big)}^{-1} \cdot \prod_{\beta, \gamma
\in \phipre} \prod_{r=1}^{f_\beta} \big( q_\beta^r - q_\beta^{-r} \big) \cdot
\prod_{s=1}^{e_\gamma} \big( q_\gamma^s - q_\gamma^{-s} \big) \cdot   \hfill
{\ }  \cr
   &  {} \hskip15pt   \cdot \prod_{(r \delta,i), (t \delta,j) \in
\phitildepim} {\left( q_i^r - q_i^{-r} \right)}^{f_{(r \delta,i)}} {\left( q_j^s
- q_j^{-s} \right)}^{e_{(s \delta,j)}} \Bigg) \cdot {\Cal I}^{\scriptscriptstyle
M}  {\ }  \cr }   \eqno (4.7)  $$
for  $ \, \gerF_\sigma = \prod_{\alpha \in \phitildep} F_\alpha^{\,(f_h)} \, $,
$ \, \gerE_\sigma = \prod_{\alpha \in \phitildep} E_\alpha^{\,(e_k)} \, $.
Similar remarks hold for the other generators of  $ \, {\Cal
A}^{\scriptscriptstyle M} \, $:  in particular, a first consequence is the
following

\vskip5pt

  $ \underline{\hbox{\sl Claim I}} $:  \, {\it  The series  $ S \left(
{F_i^\otimes}^{(f)} \right) $,  $ S \left( \left( L^\otimes_{\mu_i} ; \, c \atop
t \right) \right) $,  $ S \left( L^\otimes_{-\mu_i} \right) $  and  $ S \left(
{E^\otimes_i}^{(e)} \right) $  are  {\sl convergent}  in the  $ \left( q - \qm
\right) {\Cal I}^{\scriptscriptstyle M} $--adic  topology of  $ {\Cal
I}^{\scriptscriptstyle M} \, $;  in particular, they are  {\sl finite}  sums
modulo  $ (q-1) $.}

\vskip5pt

   In principle, one can compute all the terms of these series up to any
fixed order (in  $ \left( q - \qm \right) $):  actually, we need to know them
only up to the zeroth order.  For  $ S \big( F_i^\otimes \big) $  the first
term (call it  $ \hbox{\bf F}_1 $),  of order zero in the  $ \left( q - \qm
\right) $--adic  expansion of  $ S \big( F_i^\otimes \big) $,  corresponds to
the terms  $ \, \gerF_\sigma \cdot \phi_\sigma \cdot \gerE_\sigma \, $  in (4.6)
such that  $ \, \sum_{\alpha \in \phitildep} \big( f_\alpha + e_\alpha \big) = 1
\, $;  on the other hand, these must have degree  $ \, deg \! \left(
\gerF_\sigma \right) + deg \! \left( \gerE_\sigma \right) = -\alpha_i \, $  too,
whence it is  $ \, \gerF_\sigma = F^\otimes_i \, $  and  $ \,
\gerE^\otimes_\sigma = 1 \, $.  Now,  $ \hbox{\bf F}_1 $  takes
non-zero values only on PBW monomials of type  $ \, \ebar_i \cdot
L_\nu \, $,  ($ \nu \in M' \, $);  so let  $ V_{1,i} $  be the
free  $ {\calUzM {M'} } $--module  with basis  $ \left\{ \ebar_i
\right\} $):  direct computation shows that  $ \, \hbox{\bf F}_1 +
q_i^{-2} \cdot F^\otimes_i L^\otimes_{-\alpha_i} \, $  is zero in
$ V_{1,i}^* \, $,  therefore
  $$  \displaylines{
   {\ } \hfill   S \big( F_i^\otimes \big) \equiv - q^{-2} \cdot F_i^\otimes
L^\otimes_{-\alpha_i}   \hfill   \mod \, \left( q-\qm \right) \cdot {\Cal
I}^{\scriptscriptstyle M}  {\ }  \cr }  $$
   \indent   Similar arguments give
  $$  \displaylines{
   {\ } \hfill   S \bigg( \left( L^\otimes_{\mu_i} ; \, 0 \atop 1 \right) \bigg)
\equiv - L^\otimes_{\mu_i} \left( L^\otimes_{\mu_i} ; \, 0 \atop 1 \right)
\hfill   \mod \, \left( q-\qm \right) \cdot {\Cal I}^{\scriptscriptstyle M}
{\ }  \cr
   {\ } \hfill   S \big( E^\otimes_i \big) \equiv - q^{+2} \cdot
L^\otimes_{-\alpha_i} E^\otimes_i   \hfill   \mod \, \left( q - \qm \right)
\cdot {\Cal I}^{\scriptscriptstyle M}  {\ }  \cr }  $$
   \indent   As for the coproduct  $ \, \Delta \colon \, {D_{\scriptscriptstyle
M'}}^{\! *} \rightarrow {D_{\scriptscriptstyle M'}}^{\! *} \otimeshat
{D_{\scriptscriptstyle M'}}^{\! *} \, $,  it is the dual of the product of
$ D_{\scriptscriptstyle M'} $,  hence it is characterized by  $ \, \big\langle
\Delta(x^*), y \otimes z \big\rangle = \big\langle x^*, y \cdot z \big\rangle
\, $.  The same kind of procedure used for  $ S $  may be applied in the present
case.  Thus for instance, if
  $$  \Delta \left( {F_i^\otimes}^{(f)} \right) = \sum_{\sigma,\tau} \big(
\gerF_\sigma \cdot \phi_\sigma \cdot \gerE_\sigma \big) \otimes \big( \gerF_\tau
\cdot \phi_\tau \cdot \gerE_\tau \big)   \eqno (4.8)  $$
is the series expansion of  $ \, \Delta \left( {F_i^\otimes}^{(f)} \right) \, $
as in Proposition 4.9{\it (b)} \, (via  $ \, {\Cal I}^{\scriptscriptstyle M}
\cong {\calUMg {M'} }^* \, $),  then
  $$  \displaylines{
   {\ }  \big( \gerF_\sigma \cdot \phi_\sigma \cdot \gerE_\sigma \big) \otimes
\big( \gerF_{\sigma'} \cdot \phi_{\sigma'} \cdot \gerE_{\sigma'} \big) \in
\Bigg( \prod_{u=1}^f {\big( q_i^u - q_i^{-u} \big)}^{-1} \cdot   \hfill  \cr
   \hfill   \cdot \! \prod_{\beta, \gamma \in \phipre}
\prod_{r=1}^{f_\beta} \! \big( q_\beta^r - q_\beta^{-r} \big) \cdot \!
\prod_{s=1}^{e_\gamma} \! \big( q_\gamma^s - q_\gamma^{-s} \big) \cdot
\hskip-15pt \prod_{(r \delta,i), (t \delta,j) \in \phitildepim}
\hskip-21pt {\left( q_i^r - q_i^{-r} \right)}^{f_{(r \delta,i)}}
{\left( q_j^s - q_j^{-s} \right)}^{e_{(s \delta,j)}} \cdot
\hfill  (4.9)  \cr
   \hfill   \cdot \prod_{\beta, \gamma \in \phipre}
\prod_{r=1}^{f'_\beta} \! \big( q_\beta^r - q_\beta^{-r} \big)
\cdot \prod_{s=1}^{e'_\gamma} \! \big( q_\gamma^s - q_\gamma^{-s}
\big) \cdot \hskip-11pt \prod_{(r \delta,i), (t \delta,j)
\in \phitildepim} \hskip-11pt {\left( q_i^r - q_i^{-r}
\right)}^{f'_{(r \delta,i)}} {\left( q_j^s - q_j^{-s}
\right)}^{e'_{(s \delta,j)}} \Bigg) \cdot
{\Cal I}^{\scriptscriptstyle M}  {\ }  \cr }  $$
for  $ \, \gerF_\sigma = \prod_{\alpha \in \phitildep}
F_\alpha^{\,(f_h)} $,  $ \, \gerE_\sigma =
\prod_{\alpha \in \phitildep} E_\alpha^{\,(e_k)} $,
$ \, \gerF_{\sigma'} = \prod_{\alpha \in \phitildep}
F_\alpha^{\,(f'_h)} $,  $ \, \gerE_{\sigma'} = \prod_{\alpha
\in \phitildep} E_\alpha^{\,(e'_k)} \! $.  Similar remarks hold for
the other generators of  $ \, {\Cal A}^{\scriptscriptstyle M} \, $.
As a first consequence, we have:

\vskip5pt

  $ \underline{\hbox{\sl Claim II}} $:  \, {\it  The series  $ \Delta \left(
{F^\otimes_i}^{(f)} \right) $,  $ \Delta \left( \left( L^\otimes_{\mu_i} ; \,
c \atop t \right) \right) $,  $ \Delta \left( L^\otimes_{-\mu_i} \right) $  and
$ \Delta \left( {E^\otimes_i}^{(e)} \right) $  are  {\sl convergent}  in the
$ \left( q - \qm \right) \cdot \left( {\Cal I}^{\scriptscriptstyle M}
\otimeshat {\Cal I}^{\scriptscriptstyle M} \right) $--adic  topology of  $ {\Cal
I}^{\scriptscriptstyle M} \otimeshat {\Cal I}^{\scriptscriptstyle M} \, $;  in
particular, they are  {\sl finite}  sums modulo  $ (q-1) $.}

\vskip5pt

   Direct computation gives us the following congruences modulo  $ {\left( q -
\qm \right)}^2  \, $:
  $$  \displaylines{
   {\ }  \Delta \big( F^\otimes_i \big) \equiv F^\otimes_i
\otimes 1^\otimes + 1^\otimes \otimes F^\otimes_i + (q_i - 1) \cdot
\bigg( {L^\otimes_{\alpha_i} ; \, 0 \atop 1} \bigg) \otimes F^\otimes_i +
\hfill  \cr
   \hfill   + \, {\left( q_i - q_i^{-1} \right)}^{-1} \cdot \!\!\!\!\!\!\!
                  \sum_{\Sb  \alpha, \beta \in \phitildep  \\
                             p(\alpha) - p(\beta) = -\alpha_i  \\  \endSb}
\hskip-15pt  C^{i,+}_{\alpha,\beta} \, \big( q_\alpha - q_\alpha^{\,-1} \big)
\big( q_\beta - q_\beta^{\,-1} \big) \cdot L^\otimes_{\alpha_i}
E^\otimes_\alpha \otimes F^\otimes_\beta   \hfill   \mod \, {\left( q - \qm
\right)}^2  \cr
   \Delta \bigg( \bigg( {L^\otimes_{\mu_i} ; \, 0 \atop 1} \bigg) \bigg) \equiv
\left( L^\otimes_{\mu_i} ; \, 0 \atop 1 \right) \otimes 1^\otimes + 1^\otimes
\otimes \left( L^\otimes_{\mu_i} ; \, 0 \atop 1 \right) + (q_i - 1) \cdot \left(
L^\otimes_{\mu_i} ; \, 0 \atop 1 \right) \otimes \left( L^\otimes_{\mu_i} ; \, 0
\atop 1 \right) +   \hfill  \cr
   + \; {(2)}_{\qm}^{\,2} {(d_i)}_q^{\,-1} \cdot \!\!\!\! \sum_{\Sb
\alpha, \beta \in \phitildep  \\
                p(\alpha) - p(\beta) = 0  \\   \endSb}
\!\! (q - 1) \, C_{\alpha,\beta} \, {[d_\gamma]}_q {\big[ (\mu_i \vert \gamma)
\big]}_q \cdot L^\otimes_{\mu_i} E^\otimes_\alpha \otimes F^\otimes_\beta
L^\otimes_{\mu_i}   \hfill   \mod \, {\left( q - \qm \right)}^2  \cr }  $$
  $$  \displaylines{
   {\ }  \Delta \big( E^\otimes_i \big) \equiv 1^\otimes \otimes E^\otimes_i +
E^\otimes_i \otimes 1^\otimes + (q_i - 1) \cdot E^\otimes_i \otimes \bigg(
{L^\otimes_{\alpha_i} ; \, 0 \atop 1} \bigg) -   \hfill  \cr
   \hfill   - \, {\left( q_i - q_i^{-1} \right)}^{-1} \cdot \!\!\!\!\!\!\!
         \sum_{\Sb  \alpha, \beta \in \phitildep  \\
                    p(\alpha) - p(\beta) = +\alpha_i  \\  \endSb}
\hskip-15pt   C^{i,-}_{\alpha,\beta} \, \big( q_\alpha - q_\alpha^{\,-1} \big)
\big( q_\beta - q_\beta^{\,-1} \big) \cdot E^\otimes_\alpha \otimes
F^\otimes_\beta L^\otimes_{\alpha_i}   \hfill   \mod \, {\left( q - \qm
\right)}^2  \cr }  $$
where the  $ C^{i,\pm}_{\alpha,\beta} $'s  are given by the equations
$ \; \pi_i^- \big( [ \dot{F}_\alpha, \dot{E}_\beta ] \big) =
C^{i,-}_{\alpha,\beta} \cdot F_i \, $,  $ \, \pi_i^+ \big( [ \dot{F}_\alpha,
\dot{E}_\beta ] \big) = C^{i,+}_{\alpha,\beta} \cdot E_i \, $  ($ \, \pi_i^- :
\uqMg{Q} \twoheadrightarrow \Cq \cdot F_i \, $  and  $ \, \pi_i^+ : \uqMg{Q}
\twoheadrightarrow \Cq \cdot E_i \, $  being the canonical maps)  and the
$ C_{\alpha,\beta} $'s  by the equation  $ \, \big[ \dot{F}_\alpha,
\dot{E}_\beta \big] = C_{\alpha,\beta} \, {\, L_{p(\alpha)} -
L_{p(\alpha)}^{\,-1} \, \over \, q_\beta - q_\beta^{\,-1} \,} \, $:  here  $ \,
\dot{F}_\gamma \, $,  resp.~$ \, \dot{E}_\gamma \, $,  is the dual, resp.~minus
the dual, of  $ E_\gamma $,  resp.~$ F_\gamma $,  with respect to  $ \pi $,
resp.~$ \overline{\pi} $.

\vskip1,7truecm

 \centerline{ \bf  \S \; 5 \,  The quantum group  $ \uqMh{M} $ }

\vskip10pt

   {\bf 5.1 The quantum enveloping algebra  $ \, \uqMh{M} \, $.}  \  The
results of \S 4 can be given an axiomatic form: to this end, we introduce
a new object  $ {\uqMh M } $  which is with
respect to  $ \uh $  what  $ {\uqMg M } $  is for
$ \ug $.  Here  $ M $  is a fixed lattice as in \S 2.1.
                                                       \par
   We define  $ \hat{\hbox{\bf H}}_{\scriptscriptstyle M} $  to be the
associative  $ \Cq $--algebra  with 1 with generators  $ \; F_i \, $,
$ \, L_\mu \, $,  $ \, E_i \, $  for all  $ \, \lambda \in M , i
\in I \, $,  \; and relations (for all  $ \mu, \nu \in M $,
$ i, j \in I $,  $ i \neq j $)
  $$  \displaylines {
   \hfill   L_0 = 1 \, ,  \; \qquad  L_\mu L_\nu = L_{\mu + \nu} \, ,
\; \quad \qquad  E_i F_j - F_j E_i = 0   \hfill  \phantom{(5.1)}  \cr
   \hfill   L_\mu F_j = q^{(\alpha_j | \mu)} F_j L_\mu \, ,  \qquad
\sum_{k=0}^{1-a_{ij}} (-1)^k {\left[ {1-a_{ij} \atop k}
\right]}_{\! q_i} F_i^{\,1-\aij-k} F_j F_i^{\,k} = 0 \, \phantom{.}
\hfill  (5.1)  \cr
   \hfill   L_\mu E_j = q^{(\alpha_j | \mu)} E_j L_\mu \, ,  \qquad
   \sum_{k=0}^{1-a_{ij}} (-1)^k {\left[ {1-a_{ij} \atop k}
\right]}_{\! q_i} E_i^{\,1-\aij-k} E_j E_i^{\,k} = 0 \, .
\hfill  \phantom{(5.1)}  \cr }  $$
We also use notation  $ \, M_i := L_{\mu_i} $  ($ i \in I_\infty $),
where  $ \big\{\, \mu_i \,\big\vert\, i \in I_\infty \big\} $  is
a fixed  $ \Z $--basis  of  $ M $.
                                                     \par
   Define  $ \, \N^{\phitildep} \, $  to be the set of all functions  $ \, f:
\phitildep \rightarrow \N \, $  such that  $ \, f(\alpha) = 0 \, $  for almost
all  $ \, \alpha \in \phitildep \, $.  For any  $ \, \phi :=
{(f_\alpha)}_{\alpha \in \phitildep} \in \N^{\phitildep} \, $,  $ \, \eta :=
{(e_\alpha)}_{\alpha \in \phitildep} \in \N^{\phitildep} \, $,  set  $ \, F_\phi
:= \prod_{\alpha \in \phitildep} F_\alpha^{\,f_\alpha} \, $,  $ \, E_\eta :=
\prod_{\alpha \in \phitildep} E_\alpha^{\,e_\alpha} \, $,  where the products
are meant to be ordered, like in \S 2.  We shall also use a similar notation
when dealing with PBW monomials of  $ \gerU_{\scriptscriptstyle \pm} $  or of
$ \calU_{\scriptscriptstyle \pm} $,  e.g.~$ \, \gerF_\phi := \prod_{\alpha \in
\phitildep} F_\alpha^{\,(f_\alpha)} $,  $ \, \calE_\eta := \prod_{\alpha \in
\phitildep} \ebar_\alpha^{\,e_\alpha} $.  Moreover, for any  $ \, \tau \in
\N^{I_\infty} = \N^{n+2} \, $  set  $ \, B_\tau := u_\tau^{\,*} \, $  (the
element of the pseudobasis  $ {\Cal B}_{\scriptscriptstyle M} $  constructed in
Lemma 4.4{\it (b)}.
                                                     \par
   {\it  We define  $ {\uqMh M } $  to be the completion of
$ \hat{\hbox{\bf H}}_{\scriptscriptstyle M} $  by means of formal series
(i.e.~infinite linear combinations), with coefficients in  $ \Cq $,  in the
elements of the set}
  $$  {\Bbb B}_{\scriptscriptstyle M} := \left\{\, F_\phi \cdot B_\tau \cdot
E_\eta \;\Big\vert\; \phi, \eta \in \N^{\phitildep} \, , \; \tau \in
\N^{I_\infty} \,\right\} \, .  $$
   \indent   Thus  $ {\uqMh M } $  is the completion of  $ \hat{\hbox{\bf
H}}_{\scriptscriptstyle M} $  with respect to the topology (of
$ \hat{\hbox{\bf H}}_{\scriptscriptstyle M} $)  for which a fundamental system
of neighborhoods of  $ 0 $  is the set of vector subspaces of $ \hat{\hbox{\bf
H}}_{\scriptscriptstyle M} $  which contain almost all the elements of  $ {\Bbb
B}_{\scriptscriptstyle M} \, $,  and the set $ \, {\Bbb B}_{\scriptscriptstyle
M} \, $  is a  {\it pseudobasis}  of  $ {\uqMh M } $.  Roughly speaking,
$ \uqMh{M} $  is an algebra of (non-commutative) formal series with the (5.1)
as commutation rules.  Finally, thanks to Lemma 4.4, we can identify
$ \uqMh{M} $  with the space of formal series in the  $ F_\alpha $'s,
$ E_\alpha $'s,  ($ \alpha \in \phitildep $)  with coefficients in  $ \, {\uzM
{M'} }^* \, $.
                                                    \par
   From \S 4 we can explicitely realize  $ {\uqMh M } $  and endow it
with a Hopf structure: in fact, the definition of  $ \uqMh{M} $  is
just a presentation of  $ {\uqMg{M'} }^* $,  as the following shows:

\vskip7pt

\proclaim{Theorem 5.2}  There exists an isomorphism of topological
$ \Cq $--algebras
  $$  \nu_{\scriptscriptstyle M} \colon \, {\uqMh M } \,
{\buildrel \cong \over \llongrightarrow}
\, j_{\scriptscriptstyle M} \left( {\uqMg {M'} }^* \right)  $$
  $$  F_i \mapsto F^\otimes_i := F_i \otimes 1 \, ,  \quad  \; L_\mu \mapsto
L^\otimes_\mu := L_{-\mu} \otimes L_\mu \, ,  \quad  \; E_i \mapsto E^\otimes_i
:= 1 \otimes E_i \, .   \leqno \text{given by}  $$
   \indent  Then the pull-back of the formal Hopf structure of  $ \,
j_{\scriptscriptstyle M} \left( {\uqMg {M'} }^* \right) $  defines a formal Hopf
structure on  $ {\uqMh M } $,  so that  $ \nu_{\scriptscriptstyle M} $  and
$ {j_{\scriptscriptstyle M}}^{\! -1} \smallcirc \nu_{\scriptscriptstyle M} \, $
are isomorphisms of formal Hopf algebras.
\endproclaim

\demo{Proof}  By construction  $ \, \hat{\hbox{\bf H}}_{\scriptscriptstyle
M} \cong \um \otimes {\uzM M } \otimes \up \cong A^{\scriptscriptstyle M} \;
(\, \subseteq j_{\scriptscriptstyle M} \left( {\uqMg {M'} }^* \right) \,) \, $
as vector spaces;  now  $ \, F_i \otimes 1, \, L_{-\mu} \otimes L_\mu, \, 1
\otimes E_i \in \uMbm{M} \otimes {\uMbp {P_\infty} } \, $,  hence comparing
(5.1) and (2.1) we see that formulas above gives a well-defined isomorphism of
algebras  $ \, \nu_{\scriptscriptstyle M} \colon \, \hat{\hbox{\bf
H}}_{\scriptscriptstyle M} @>{\cong}>> A^{\scriptscriptstyle M} \, $.
Moreover,  $ A^{\scriptscriptstyle M} $  contains a pseudobasis  $ \hbox{\bf
B}_{\scriptscriptstyle M} $  of  $ j_{\scriptscriptstyle M} \left( {\uqMg {M'}
}^* \right) $  (cf.~Lemma 4.4, Proposition 4.6, and Remark 4.7)  such that  $ \,
\nu_{\scriptscriptstyle M}({\Bbb B}_{\scriptscriptstyle M}) = \hbox{\bf
B}_{\scriptscriptstyle M} \, $,  hence  $ \nu_{\scriptscriptstyle M} $
continuosly extends, in a unique way, to an isomorphism of topological algebras
$ \, \nu_{\scriptscriptstyle M} \colon \, \uqMh{M} @>\cong>>
j_{\scriptscriptstyle M} \left( {\uqMg {M'} }^* \right) \, $,  q.e.d.
$ \square $
\enddemo

\vskip7pt

   $ \underline{\hbox{\sl Remark}} $:  \, For "restricted" specializations at
roots of 1 the algebra  $ \uqMh{M} $  is too big; in order to get a reasonable
$ R $--integer  form, it turns out to be necessary to impose two bounding
conditions.  The first condition is on the size of the "toral part" of an
element in  $ \uqMh{M} $:  we tackle it in \S\S 5.4--5 below; the second one is
on the behaviour of an element in  $ \uqMh{M} $  (as a series), and is inspired
by (4.7) and (4.9): we deal with it in \S\S 5.7--9.
                                                \par
   Hereafter, we freely use the term  {\sl pseudobasis}  to mean a topological
basis of a topological module, so that any element in the module has a unique
expansion as a series in the elements of the pseudobasis.

\vskip7pt

\proclaim{Definition 5.3}  We define the subset  $ \, \Omega_{\scriptscriptstyle
M} \, $  of  $ \uqMh{M} $  to be
  $$  \Omega_{\scriptscriptstyle M} := \Big\{\, x = \sum\nolimits_\sigma
F_\sigma \cdot \varPhi_\sigma \cdot E_\sigma \in {\uqMh M } \;\Big\vert\;
\varPhi_\sigma \in {\uzM M }, \; \forall\, \sigma \,\Big\}  $$
(where  $ \, x = \sum_\sigma F_\sigma \cdot \varPhi_\sigma \cdot E_\sigma \, $
is the expansion of  $ x $  as a series with coefficients in
$ \, {\uzM {M'} }^* \, $).
\endproclaim

\vskip7pt

\proclaim{Lemma 5.4}                          \hfill\break
   \indent   (a)  The set of all products  $ \, F \cdot L \cdot E \, $  in which
$ F $,  resp.~$ L $,  resp.~$ E $,  is any element of a fixed PBW basis of
$ \um $,  resp.~$ \uzM{M} $,  resp.~$ \up $,  is a pseudobasis of  $ \,
\Omega_{\scriptscriptstyle M} \, $.
                                              \hfill\break
   \indent   (b)  $ \, \Omega_{\scriptscriptstyle M} \, $  is a formal Hopf
subalgebra of  $ {\uqMh M } $.
\endproclaim

\demo{Proof}  Claim  {\it (a)}  follows from definitions.  As for
{\it (b)},  it is clear that  $ \Omega_{\scriptscriptstyle M} $
is a subalgebra of  $ \uqMh{M} $.  Now we show that it is also
closed for the antipode and the coproduct.
                                                   \par
   Let  $ \, x = \sum_\tau F'_\tau \cdot \varPhi'_\tau \cdot E'_\tau
\in \Omega_{\scriptscriptstyle M} \, $:  then $ \, \varPhi'_\tau =
\sum_{\mu \in M } c_{\tau,\mu} L_\mu \, $  with  $ \, c_{\tau,\mu}
\neq 0 \, $  for finitely many  $ \mu $.  Let  $ \, S(x) =
\sum_\sigma F_\sigma \cdot \varPhi_\sigma \cdot E_\sigma \, $:
for any fixed  $ \bar\sigma $,  we must prove that  $ \,
\varPhi_{\bar\sigma}
                           \in {\uzM M } $\break
$ \left(\, \subseteq {\uzM {M'} }^* \,\right) \, $,  so that
$ \, S \big( \Omega_{\scriptscriptstyle M} \big) =
\Omega_{\scriptscriptstyle M} \, $;  to this end, we identify
$ \, \uqMh{M} \cong {\uqMg {M'} }^* \, $  (cf.~Theorem 5.2).
By Lemma 2.5 there exists two PBW monomials  $ \,
\dot{\!\calE}_{\bar\sigma} $  and  $ \,\dot{\!\calF}_{\bar\sigma} $
such that
  $$  \left\langle S(x) \, , \, \,\dot{\!\calE}_{\bar\sigma} \cdot y \cdot
\,\dot{\!\calF}_{\bar\sigma} \right\rangle = \left\langle F_\sigma \cdot
\varPhi_\sigma \cdot E_\sigma \, , \, \,\dot{\!\calE}_{\bar\sigma} \cdot y \cdot
\,\dot{\!\calF}_{\bar\sigma} \right\rangle = \big\langle F_{\bar\sigma} \ , \,
\,\dot{\!\calE}_{\bar\sigma} \big\rangle \cdot \big\langle E_{\bar\sigma} \, ,
\, \,\dot{\!\calF}_{\bar\sigma} \big\rangle \cdot \varPhi_{\bar\sigma} \big( y
\cdot L_\alpha \big)  $$
for all  $ \, y \in {\uzM {M'} } \, $,  with  $ \, \alpha := s \big(
F_{\bar\sigma} \big) + s \big( E_{\bar\sigma} \big) \, $  and  $ \,
c_{\bar\sigma} := \big\langle F_{\bar\sigma} \, , \,
\,\dot{\!\calE}_{\bar\sigma} \big\rangle \cdot \big\langle E_{\bar\sigma} \, ,
\, \,\dot{\!\calF}_{\bar\sigma} \big\rangle \neq 0 \, $:  in other words,
$ \, \varPhi_{\bar\sigma} = {c_{\bar\sigma}}^{\! -1} \cdot \Big( \big(
L_{-\alpha} \cdot \,\dot{\!\calF}_{\bar\sigma} \big) \triangleright S(x)
\triangleleft \,\dot{\!\calE}_{\bar\sigma} \Big){\Big\vert}_{{\uzM {M'} }} \, $
(where  $ \triangleleft $  and  $ \triangleright $  denote standard left and
right action, cf.~[DL], \S 1.4),  hence we have to study  $ \, \left\langle S(x)
\, , \, \,\dot{\!\calE}_{\bar\sigma} \cdot y \cdot L_{-\alpha}
\,\dot{\!\calF}_{\bar\sigma} \right\rangle \, $  as a function of  $ \, y \in
{\uzM {M'} } \, $;  by linearity we can assume  $ \, y = L_\nu \, $,  $ \, \nu
\in M' \, $.  By definition,  $ \, \left\langle S(x) \, , \,
\,\dot{\!\calE}_{\bar\sigma} \cdot y \cdot L_{-\alpha}
\,\dot{\!\calF}_{\bar\sigma} \right\rangle = \Big\langle x \, ,
\, S \big( \,\dot{\!\calE}_{\bar\sigma} \cdot y L_{-\alpha} \cdot
\,\dot{\!\calF}_{\bar\sigma} \big) \Big\rangle \, $;  to compute the
latter we have to "straighten"  $ S \big( \,\dot{\!\calE}_{\bar\sigma}
\cdot y \cdot L_{-\alpha} \,\dot{\!\calF}_{\bar\sigma} \big) $,
i.e.~to write it in terms of a PBW basis of  $ \, \um \otimes
{\uzM {M'} } \otimes \up \, $.
                                                  \par
  Since  $ \, S \big( \,\dot{\!\calE}_{\bar\sigma} \cdot y \cdot
L_{-\alpha} \,\dot{\!\calF}_{\bar\sigma} \big) = S \big( L_{-\alpha}
\,\dot{\!\calF}_{\bar\sigma} \big) \cdot S(y) \cdot S \big(
\,\dot{\!\calE}_{\bar\sigma} \big) \, $,  let us consider the
various factors.  First,  $ S \big( L_{-\alpha}
\,\dot{\!\calF}_{\bar\sigma} \big) \in {\uMbm {M'} } \, $,
and  $ S \big( L_{-\alpha} \,\dot{\!\calF}_{\bar\sigma} \big) $
{\it does not depend on}  $ y $.  Second,  $ S(y) = S \big( L_\nu
\big) = L_{-\nu} \, $.  Third,  $ S \big( \,\dot{\!\calE}_{\bar\sigma}
\big) \in {\uMbp {M'} } \, $,  and  $ S \big(
\,\dot{\!\calE}_{\bar\sigma} \big) $  {\it does not depend on}  $ y $.
                                                  \par
   Now we straighten the product.  Commuting  $ S \big( L_{-\alpha}
\,\dot{\!\calF}_{\bar\sigma} \big) $  and  $ \, S(y) = L_{-\nu} \, $
produces a coefficient  $ \, q^{-(\nu \vert \beta_{\bar\sigma})} = {\big\langle
L_{-\beta_{\bar\sigma}}, L_\nu \big\rangle}_{\! \overline{\pi}} $  in which
$ \, \beta_{\bar\sigma} \in Q_- \, $  is the weight of  $ S \big(
\,\dot{\!\calF}_{\bar\sigma} \big) $.  Straightening the product  $ \, S \big(
L_{-\alpha} \,\dot{\!\calF}_{\bar\sigma} \big) \cdot S \big(
\,\dot{\!\calE}_{\bar\sigma} \big) \, $  produces a sum  $ \, \sum_k x_k \, $
of terms which  {\it do not depend on}  $ y \, $.  Straightening the product
$ \, S(y) = L_{-\nu} \cdot \sum_k x_k \, $  produces for each term  $ x_k $  a
coefficient  $ \, q^{-(\nu \vert \gamma_{ \bar\sigma, k})} = {\big\langle
L_{-\gamma_{ \bar\sigma, k}}, L_\nu \big\rangle}_{\! \overline{\pi}} \, $,
where  $ \, \gamma_{\bar\sigma, k} \in Q_+ \, $  is the weight of the "positive"
part  $ x_k^+ $  of  $ x_k $  (with respect to the triangular decomposition).
                                                  \par
   Therefore  $ \, \left\langle x \, , \, S \big( \,\dot{\!\calE}_{\bar\sigma}
\cdot y \cdot L_{-\alpha} \,\dot{\!\calF}_{\bar\sigma} \big) \right\rangle \, $
depends on  $ y $  according to the functions  $ L_{-\beta_{ \bar\sigma}} $,
$ L_{-\gamma_{ \bar\sigma, k}} $,  and  $ \, \varPhi'_\tau \smallcirc S \, $:
to be precise,  $ \, \varPhi_{\bar\sigma} = \big( L_{-\alpha}
\,\dot{\!\calF}_{\bar\sigma} \triangleright S(x) \triangleleft
\,\dot{\!\calE}_{\bar\sigma} \big) \Big\vert_{{\uzM {M'} }} \, $  is a
linear combination of functions of type  $ \, L_{-\beta_{\bar\sigma}} \cdot
\big( \varPhi'_\tau \smallcirc S \big) \cdot L_{-\gamma_{\bar\sigma, k}} =
\sum_{\mu \in M} c_{\tau,\mu} L_{- \mu - \beta_{\bar\sigma} -
\gamma_{\bar\sigma, k}} \, $,  so  $ \, \varPhi_{\bar\sigma} \in {\uzM M }
\, $,  q.e.d.  An entirely analogous procedure   --- slightly simpler indeed
---   works for comultiplication, thus proving that  $ \, \Delta \big(
\Omega_{\scriptscriptstyle M} \big) \subseteq \Omega_{\scriptscriptstyle M}
\otimeshat \Omega_{\scriptscriptstyle M} \, $.  The claim follows.
$ \square $
\enddemo
%
%
%
%

\vskip7pt

   Now we start introducing integer forms of  $ \uqMh{M} $  and proving their
first properties.

\vskip7pt

\proclaim{Definition 5.5}  We define  $ \hat{\Cal H}_{\scriptscriptstyle M} $
to be the  $ R $--subalgebra  of  $ {\uqMh M } $
                                        generated by\break
$ \, \big\{\, \fbar_\alpha, \, L_\mu, \, \ebar_\alpha \,\big\vert\, \alpha \in
\phitildep \,;  \; \mu \in M \,\big\} \, $,  and  $ \, \calUMh{M} \, $  to be
its closure in  $ \uqMh{M} $.
\endproclaim

\vskip7pt

\proclaim{Theorem 5.6}  $ \calUMh{M} $  is an  $ R $--integer  form (in
topological sense) of  $ \uqMh{M} $,  as a formal Hopf algebra, with
$ R $--pseudobasis
  $$  {\widetilde{\Bbb B}}_{\scriptscriptstyle M} :=
\Big\{\, Y_{\scriptscriptstyle \phi, \tau, \eta} \,\Big\vert\, \tau \in
\N^{I_\infty} \, ;  \, \phi, \eta \in \N^{\phitildep} \,\Big\} = \Big\{\,
\calF_\phi \cdot B_\tau \cdot \calE_\eta \,\Big\vert\, \tau \in \N^{I_\infty}
\, ; \, \phi, \eta \in \N^{\phitildep} \Big\}   \eqno (5.2)  $$
where  $ \, Y_{\scriptscriptstyle \phi, \tau, \eta} := \calF_\phi \cdot B_\tau
\cdot \calE_\eta \, $;  in particular  $ \, \nu_{\scriptscriptstyle M}
\Big( {\calUMh M } \Big) = j_{\scriptscriptstyle M} \left( {\gerUMg {M'\!} }^*
\right) =: {\frak I}^{\scriptscriptstyle M} \, $.
\endproclaim

\demo{Proof}  By construction  $ \, {\widetilde{\Bbb B}}_{\scriptscriptstyle M}
\subseteq \calUMh{M} \, $,  so the claim follows directly from \S 5.1.
$ \square $
\enddemo

\vskip7pt

\proclaim{Lemma 5.7}  Let  $ \, \Omega^\vee_{\scriptscriptstyle M} :=
\Omega_{\scriptscriptstyle M} \cap \nu_{\scriptscriptstyle M}^{\,-1}
({\Cal I}^{\scriptscriptstyle M}) \, $;  then
$ \, \Omega^\vee_{\scriptscriptstyle M} \, $  is the set of all  $ \, x \in
\uqMh{M} \, $  which have a series expansion  $ \, x = \sum_{\gerF, \gerE} \gerF \cdot \phi_{\gerF, \gerE} \cdot \gerE \, $  where the  $ \, \gerF \in \gerUm
\, $  and  $ \, \gerE \in \gerUp \, $  are PBW monomials and  $ \, \phi_{\gerF,
\gerE} \in \gerUzM{M} \, $  for all  $ \gerF $,  $ \gerE $.
\endproclaim

\demo{Proof}  Trivial from definitions and Proposition 4.9{\it (b)\/}.
$ \square $
\enddemo

\vskip7pt

\proclaim{Definition 5.8}
                                    \hfill\break
   \indent  (a) We define  $ \, \hat{\frak H}_{\scriptscriptstyle M} \, $
to be the  $ R $--subalgebra  of  $ \uqMh{M} $
                               generated by the set\break
$ \, \Big\{\, F_i^{\,(f)}, \, \left( M_j ; \, c \atop t \right), \, M_j^{-1},
\, E_i^{\,(e)} \;\Big\vert\; f, c, t, e \in \N \, ; \; j \in I_\infty \,\Big\}
\, $.
                                    \hfill\break
   \indent   (b) We define  $ \, \gerUMh{M} \, $  to be the subset of all
elements  $ x $  in  $ \Omega^\vee_{\scriptscriptstyle M} $  whose series
expansion  $ \; x = \sum_{\varphi, \eta \in \N^{\phitildep}} \gerF_\varphi \cdot
\phi_{\varphi,\eta} \cdot \gerE_\eta \; $  (with notation of \S 5.1) is such
that, for some  $ \, p(x) \in \Z \left[ q, \qm \right] \, $,
  $$  \eqalign{
   {\ }  \phi_{\varphi,\eta} \in {p(x)}^{-1} \cdot \Bigg( \prod_{\beta, \gamma
\in \phipre} \prod_{r=1}^{f_\beta}  &  \big( q_\beta^r - q_\beta^{-r} \big)
\cdot \prod_{s=1}^{e_\gamma} \big( q_\gamma^s - q_\gamma^{-s} \big) \cdot
\hfill  {\ }  \cr
   \cdot \prod_{(r \delta,i), (s \delta,j) \in \phitildepim}  &
\!\!\! {\left( q_i^r - q_i^{-r} \right)}^{f_{(r \delta,i)}} {\left( q_j^s -
q_j^{-s} \right)}^{e_{(s \delta,j)}} \Bigg) \cdot \gerUzM{M}  {\ }  \cr }
\eqno (5.3)  $$
for almost all the  $ \, \phi = {(f_\alpha)}_{\alpha \in \phitildep}, \eta =
{(e_\alpha)}_{\alpha \in \phitildep} \in \N^{\phitildep} $.
\endproclaim

\vskip7pt

  {\sl  $ \underline{\text{Remark}} $:}  \, indeed, what we do in
Definition 5.8 is a "refinement" of the construction provided in
[Dr], \S 7, to locate inside a quantized formal series Hopf algebra
--- whose semiclassical limit is a formal Poisson group
$ G^\infty $  ---   a quantized universal enveloping algebra
--- whose semiclassical limit yields the Lie bialgebra dual to
that of  $ G^\infty \, $;  see also [Ga3].  This will permit us
to specialize at roots of 1, which is not possible in [Dr].

\vskip7pt

\proclaim{Theorem 5.9} \; (a) \,  $ \gerUMh{M} $  is a topological Hopf
subalgebra (over  $ R \, $)  of  $ \uqMh{M} $.
                                    \hfill\break
   \indent   (b) \,  $ \gerUMh{M} $  is an  $ R $--integer  form of
$ \uqMh{M} $  and  $ \Omega_{\scriptscriptstyle M} $  (as topological Hopf
algebras).
\endproclaim

\demo{Proof}  By construction  $ \gerUMh{M} $  is
an  $ R $--subalgebra  of  $ \uqMh{M} $  and of
$ \Omega_{\scriptscriptstyle M} $;  moreover,
Theorem 5.2 and Proposition 4.9{\it (b)}  ensure
that  $ \Omega^\vee_{\scriptscriptstyle M} $  is
an  $ R $--integer  form (in topological sense) of
$ \Omega_{\scriptscriptstyle M} $  (as an algebra), hence also
$ \gerUMh{M} $  is.  Furthermore, Proposition 4.9{\it (b)}  and Lemma 5.4 imply
that  $ \Omega^\vee_{\scriptscriptstyle M} $  is a (topological) Hopf subalgebra
of  $ \Omega_{\scriptscriptstyle M} $.  Now the analysis in \S 4.13 via
$ \nu_{\scriptscriptstyle M}^{\,-1} $  gives  $ \, S \big( \gerUMh{M} \big) =
\gerUMh{M} \, $  and  $ \, \Delta \big( \gerUMh{M} \big) \subseteq \gerUMh{M}
\otimeshat \gerUMh{M} \, $,  thus  $ \gerUMh{M} $  is also a topological Hopf
subalgebra.  The claim follows.   $ \square $
\enddemo

\vskip7pt

  {\sl  $ \underline{\text{Remark}} $:}  \, in addition, the analysis
in \S 4.13 shows   --- via  $ \, {\nu_{\scriptscriptstyle M}}^{\! -1}
\, $  ---   also that  $ \gerUMh{M} $  contains the minimal topological
Hopf subalgebra (over  $ R \, $)  of  $ \uqMh{M} $  (or of
$ \Omega_{\scriptscriptstyle M} \, $)  which contains  $ \hat{\frak
H}_{\scriptscriptstyle M} \, $:  this follows from formulas (4.7) and (4.9) and
condition (5.3) (which is modeled on (4.7) indeed).  In fact for our purposes we
might also take this minimal algebra to play the role of  $ \gerUMh{M} $:  the
crux fact in any case is that (5.3) holds.

\vskip7pt

   {\bf 5.10  Presentation of  $ \gerUMh{M} $.}  \;  From the similar result
available for  $ \gerUMg{M} $  (cf.~[DL], \S 3.4: it deals with the finite case,
but it is the like) we get a presentation of  $ \gerUMh{M} $  by (topological)
generators and relations.  The algebra  $ \hat{\frak H}_{\scriptscriptstyle M} $
of \S 5.8 is the associative  $ R $--algebra  with 1 with
generators
  $$  F_i^{\,(s)} \, ,  \quad  M_j \, ,  \quad  M_j^{-1} \, ,  \quad  \left(
M_j \, ; \, c \atop t \right) \, ,  \quad  E_i^{\,(r)}  $$
($ i \in I $;  $ j \in I_\infty $;  $ c \in \Z $,  $ t, r, s \in \N $;  here we
set  $ \, M_j := L_{\mu_j} \, $),  and relations
  $$  \displaylines{
   M_j M_j^{-1} = 1 = M_j^{-1} M_j \, , \; \qquad  M_j^{\,\pm 1} M_j^{\,\pm
1} = M_j^{\,\pm 1} M_j^{\,\pm 1}  \cr }  $$
  $$  \displaylines{
   M_j^{\,\pm 1} \left( M_j \, ; \, c  \atop t \right) =
\left( M_j \, ; \, c \atop t \right) M_j^{\,\pm 1} \, ,  \quad
\left( M_j \, ; \, c \atop 0 \right) = 0 \, , \quad  (q_i - 1)
\left( M_j \, ; \, 0 \atop 1 \right) = M_j - 1  \cr
   \left( M_j \, ; \, c \atop t \right) \left( M_j \, ; \, c-t \atop s \right) =
{\left( t+s \atop t \right)}_{\! q} \left( M_j \, ; \, c \atop t+s \right)  \cr
   \left( M_j \, ; \, c+1 \atop t \right) - q^t \left( M_j \, ; \, c \atop t
\right) = \left( M_j \, ; \, c \atop t-1 \right) \, ,  \; \qquad \forall \, t
\geq 1  \cr
   \left( M_j \, ; \, c \atop t \right) = \sum_{p \geq 0}^{p \leq c, t}
q^{(c-p)(t-p)} \left( c \atop p \right)_{\! q} \left( M_j \, ; \, 0 \atop t-1
\right) \, ,   \; \qquad \forall \, c \geq 0  \cr
   \left( M_j \, ; \, -c \atop t \right) = \sum_{p=0}^t {(-1)}^p q^{-t(c+p) +
p(p+1)/2} {\left( p+c-1 \atop p \right)}_{\! q} \left( M_j \, ; \, 0 \atop t-p
\right) \, ,   \; \qquad \forall \, c \geq 1  \cr
   \left( M_j \, ; \, c+1 \atop t \right) - \left( M_j \, ; \, c \atop t \right)
= q^{c-t+1} M_j \left( M_j \, ; \, c \atop t-1 \right) \, , \; \qquad \forall \,
t \geq 1  \cr
   M_j E_j^{\,(p)} = q^{+p \, (\alpha_j \vert \mu_i)} E_j^{\,(p)} M_j \, , \; \;
 M_j F_j^{\,(p)} = q^{-p \, (\alpha_j \vert \mu_i)} F_j^{\,(p)} M_j  \cr
   \left( M_j \, ; \, c \atop t \right) E_j^{\,(p)} = E_j^{\,(p)} \left( M_j \,
; \, c + p \, (\alpha_j \vert \mu_i) \atop t \right)  \cr
   \left( M_j \, ; \, c \atop t \right) F_j^{\,(p)} = F_j^{\,(p)} \left( M_j \,
; \, c - p \, (\alpha_j \vert \mu_i) \atop t \right)  \cr
   E_i^{\,(r)} E_i^{\,(s)} = {\left[ r+s \atop r \right]}_{q_i} E_i^{\,(r+s)} \,
, \qquad  F_i^{\,(r)} F_i^{\,(s)} = {\left[r+s \atop r \right]}_{q_i}
F_i^{\,(r+s)}  \cr
   \sum_{r+s=1-\aij} {(-1)}^s E_i^{\,(r)} E_j E_i^{\,(s)} = 0 \, ,  \;
\sum_{r+s=1-\aij} {(-1)}^s F_i^{\,(r)} F_j F_i^{\,(s)} = 0 \, ,  \quad  \forall
\; i \neq j  \cr
   E_i^{\,(0)} = 1 \, ,  \quad  E_i^{\,(r)} F_j^{\,(s)} = F_j^{\,(s)}
E_i^{\,(r)} \, ,  \; \quad  F_i^{\,(0)} = 1  \cr }  $$
   \indent   Then  $ \gerUMh{M} $  is the completion of  $ \hat{\frak
H}_{\scriptscriptstyle M} $  obtained by taking formal series in the
PBW monomials of  $ \gerUm $  and  $ \gerUp $,  with coefficients
in  $ \gerUzM{M} $,  which satisfy the "growth condition" in
Definition 5.8.  Finally, formulas in \S 4.13 yield
%
%
the following (with  $ \, K_i := L_{\alpha_i} \, $):
  $$  \displaylines{
   {\ }  \Delta \big( F_i \big) \equiv F_i \otimes 1 + 1 \otimes F_i + (q_i - 1)
\cdot \bigg( {K_i \, ; \, 0 \atop 1} \bigg) \otimes F_i +   \hfill  \cr
   \hfill   + \, {\left( q_i - q_i^{-1} \right)}^{-1} \cdot \!\!\!\!\!\!\!
                  \sum_{\Sb  \alpha, \beta \in \phitildep  \\
                             p(\alpha) - p(\beta) = -\alpha_i  \\  \endSb}
\hskip-15pt  C^{i,+}_{\alpha,\beta} \, \big( q_\alpha - q_\alpha^{\,-1} \big)
\big( q_\beta - q_\beta^{\,-1} \big) \cdot K_i E_\alpha \otimes F_\beta
\hfill   \mod \, {\left( q - \qm \right)}^2  \cr
   \Delta \bigg( \bigg( {L_{\mu_i} \, ; \, 0 \atop 1} \bigg) \bigg) \equiv
\left( L_{\mu_i} \, ; \, 0 \atop 1 \right) \otimes 1 + 1 \otimes \left(
L_{\mu_i} \, ; \, 0 \atop 1 \right) + (q_i - 1) \cdot \left( L_{\mu_i} \, ; \, 0
\atop 1 \right) \otimes \left( L_{\mu_i} \, ; \, 0 \atop 1 \right) +   \hfill
\cr
   \hfill   + \; {(2)}_{\qm}^{\,2} {(d_i)}_q^{\,-1} \cdot \!\!\!\! \sum_{\Sb
\alpha, \beta \in \phitildep  \\
                p(\alpha) - p(\beta) = 0  \\   \endSb}
\!\! (q - 1) \, C_{\alpha,\beta} \, {[d_\gamma]}_q {\big[ (\mu_i \vert \gamma)
\big]}_q \cdot L_{\mu_i} E_\alpha \otimes F_\beta L_{\mu_i}   \hfill   \mod \,
{\left( q - \qm \right)}^2  \cr }  $$
  $$  \displaylines{
   {\ }  \Delta \big( E_i \big) \equiv 1 \otimes E_i + E_i \otimes 1 + (q_i - 1)
\cdot E_i \otimes \bigg( {K_i \, ; \, 0 \atop 1} \bigg) -   \hfill  \cr
   \hfill   - \, {\left( q_i - q_i^{-1} \right)}^{-1} \cdot \!\!\!\!\!\!\!
         \sum_{\Sb  \alpha, \beta \in \phitildep  \\
                    p(\alpha) - p(\beta) = +\alpha_i  \\  \endSb}
\hskip-15pt   C^{i,-}_{\alpha,\beta} \, \big( q_\alpha - q_\alpha^{\,-1} \big)
\big( q_\beta - q_\beta^{\,-1} \big) E_\alpha \otimes F_\beta K_i   \hfill
\mod \, {\left( q - \qm \right)}^2  \cr
   \hfill   S \big( F_i \big) \equiv - q_i^{-2} \cdot F_i K_i^{-1} \; ,  \qquad
S \big( E_i \big) \equiv - q_i^{+2} \cdot K_i^{-1} E_i   \hfill   \mod \,
{\left( q - \qm \right)}^{\phantom{2}}  \cr
   \hfill   S \left( \left( M_j \, ; \, 0 \atop 1 \right) \right) \equiv -
M_j^{-1} \cdot \left( M_j \, ; \, 0 \atop 1 \right)   \hfill   \mod \, {\left( q
- \qm \right)}^{\phantom{2}}  \cr
   \hfill \; \phantom{.}   \epsilon \big( F_i \big) = 0 \; ,  \qquad  \epsilon
\left( \left( M_j \, ; \, 0 \atop 1 \right) \right) = 0 \; ,  \qquad  \epsilon
\big( E_i \big) = 0  \; .   \hfill   \phantom{\mod \,  \left( q - \qm \right)}
\cr }  $$

\vskip7pt

   {\bf  5.11  Quantum Poisson pairing.}  \;  Since  $ \, {j_{\scriptscriptstyle
M}}^{\! -1} \smallcirc \nu_{\scriptscriptstyle M} \colon \, {\uqMh M }
{\buildrel \cong \over \longrightarrow} {\uqMg {M'} }^* \, $,  via evaluation we
get a perfect Hopf pairing
  $$  \pi_q^{\scriptscriptstyle M} \colon \, {\uqMh M } \otimes
{\uqMg {M'} } \llongrightarrow \Cq  $$
defined by \  $ \pi_q^{\scriptscriptstyle M} (h,g) := \left\langle
{j_{\scriptscriptstyle M}}^{\! -1} \big( \nu_{\scriptscriptstyle M}
(h) \big), g \right\rangle $  \  for all  $ \, h \in
{\uqMh M } $,  $ g \in {\uqMg {M'} } \, $.
                                                      \par
   We call  $ \pi_q^{\scriptscriptstyle M} $  {\it quantum Poisson pairing}.
                                                      \par
   By previous analysis, the integer forms of quantum enveloping algebras are
$ R $--dual  of each other (cf.~\S 2.7) with respect to
$ \pi_q^{\scriptscriptstyle M} \, $;  so the latter restrict to perfect Hopf
pairings
  $$  \pi_{q, \widehat{H}_{\scriptscriptstyle M}} : \, {\gerUMh
{M'} } \otimes {\calUMg M } \loongrightarrow R \, ,  \qquad
\pi_{q, \widehat{G}_{\scriptscriptstyle M}^\infty} : \, {\calUMh M }
\otimes {\gerUMg {M'} } \loongrightarrow R \; .  $$

\vskip1,7truecm

 \centerline{ \bf  \S \; 6 \,  Specialization at roots of 1 }

\vskip10pt

   {\bf 6.1  The case  $ \, q \rightarrow 1 \, $:  specialization of
$ \gerUMh{M} $  to  $ \uh $  and consequences.}  \;  To begin with, set
  $$  \gerUunoMh{M} := \, {\gerUMh M } \Big/ (q-1) \, {\gerUMh M }
\, \cong \, {\gerUMh M } \otimes_R \C  $$
and let  $ \, p_1 \colon \, {\gerUMh M } \rightarrow
\gerUunoMh{M} \, $  be the canonical projection; then set  $ \, \text{f}_i :=
p_1 \left( F_i^{\,(1)} \right) $,  $ \, \text{m}_i := p_1 \left( \left( M_j ;
\, 0 \atop 1 \right) \right) $,  $ \, \text{e}_i := p_1 \left( E_i^{\,(1)}
\right) $,  (where  $ \, M_j := L_{\mu_j} \, $)  for all  $ \, i \in I $,  $ j
\in I_\infty \, $.

\vskip7pt

\proclaim{Theorem 6.2}  For  $ \, q \rightarrow 1 \, $,
$ \, {\gerUMh M } $  specializes to the Poisson Hopf coalgebra
$ \uh \, $;  in other words,
there exists an isomorphism of Poisson Hopf coalgebras
  $$  \gerUunoMh{M} \cong \uh \, .  $$
\endproclaim

\demo{Proof}  The proof mimick that for  $ \, \gerUunoMg{M} \cong \ug \, $.
From the very definition of  $ {\gerUMh M } $  we get  $ \, \gerUunoMh{M} =
\hat{\frak H}_{\scriptscriptstyle M}{\Big\vert}_{q=1} :=
\hat{\frak H}_{\scriptscriptstyle M} \Big/ (q-1) \, \hat{\frak
H}_{\scriptscriptstyle M} \, $,  hence we are reduced to study
$ \, \hat{\frak H}_{\scriptscriptstyle M}{\Big\vert}_{q=1} \, $;
 \eject
\noindent   moreover the presentation of
$ \hat{\frak H}_{\scriptscriptstyle M} $  provides
one of  $ \, \hat{\frak H}_{\scriptscriptstyle
M}{\Big\vert}_{q=1} \, $.  Now the definition of
$ \hat{\frak H}_{\scriptscriptstyle M}{\Big\vert}_{q=1} $
and the explicit form of the PBW basis of  $ {\gerUzM M } $
(cf.~\S 2.6) imply that the elements  $ \, F_i^{\,(r)} $,
$ \left( M_j ; \, 0 \atop t \right) $,  $ M_j^{-1} $,
$ E_i^{\,(s)} \, $  ($ \, i \in I \, , j \in I_\infty \, ; \, r, t, s \in \N
\, $)  are enough to generate  $ \hat{\frak H}_{\scriptscriptstyle M} \, $;
finally, straightforward computation gives  $ \; p_1 \left( F_i^{\,(r)}\right) =
{{\ \text{f}_i^{\,r} \ } \over {\ r! \ }} \, $,  $ \; p_1 \left( \left( M_j ;
\, 0 \atop t \right) \right) = {\text{m}_j \choose t} \, $,  $ \; p_1 \left(
M_j^{-1} \right) = 1 \, $,  $ \; p_1 \left( E_i^{\,(s)} \right) = {{\
\text{e}_i^{\,s} \ } \over {\ s! \ }} \; $  (where  $ \, {\text{m}_i \choose t}
:= {\text{m}_i (\text{m}_i - 1) (\text{m}_i - 2) \cdots (\text{m}_i - t + 1)
\over t!} \, $),  hence  $ \, \gerUunoMh{M} = \hat{\frak H}_{\scriptscriptstyle
M} {\Big\vert}_{q=1} \, $  is generated by the  $ \text{f}_i $'s,  the
$ \text{m}_i $'s,  and the  $ \text{e}_i $'s,  with some relations.
                                                 \par
  When  $ \, M = Q_\infty \, $  this presentation is exactly the same
of  $ \uh $  (cf.~(1.1)), with  $ \, \hbox{h}_i = \hbox{m}_i \, $;
in addition, comparing (1.2) with formulas in \S 5.10  (for  $ q=1 $)
shows that also the Hopf structure is the same.  In particular
$ \gerUunoMh{Q_\infty} $  is cocommutative, hence has a canonical
co-Poisson structure, given by  $ \; \delta := {\, \Delta -
\Delta^{op} \, \over \, q-1 \,} \; $,  described by formulas
--- deduced from those in \S 5.10 ---   which do coincide with
(1.3), as a straightforward checking shows.
                                                 \par
   Finally, for  $ \, M \neq Q_\infty \, $  we prove that  $ \, \gerUunoMh{M}
\cong \gerUunoMh{Q_\infty} \, $  as Poisson Hopf coalgebras: since
$ \, \gerUMh{M} \supseteq \gerUMh{Q_\infty} \, $  by definition, it is
enough to check that  $ \, \hat{\frak H}_{\scriptscriptstyle M}
{\Big\vert}_{q=1} = \hat{\frak H}_{\scriptscriptstyle Q_\infty}
{\Big\vert}_{q=1} \, $  as  $ \C $--vector  spaces.  Assume we are in the simply
laced case.  Since  $ \, M_i := L_{\mu_i} \, $  and  $ \, K_j := L_{\alpha_j}
\, $,  it is  $ \, K_j := \prod_{i=1}^n M_i^{\,c_{ij}} \, $,  where
$ \, c_{ij} \in \Z \, $  are such that  $ \, \alpha_j = \sum_{i=1}^n c_{ij}
\mu_i \, $.  Then  $ \; \left( K_j ; \, 0 \atop t \right) {\Big\vert}_{q=1} =
\sum_{i=1}^n c_{ij} \cdot \left( M_i \, ; \, 0 \atop t \right) {\Big\vert}_{q=1}
\; $  so that  $ \, {\gerUzM M }{\Big\vert}_{q=1} = {\gerUzM
{Q_\infty} }{\Big\vert}_{q=1} \, $  follows, whence  $ \, \hat{\frak
H}_{\scriptscriptstyle M}{\Big\vert}_{q=1} = \hat{\frak H}_{\scriptscriptstyle
Q_\infty} {\Big\vert}_{q=1} \, $,  q.e.d.  In the other cases  $ \, M = P_\infty
\, $,  and this argument still works,  {\sl mutatis mutandis},  because  $ \,
\alpha_j = \sum_{i=1}^n \aij \, \omega_i \, $,  hence  $ \, K_j := \prod_{i=1}^n
L_{\alpha_i}^{\,\aij} \, $,  so that  $ \; \left( K_j ; \, 0 \atop t \right)
{\Big\vert}_{q=1} = \sum_{i=1}^n a_{ji} \cdot \left( L_i \, ; \, 0 \atop t
\right) {\Big\vert}_{q=1} \; $  and we are done again.   $ \square $
\enddemo

\vskip7pt

   Theorem 6.2 has an interesting corollary, namely  $ \,
\calUMg{M} @>{\, q \rightarrow 1 \,}>> F \big[
\widehat{H}_{\scriptscriptstyle M} \big] $  (cf.~\S 3.4).  The
original proof in [BK] is lenghty involved, and requires hard
computations; on the contrary, we can deduce this result as an
easy consequence of the previous one:

\vskip7pt

\proclaim{Theorem 6.3}  The Hopf algebra  $ {\calUMg M } $  specializes to
the Poisson Hopf algebra  $ F \big[ \widehat{H}_{\scriptscriptstyle M} \big] $
for  $ \, q \rightarrow 1 \, $,  in other words, there exists an isomorphism of
Poisson Hopf algebras
  $$  \calUunoMg{M} := \, \calUMg{M} \Big/ (q-1) \, \calUMg{M} \, \cong \, F
\big[ \widehat{H}_{\scriptscriptstyle M} \big] \; .  $$
\endproclaim

\demo{Proof}  Since  $ \, {\calUMg M } \, $  is perfectly paired with
$ \, {\gerUMh M' } \, $,  we have that  $ \, \calUunoMg{M} \, $  is
perfectly paired with  $ \, \gerUunoMh{M'} \cong \uh \, $:  the latter
is cocommutative, hence the former is commutative.  Then $ \calUunoMg{M} $
is a commutative Hopf algebra over  $ \C $,  with countably many generators:
hence it is the algebra of regular functions of a complex affine proalgebraic
group, say  $ \widehat{H}' \, $;  moreover  $ \, \calUunoMg{M} = F \big[
\widehat{H}' \big] \, $  inherits from  $ \calUMg{M} $  a Poisson structure, so
$ \widehat{H}' $  is a Poisson (proalgebraic) group.  It is clear from the
presentation of  $ \calUMg{M} $  that  $ \, F \big[ \widehat{H}' \big] \, \big(
= \calUunoMg{M} \,\big) \cong F \big[ \widehat{H}_{\scriptscriptstyle M} \big]
\, $  as Hopf algebras, hence  $ \, \widehat{H}' =
\widehat{H}_{\scriptscriptstyle M} \, $  as proalgebraic
groups (the non-trivial part, as in [BK], is about the Poisson
structure!).  Now the Hopf pairing among  $ \, \gerUunoMh{M'} \cong
\uh \, $  and  $ \, \calUunoMg{M} = F \big[ \widehat{H}' \big] =
F \big[ \widehat{H}_{\scriptscriptstyle M} \big] \, $  is compatible
with Poisson and co-Poisson structures, that is  $ \, \big\langle h, \{f,g\}
\big\rangle = \big\langle \delta(h), f \otimes g \big\rangle \, $,  where
$ \delta $  is the Poisson cobracket of  $ \, \gerUunoMh{M'} = \uh \, $  and
$ \{\ ,\ \} $  is either the Poisson bracket  $ {\{\ ,\ \}}_\star $  of
$ \widehat{H}_{\scriptscriptstyle M} $  or the Poisson bracket
$ {\{\ ,\ \}}_\circ $  of  $ \widehat{H}' \, $:  since the pairing is perfect,
we must have  $ \, {\{\ ,\ \}}_\star = {\{\ ,\ \}}_\circ \, $,  q.e.d.
$ \square $
\enddemo

\vskip7pt

   {\bf  6.4  The case  $ \, q \rightarrow 1 \, $:  specialization of
$ \calUMh{M} $  to  $ F^\infty \big[ \widehat{G}_{\scriptscriptstyle M}
\big] $.}  \;  We are going to show that  $ \calUMh{M} $  is a quantization of
$ F^\infty \big[ \widehat{G}_{\scriptscriptstyle M} \big] $;  such a result can
be seen as (Poisson) dual counterpart of  $ \, \calUMg{M} \, @>{\, q \rightarrow
1 \,}>> \, F \big[ \widehat{H}_{\scriptscriptstyle M} \big] \, $  (cf.~Theorem
6.3).  As usual, we set
  $$  \calUunoMh{M} := \, {\calUMh M } \Big/ (q-1) \, {\calUMh M } \, \cong \,
{\calUMh M } \otimes_R \C \, .  $$

\vskip7pt

\proclaim{Theorem 6.5}  The formal Hopf algebra
$ {\calUMh M } $  specializes to the formal Poisson Hopf algebra
$ F^\infty \left[ \widehat{G}_{\scriptscriptstyle M} \right] $  for  $ \, q
\rightarrow 1 \, $;  in other words, there exists an isomorphism of formal
Poisson Hopf algebras
  $$  \calUunoMh{M} \cong F^\infty \big[ \widehat{G}_{\scriptscriptstyle M}
\big] \, .  $$
\endproclaim

\demo{Proof}  Recall that  $ F^\infty \big[ \widehat{G}_{\scriptscriptstyle M}
\big] $  is isomorphic to the linear dual of  $ \ug $,  that is  $ \, F^\infty
\big[ \widehat{G}_{\scriptscriptstyle M} \big] \cong {\ug}^* \, $.  On the other
hand, we have an isomorphism  $ \, {j_{\scriptscriptstyle M}}^{\! -1} \smallcirc
\nu_{\scriptscriptstyle M} : \, {\uqMh M } @>{\;\cong\;}>> {\uqMg {M'} }^* \, $
of formal Hopf algebras, and Theorem 5.6 ensures that it restricts to
  $$  {j_{\scriptscriptstyle M}}^{\! -1} \smallcirc
\nu_{\scriptscriptstyle M} : \, \calUMh{M} \, {\buildrel \cong
\over \llongrightarrow} \, {\gerUMg {M'} }^* \, .   \eqno (6.1)  $$
   \indent   When  $ \, q \longrightarrow 1 \, $,  we have that
$ {\gerUMg {M'} } $  specializes to  $ \ug $,
                              therefore (6.1) implies\break
$ \, \calUunoMh{M} \cong {\gerUMg {M'} }^* \otimes_R \C = {\gerUunoMg{M'} }^*
\cong {\ug}^* = F^\infty \big[ \widehat{G}_{\scriptscriptstyle M} \big] =
F^\infty \big[ \widehat{G}_{\scriptscriptstyle M} \big] \, $,  q.e.d.
$ \square $
\enddemo

\vskip7pt

   {\bf 6.6  The case  $ \, q \rightarrow \varepsilon \, $:
quantum Frobenius morphisms.}  \;  Let  $ \varepsilon $  be a primitive
$ \ell $--th  root of 1 in  $ k $,  with  $ \ell $  {\it odd\/}  satisfying
(3.8) and such that  $ \, \ell > d:= \max_i \{d_i\} \, $,  and set
  $$  \gerUepsilonMh{M} := \, {\gerUMh M } \Big/ (q - \varepsilon) \,
{\gerUMh M } \, \cong \, {\gerUMh M } \otimes_R \C  $$
   \indent   The next result is the analogue of (3.5) for  $ {\uqMh M } $.

\vskip7pt

\proclaim{Theorem 6.7}  There exists a continuous epimorphism of
formal Hopf algebras
  $$  \gerFrh \colon \, \gerUepsilonMh{M} \llongtwoheadrightarrow
\gerUunoMh{M} \cong \uh  $$
defined by (for all  $ \, i \in I $,  $ j \in I_\infty \, $)
  $$  {} \!\!\!  \gerFrh \colon
 \cases
   {F_i}^{(s)} \Big\vert_{q=\varepsilon} \!\!\!\! \mapsto
\! {F_i}^{(s / \ell)} \Big\vert_{q=1},  \, \left( M_j ; \,
\, 0 \atop s \right) \! \bigg\vert_{q=\varepsilon} \!\!\!\! \mapsto \!
\left( M_j ; \, \, 0 \atop s / \ell \right) \! \bigg\vert_{q=1},
\, {E_i}^{(s)} \Big\vert_{q=\varepsilon} \!\!\!\! \mapsto \!
{E_i}^{(s / \ell)} \Big\vert_{q=1}  \;\; \hbox{ if } \; \ell
\Big\vert s  \\
   {F_i}^{(s)} \Big\vert_{q=\varepsilon} \!\!\! \mapsto 0 \, ,  \quad
\left( M_j ; \, \, 0 \atop s \right) \! \bigg\vert_{q=\varepsilon}
\!\!\! \mapsto 0 \, ,  \quad {E_i}^{(s)} \Big\vert_{q=\varepsilon}
\!\!\! \mapsto 0  \quad  \hbox{ \, otherwise \, }  \\
   {M_i}^{-1} \Big\vert_{q=1} \!\!\! \mapsto 1  \\
 \endcases  $$
which is adjoint of  $ \calFrg $  (cf.~(3.9)) with respect to the
quantum Poisson pairings, i.e.
  $$  \qquad  \pi_{1, H_{\scriptscriptstyle {M'}}} \Big( \gerFrh(h)
\, , \, g \Big) = \pi_{\varepsilon, H_{\scriptscriptstyle {M'}}}
\Big( h \, , \, \calFrg(g) \Big)   \qquad \quad  \forall \;\;
h \in {\gerUepsilonMh M }, \; g \in \calUepsilonMg{M'} \, .  $$
\endproclaim

\demo{Proof}  The formulas above uniquely determine a continuous
Hopf algebra epimorphism  $ \, \gerFrh \, $,  if any, because
$ \, F_i^{\,(s)} \Big\vert_{q=\varepsilon} $,  $ \left( M_j ;
\, 0 \atop s \right) \! \bigg\vert_{q = \varepsilon} $,
$ M_j^{-1}\Big\vert_{q=\varepsilon} $,  $ E_i^{\,(s)}
\Big\vert_{q=\varepsilon} \, $  are topological generators
of the algebra  $ \gerUepsilonMh{M} $.  Now consider the embedding
$ \; \calFrg : \, F \big[ \widehat{H}_{\scriptscriptstyle M'} \big]
\! \cong \! \calUunoMg{M'} \longhookrightarrow \calUepsilonMg{M'} \; $
of Hopf algebras (cf.~(3.9)): its dual is a (continuous) epimorphism
of formal Hopf algebras  $ \, {\calUepsilonMg{M'} }^*
\longtwoheadrightarrow {\calUunoMg{M'} }^* \, $.  On the other
hand we have an embedding  $ \, \gerUepsilonMh{M} \longhookrightarrow
{\calUepsilonMg{M'} }^* \, $  provided by the specialized quantum
Poisson pairing  $ \; \pi_{\varepsilon,H_{\scriptscriptstyle {M'}}}
: \, \gerUepsilonMh{M} \otimes \calUepsilonMg{M'} \longrightarrow
\C \; $:  therefore composition yields a morphism  $ \, \gerFrh :
\, \gerUepsilonMh{M} \llongrightarrow {\calUunoMg{M'} }^* \, $.
Furthermore, the very construction gives  $ \, \big\langle \gerFrh(h)
\, , \, g \big\rangle = \pi_{1, H_{\scriptscriptstyle {M'}}} \!
\big( \gerFrh(h) \, , \, g \big) = \pi_{\varepsilon,
H_{\scriptscriptstyle {M'}}} \! \big( h \, , \, \calFrg(g) \big)
\, $,  hence  $ \gerFrh $  is adjoint of  $ \calFrg(g) $,  and is
described by the formulas above.  Finally, we have to prove that
$ \gerFrh $  has image  $ \gerUunoMh{M} $.
                                               \par
  The problem is that  $ \gerUepsilonMh{M} $  is made of series,
say  {\sl infinite}  linear combinations   --- of PBW monomials ---
whereas  $ \gerUunoMh{M} $  is made only of  {\sl finite}  linear
combinations   --- of (the same type of) PBW monomials: so we must
show that these  {\sl infinite}  linear combinations indeed are
mapped by  $ \gerFrh $  to  {\sl finite}  ones.
                                               \par
  Recall the definition of  $ \gerUMh{M} $  (see Definition
5.8{\it (b)\/}),  in particular condition (5.3): when  $ \, q =
\varepsilon \, $  (a root of one of odd order  $ \ell $  satisfying
(3.8))  the right-hand-side in (5.3) vanishes as soon as  $ \varphi $
and  $ \eta $  are such that  $ \, f_\beta \geq \ell \, $  or
$ \, e_\gamma \geq \ell \, $  for some  $ \, \beta, \gamma
\in \phipre $  (for in this case one of the factors in the
upper line is zero); and, of course, almost all the PBW
monomials  $ \gerF_\varphi $  and  $ \gerE_\eta $  which
appear in the summands  $ \, \gerF_\varphi \cdot \phi_{\varphi,\eta}
\cdot \gerE_\eta \, $  of the series expansion of an element of
$ \gerUepsilonMh{M} $  do satisfy one of (or both) the conditions
$ \, f_\beta \geq \ell $,  $ e_\gamma \geq \ell \, $  for all
$ \, \beta, \gamma \in \phipre $.  Now let  $ \, (r \delta,i) \in \phitildepim \, $;  let  $ \,\dot{\fbar}_{(r
\delta,i)} \, $  be the element of  $ \calUm $  which is dual (with respect to
$ \pi $  or to  $ \overline{\pi} $)  of  $ \, {\, r \, \over \, {[r]}_{q_i} \,}
\, E_{(r \delta,i)} \, $:  by  {\sl Claim {\it (a)}}  in \S 2.7 we know that
$ \,\dot{\fbar}_{(r \delta,i)} \, $  is an  $ R $--linear  combination of the
$ \, \fbar_{(r \delta,j)} \, $'s,  $ \, j \in I_0 \, $.  Now, by (3.11) and the
{\sl Claim}  in \S 2.3 we have that  $ \calFrg $  maps  $ \dot{\fbar}_{(r
\delta, i)} $  to  $ \ell \,\dot{\fbar}_{(r \ell \delta,i)} \, $:  the dual of
the latter in  $ \gerUepsilonMh{M} $  is  $ \, \ell^{-1} \widehat{E}_{(r \ell
\delta,i)} = \ell^{-1} \, {\, r \ell \, \over \, {[r \ell]}_{q_i} \,} \, E_{(r
\ell \delta,i)} \, $.  Therefore  $ \calFrg $  maps a PBW monomial in the
$ \fbar_\gamma $'s  ($ \gamma \in \phipre $)  and the  $ \dot{\fbar}_{(r
\delta,i)} $  ($ (r \delta,i) \in \phitildepim $),  call it
$ \,\dot{\!\calF}_\varphi $,  to a similar PBW monomial
$ \,\dot{\!\calF}_{\varphi'} $  such that  $ \, \varphi' = {(f_\alpha)}_{\alpha
\in \phitildep} \, $  with  $ \, f_{(r \delta,i)} = 0 \, $  for all  $ \, i \in
I_0 \, $  and  $ \, r \in \N \setminus \ell \, \N \, $.  Similarly occurs when
reverting the roles of the  $ F $'s  and the  $ E $'s.  This fact clearly
implies that, given a PBW monomial  $ g $  of  $ {\calUunoMg M' } $,  we have
$ \, \big\langle \gerFrh(h) \, , \, g \big\rangle = \pi_{1,
H_{\scriptscriptstyle {M'}}} \! \big( \gerFrh(h) \, , \, g \big) =
\pi_{\varepsilon, H_{\scriptscriptstyle {M'}}} \! \big( h \, , \, \calFrg(g)
\big) = 0 \, $  for all PBW monomials  $ \, h = \gerF_\varphi \cdot u_\tau \cdot
\gerE_\eta \, $  of  $ \hat{\frak H}_{\scriptscriptstyle
M}{\big\vert}_{q=\varepsilon} $  ($ \subset \gerUepsilonMh{M} $)  such that
$ \, f_{(r \delta,i)} \neq 0 \, $  or  $ \, e_{(r \delta,i)} \neq 0 \, $  for
some  $ \, i \in I_0 \, $  and  $ \, r \in \N \setminus \ell \, \N \, $:  thus
$ \gerFrh $  maps to zero this kind of PBW monomials of  $ \gerUepsilonMh{M} $.
As for the other ones, we can read the condition  $ \, f_{(r \delta,i)} = 0
\, $  or  $ \, e_{(r \delta,i)} = 0 \, $  for all  $ \, i \in I_0 \, $  and
$ \, r \in \N \setminus \ell \, \N \, $  as
  $$  f_{(r \delta,i)} \neq 0 \; =\joinrel=\joinrel\Rightarrow \; \ell
\,\big\vert\, r \; ,  \qquad  e_{(r \delta,i)} \neq 0 \;
=\joinrel=\joinrel\Rightarrow \; \ell \,\big\vert\, r \; ; $$
then in (5.3) the right-hand-side vanishes again (this time because,
in this case, one of the factors in the lower line is zero).  So the
upshot is that  $ \gerFrh $  kills almost all the summands in the
series expansion of any element of  $ \gerUepsilonMh{M} $,  which
solves our problem.   $ \square $\break
\enddemo

\vskip7pt

\proclaim{Corollary 6.8} The images of imaginary root vectors under the quantum
Frobenius morphisms  $ \, \gerFrh : \gerUepsilonMh{M} \llongtwoheadrightarrow
\gerUunoMh{M} \, $  and  $ \, \gerFrg : \gerUepsilonMg{M}
\llongtwoheadrightarrow \gerUunoMg{M} \, $  are given by (for all  $ \, (r
\delta,i) \in \phitildepim \, $)
  $$  \gerFrh \, , \, \gerFrg :  \cases
        \widehat{F}_{(r \delta,i)} \Big\vert_{q=\varepsilon}
\mapsto \; \ell \, \widehat{F}_{(r/\ell \, \cdot \, \delta, i)}
\Big\vert_{q=\varepsilon} \, ,  \quad  \widehat{E}_{(r \delta,i)}
\Big\vert_{q=\varepsilon} \mapsto \; \ell \,
\widehat{E}_{(r/\ell \, \cdot \, \delta, i)}
\Big\vert_{q=\varepsilon}  \hskip5pt  \hbox{ if }
\hskip9pt  \ell \, \Big\vert \, r  \\
        \widehat{F}_{(r \delta,i)} \Big\vert_{q=\varepsilon} \mapsto \; 0 \; ,
\quad  \hskip55,7pt  \widehat{E}_{(r \delta,i)} \Big\vert_{q=\varepsilon}
\mapsto \; 0  \qquad  \hskip35,3pt  \hbox{ otherwise }  \\
                                  \endcases  $$
\endproclaim

\demo{Proof}  For  $ \gerFrh $  the result follows from (3.11) by duality
(cf.~the proof of Theorem 6.7); since on quantum unipotent subalgebras
$ \gerFrg $  coincides with  $ \gerFrh $  (they are defined by the same
formulas) the result follows for  $ \gerFrg $  too (otherwise, we can get it
again by duality from (6.4) below, with the same argument used for
$ \gerFrh $).   $ \square $
\enddemo

\vskip7pt

   Similar arguments to those used to prove Theorem 6.7 provide a proof of the
next result, which is the analogue of (3.9); as usual, we set
  $$  \calUepsilonMh{M} := \, \calUMh{M} \Big/ (q - \varepsilon) \, \calUMh{M}
\, \cong \, \calUMh{M} \otimes_R \C \, ;  $$
moreover we define the set
  $$  \N^{\phitildep}[\ell] := \Big\{\, \chi = {(x_\alpha)}_{\alpha
\in \phitildep} \! \in \! \N^{\phitildep} \;\Big\vert\; x_\gamma \!
\in \! \ell \, \N \;\, \forall \, \gamma \! \in \! \phipre, \;\,
x_{(r \delta,i)} = 0 \;\; \forall \, i \! \in \! I_0, \, r \!
\in \! \N \setminus \ell \, \N \,\Big\}  $$

\vskip7pt
\proclaim{Theorem 6.9} \; (a) There exists a unique continuous monomorphism of
formal Hopf algebras
  $$  \calFrh \colon \, F^\infty \big[ \widehat{G}_{\scriptscriptstyle M} \big]
\cong \calUunoMh{M} \llonghookrightarrow \calUepsilonMh{M}   \eqno (6.2)  $$
defined (for all  $ \, \alpha \in \phipre \, $,  $ \mu \in M \, $)  by
  $$  \calFrh :  \; \quad  \fbar_\alpha \Big\vert_{q=1}
\mapsto \fbar_\alpha^{\,\ell} \Big\vert_{q=\varepsilon}
\, , \; \quad L_\mu \Big\vert_{q=1} \mapsto
L_\mu^{\,\ell}\Big\vert_{q=\varepsilon} \, ,  \; \quad \ebar_\alpha
\Big\vert_{q=1} \mapsto \ebar_\alpha^{\,\ell}
\Big\vert_{q=\varepsilon}   \eqno (6.3)  $$
and enjoying
  $$  \calFrg : \, \; \fbar_{(r \delta,i)} \Big\vert_{q=1} \mapsto \ell \,
\fbar_{(r \ell \delta,i)} \Big\vert_{q=\varepsilon} \, ,  \;\, \ebar_{(r
\delta,i)} \Big\vert_{q=1} \mapsto \ell \, \ebar_{(r \ell \delta,i)}
\Big\vert_{q=\varepsilon}  \quad  \forall \; (r \delta,i) \in \phitildepim
\eqno (6.4)  $$
which is adjoint of  $ \gerFrg $  (cf.~(3.5)) with respect to quantum Poisson
pairings, that is
  $$  \qquad  \pi_{\varepsilon,\widehat{G}_{\scriptscriptstyle M}} \! \Big(
\calFrh(h) \, , \, g \Big) = \pi_{1,\widehat{G}_{\scriptscriptstyle M}} \! \Big(
h \, , \, \gerFrg(g) \Big)   \quad \qquad  \forall \;\; h \in \calUunoMh{M}
\, , \; g \in \gerUepsilonMg{M'} \, .  $$
                                            \hfill\break
   \indent   (b) The image  $ \, Z_0 \; \left( \, \cong_{\calFrh} \calUunoMh{M}
\, \right) \, $  of  $ \calFrh $  is a formal Hopf subalgebra contained in the
centre of  $ \calUepsilonMh{M} $.
                                            \hfill\break
   \indent   (c) The set  $ \Big\{\, \calF_\phi \cdot B_\tau
\cdot \calE_\eta \,\Big\vert\, \phi, \eta \in \N^{\phitildep}[\ell]
\, ,  \tau \in \ell \, \N^{I_\infty} \Big\} $  is a pseudobasis of
$ Z_0 $  over  $ \C $.
                                            \hfill\break
   \indent   (d) The set  $ \, \Big\{\, \calF_\phi \cdot B_\tau \cdot \calE_\eta
\; \Big\vert \; \phi, \eta \in \N^{\phitildep} \setminus \N^{\phitildep}[\ell]
\, ,  \tau \in \N \setminus \ell \, \N^{I_\infty} \,\Big\} \, $  is a
pseudobasis of  $ \calUepsilonMh{M} $  over  $ Z_0 \, $;  therefore also the set
of ordered PBW monomials  $ \, \Big\{\, \calF_\phi \cdot L_\mu \cdot \calE_\eta
\; \Big\vert \; \phi, \eta \in \N^{\phitildep} \setminus \N^{\phitildep}[\ell]
\, , \mu \in \N \setminus \ell \, \N^{I_\infty} \,\Big\} \, $  is a
$ Z_0 $--basis  of  $ \calUepsilonMh{M} $,  so  $ \calUepsilonMh{M} $  is a
free  $ Z_0 $--module.
\endproclaim

\demo{Proof}  Like in [BK], Lemma 3.2.2, we exploit the fact that  $ \,
\fbar_\alpha \Big\vert_{q=1} $,  $ L_\mu \Big\vert_{q=1} $,  $ \ebar_\alpha
\Big\vert_{q=1} $  ($ \alpha \in \phipre $,  $ \mu \in M \, $)  are generators
of  $ \calUunoMh{M} $  as a formal Poisson Hopf algebra (the imaginary root
vectors being obtained from them by means of the Poisson bracket) to get sure
that the formulas above uniquely determine a continuous monomorphism
$ \calFrh $,  if any.  Now consider  $ \, \gerFrg \colon \, \gerUepsilonMg{M'}
\longtwoheadrightarrow \gerUunoMg{M'} \cong \ug \, $  (cf.~(3.5)): its dual is a
formal Hopf monomorphism  $ \, {\gerUunoMg{M'} }^* \longhookrightarrow
{\gerUepsilonMg{M'} }^* \, $;  composing the latter with the isomorphisms  $ \,
\calUunoMh{M} @>\cong>> {\gerUunoMg{M'} }^* \, $,  $ \, {\gerUepsilonMg{M'} }^*
@>\cong>> \calUepsilonMh{M} \, $  (given by specialized quantum Poisson
pairings) provides a monomorphism  $ \, \calFrh \colon \, \calUunoMh{M}
\longhookrightarrow \calUepsilonMh{M} \, $,  which by construction is
continuous; moreover,
  $$  \quad  \Big\langle \calFrh(h) \, , \, g \Big\rangle =
\pi_{\varepsilon,\widehat{G}_{\scriptscriptstyle M}} \! \Big(
\calFrh(h) \, , \, g \Big) = \pi_{1, \widehat{G}_{\scriptscriptstyle M}} \!
\Big( h \, , \, \gerFrg(g) \Big)   \;\;\;  \forall \;\; h \in \calUunoMh{M} \, ,
\; x \in \gerUepsilonMg{M'}  $$
hence  $ \calFrh $  is described by formulas above, as one sees at
once just proceeding as for Theorem 6.7, by comparing PBW monomials
--- or simply root vectors ---   on one hand and their duals on the
other: this also shows that (6.4) holds.  So claim  {\it (a)\/}  is
proved.
                                                  \par
   Claim  {\it (b)\/}  follows from the analogous result for
$ \calUepsilonMg{M} $  (cf.~[BK], \S\S 2--3)  and comparison between
$ \calUepsilonMg{M} $  and  $ \calUepsilonMh{M} $.  Claim  {\it (c)\/}
follows from (6.3) and (6.4).  Finally, the span of  $ \,
\big\{\, B_\tau \;\big\vert\; \tau \in \ell \, \N^{I_\infty} \,\big\}
\, $  inside  $ \calUepsilonMh{M} $  equals the span of
$ \, \big\{\, L_\mu \;\big\vert\; \mu \in \ell \, \N^{I_\infty}
\cong \ell \, M_+ \,\big\} \, $;  this and the explicit
description of the pseudobasis of  $ \calUepsilonMh{M} $
give claim  {\it (d)}.   $ \square $
\enddemo

\vskip7pt

  We call also  $ \gerFrh $  and  $ \calFrh $  {\bf quantum Frobenius
morphisms}.

\vskip7pt

   {\bf 6.10  Specializations of quantum Poisson pairings.}  \;
From \S\S 6.3--5 we get that the Hopf pairings  $ \, \pi_{q,
\widehat{H}_{\scriptscriptstyle M}} \colon \, {\gerUMh M' }
\otimes \calUMg{M} \loongrightarrow R \, $,  $ \, \pi_{q,
\widehat{G}_{\scriptscriptstyle M}} \colon \, \calUMh{M} \otimes
{\gerUMg M' } \loongrightarrow R \, $  (cf.~\S 5.11)
respectively specialize to the Hopf pairings  $ \,
\pi_{\widehat{H}_{\scriptscriptstyle M}} \colon \, \uh \otimes
F \big[ \widehat{H}_{\scriptscriptstyle M} \big] \loongrightarrow
\C \, $,  $ \, \pi_{\widehat{G}_{\scriptscriptstyle M}} \colon
\, F^\infty \big[ \widehat{G}_{\scriptscriptstyle M} \big]
\otimes \ug \loongrightarrow \C \, $;  in other words,
$ \, \pi_{q,\widehat{H}_{\scriptscriptstyle M}} \big( {\hat h},
{\tilde g} \big){\big\vert}_{q=1} =  \pi_{\widehat{H}_{\scriptscriptstyle M}}
\big({\hat h}{\big\vert}_{q=1}, {\tilde g}{\big\vert}_{q=1} \big) \, $,
$ \, \pi_{q, \widehat{G}_{\scriptscriptstyle M}} \big( {\tilde h}, {\hat g}
\big){\big\vert}_{q=1} =  \pi_{\widehat{G}_{\scriptscriptstyle M}} \big( {\tilde
h}{\big\vert}_{q=1}, {\hat g}{\big\vert}_{q=1} \big) \, $.  Thus the quantum
Poisson pairing is a quantization of the classical Hopf pairing on both our
Poisson groups   --- maybe formal ---   dual of each other.  In addition
we show that it can also be thought of as a quantization of the
classical Poisson pairing  $ \, \pi_{\Cal P} : \, \gerh \otimes \gerg
\rightarrow \C \, $,  and of  {\sl new}  pairings between function algebras.  We
use notations  $ \, [\ \, ,\ ] := m - m^{op} \, $,  $ \, \nabla := \Delta -
\Delta^{op} \, $  (superscript "$ op $" denoting opposite operation).
                                                      \par
   First of all, we define a suitable grading on  $ \gerUMg{Q} $  (as an
$ R $--module)  by
  $$  deg \left( \, \prod_{\alpha \in \phitildep} \! E_\alpha^{\,(e_\alpha)} \!
\cdot \prod_{j \in I_\infty} \! \left( K_j ; 0 \atop t_j \right) \!
K_j^{-Ent(t_j/2)} \cdot \prod_{\alpha \in \phitildep} \! F_\alpha^{\,(f_\alpha)}
\! \right) \! := \sum_{\alpha \in \phitildep} \big( e_\alpha + f_\alpha \big) +
\sum_{j \in I_\infty} t_j  $$
and linear extension.  Now let  $ \, R =: \gerU_0 \subset \gerU_1 \subset
\cdots \subset \gerU_h \subset \cdots (\, \subset {\gerUMg Q }) \, $  be the
associated filtration, and set  $ \; \partial(x) := h \; $  for all  $ \, x \in
\gerU_h \setminus \gerU_{h-1} \, $;  notice that a similar notion of degree
exists for  $ \ug $,  defined by means of the filtration  $ \, U_0 \subset U_1
\subset \cdots \subset U_N \subset \cdots \subset \ug \, $  induced by the
canonical filtration of  $ T \left( \gerg \right) $  (the tensor algebra on
$ \gerg $),  and similarly for  $ \, \ug \otimes \ug \, $.  Then  $ \,
{(q-1)}^{\partial (g)} g \in \calUMg{Q} \, $  for all  $ \, g \in \gerUMg{Q}
\, $,  thus
  $$  {} \qquad  \pi_{q,{\Cal P}} (h,g) := {(q-1)}^{+\partial (g)} \cdot \pi_q
(h,g)   \eqno  \forall \;\; h \in {\gerUMh Q }, \, g \in \gerUMg{Q}  \qquad  $$
defines a perfect pairing  $ \, \pi_{q,{\Cal P}} \colon \, {\gerUMh Q } \times
{\gerUMg Q } \longrightarrow R \, $,  which can be specialized at  $ q = 1 $.

\vskip7pt

\proclaim{Theorem 6.11}  \  $ \pi_{q,{\Cal P}} : \gerUMh{Q} \times \gerUMg{Q}
\loongrightarrow R $  \  specializes at  $ \, q = 1 \, $  to a pairing
  $$  \pi_{\Cal P} : \uh \times \ug \llongrightarrow \C  $$
which extends the Lie bialgebra pairing  $ \; \langle \ , \ \rangle : \, \gerh
\otimes \gerg \longrightarrow \C \; $  (cf.~\S 1.2) and is such that
  $$  \eqalign{
   \pi_{\Cal P} (\alpha \cdot x + \beta \cdot y \, , \, z)  &  = \alpha \cdot
\pi_{\Cal P} (x,z) + \beta \cdot \pi_{\Cal P} (y,z)  \cr
   \pi_{\Cal P} (x \, , \, \alpha \cdot u + \beta \cdot v)  &  = \alpha \cdot
\pi_{\Cal P} (x,u) + \beta \cdot \pi_{\Cal P} (x,v)  \cr
   \pi_{\Cal P} \big( x \cdot y \, , \, z \big) = \pi_{\Cal P} \big( x \otimes y
\, , \, \Delta(z) \big) \, ,  &  \qquad  \pi_{\Cal P} \big( x \, , \, z \cdot w
\big) = \pi_{\Cal P} \big( \Delta(x) \, , \, z \otimes w \big)  \cr
   \pi_{\Cal P} \big( [x,y] \, , \, z \big) = \pi_{\Cal P} \big( x \otimes y \,
, \, \delta(z) \big) \, ,  &  \qquad  \pi_{\Cal P} \big( x \, , \, [z,w] \big) =
\pi_{\Cal P} \big( \delta(x) \, , \, z \otimes w \big)  \cr }   \eqno (6.5)  $$
for all  $ \, \alpha, \beta \in \C $,  $ x, y \in \uh $,  $ z, w, u, v \in \ug
\, $  such that  $ \, \partial(\alpha \cdot u + \beta \cdot v) = \partial(u)
= \partial(v) \, $.
\endproclaim
 \eject

\demo{Proof}  Let  $ \, x \in \uh $,  $ z \in \ug $,  and pick  $ \, x'
\in \gerUMh{Q} $,  $ z' \in \gerUMg{Q} $,  such that  $ \, x =
x'{\big\vert}_{q=1} $,  $ z = z'{\big\vert}_{q=1} $.  By definition,
$ \pi_{\Cal P} (x,z) $  is given by
  $$  \pi_{\Cal P} (x,z) := {\pi_{q,{\Cal P}} \big( x', z'
\big)}{\Big\vert}_{q=1} = {\Big( {(q-1)}^{\partial(z')} \cdot \pi_q \big( x', z'
\big) \Big)}{\Big\vert}_{q=1} \; ;  $$
in particular, we can select  $ x' $  and  $ z' $  such that  $ \, \partial
\big( x' \big) = \partial(x) \, $,  $ \, \partial \big( z' \big) = \partial(z)
\, $.  Now, the first two lines in (6.5) follow directly from similar
properties for  $ \, \pi_{q,{\Cal P}} \, $,  which are implied by
definitions.  Second, definitions and Leibnitz' and co-Leibnitz' rules imply
  $$  \displaylines{
   \partial(x \cdot y) = \partial(x) + \partial(y) = \partial(x \otimes y)
\, ,  \qquad  \partial \big( \Delta(x) \big) = \partial(x)  \cr
   \partial \big( \delta(x) \big) = \partial(x) + 1 \, ,  \qquad  \partial \big(
[x,y] \big) = \partial(x) + \partial(y) - 1  \cr }  $$
for all  $ \, x, y \in \ug \, $  provided that  $ \, [x,y] \neq 0 \, $  (to be
complete we may set  $ \, \partial(0) := - \infty \, $);  using these identities
and Leibnitz' and co-Leibnitz' rules we easily reduce to prove that the
remaining identities in (6.5) do hold in degree 1, i.e.~for  $ \, x, y \in \hhat
\, $  and  $ \, z, w \in \ghat \, $:  but this again follows from definition.
Finally to prove that  $ \pi_{\Cal P} $  is an extension of the classical
Poisson pairing it is enough to perform a computation on Chevalley generators,
which is completely straightforward.   $ \square $
\enddemo

\vskip7pt

   {\bf 6.12  The pairing  $ \, F^\infty \big[ \widehat{G}_{\scriptscriptstyle
M'} \big] \times F \big[ \widehat{H}_{\scriptscriptstyle M} \big]
\longrightarrow \C \, $.}  \; The construction in \S 6.10 can be reversed as
follows.  Define a grading on  $ \calUMg{M} $  (as a  $ R $--module)  by
  $$  deg \left( \, \prod_{\alpha \in \phitildep} \ebar_\alpha^{\,e_\alpha}
\cdot \prod_{j \in I_\infty} \! {\left( M_j^{\pm 1} - 1 \right)}^{m_j} \cdot
\prod_{\alpha \in \phitildep} \fbar_\alpha^{\,f_\alpha} \right) := \sum_{\alpha
\in \phitildep} \big( e_\alpha + f_\alpha \big) + \sum_{j \in I_\infty} m_j  $$
and linear extension; then let  $ \, R =: \calU_0 \subset \calU_1 \subset \cdots
\subset \calU_h \subset \cdots (\, \subset {\calUMg M }) \, $  be the associated
filtration, and set  $ \; \partial(x) := h \; $  for all  $ \, x \in \calU_h
\setminus \calU_{h-1} \, $  ($ h \in \N $).  Finally define
  $$  {} \qquad  \pi_q^{\Cal P} (h,g) := {(q-1)}^{-\partial (g)} \cdot
\pi_q(h,g)   \eqno  \forall \;\; h \in {\calUMh M' }, \, g \in {\calUMg M }
\; ;  \qquad  $$
this yields a perfect pairing  $ \, \pi_q^{\Cal P} : {\calUMh M' } \times
{\calUMg M } \longrightarrow R \, $,  to be specialized at  $ \, q = 1 \, $.

\vskip7pt

\proclaim{Teorema 6.13} \  $ \pi_q^{\Cal P} : \calUMh{M'} \times \calUMg{M}
\loongrightarrow R \; $  specializes at  $ \, q = 1 \, $  to a pairing
  $$  \pi^{\Cal P} : \, F^\infty \big[ \widehat{G}_{\scriptscriptstyle M'} \big]
\times F \big[ \widehat{H}_{\scriptscriptstyle M} \big] \llongrightarrow \C  $$
such that
  $$  \eqalign{
   \pi_\tau^{\Cal P} (\alpha \cdot x + \beta \cdot y \, , \, z)  &  = \alpha
\cdot \pi_\tau^{\Cal P} (x,z) + \beta \cdot \pi_\tau^{\Cal P} (y,z)  \cr
   \pi_\tau^{\Cal P} (x \, , \, \alpha \cdot u + \beta \cdot v)  &  = \alpha
\cdot \pi_\tau^{\Cal P} (x,u) + \beta \cdot \pi_\tau^{\Cal P} (x,v)  \cr
   \pi_\tau^{\Cal P} \big( x \cdot y \, , \, z \big) = \pi_\tau^{\Cal P} \big( x
\otimes y \, , \, \Delta(z) \big) \, ,  &  \qquad  \pi_\tau^{\Cal P} \big( x \,
, \, z \cdot w \big) = \pi_\tau^{\Cal P} \big( \Delta(x) \, , \, z \otimes w
\big)  \cr
   \pi_\tau^{\Cal P} \big( \{x,y\} \, , \, z \big) = \pi_\tau^{\Cal P} \big( x
\otimes y \, , \, \nabla(z) \big) \, ,  &  \qquad  \pi_\tau^{\Cal P} \big( x \,
, \, \{z,w\} \big) = \pi_\tau^{\Cal P} \big( \nabla(x) \, , \, z \otimes w
\big)  \cr }   \eqno (6.6)  $$
for all  $ \, \alpha, \beta \in \C $,  $ x, y \in F^\infty \big[
\widehat{G}_{\scriptscriptstyle M'} \big] $,  $ z, u, v \in F \big[
\widehat{H}_{\scriptscriptstyle M} \big] \, $  such that  $ \, \partial(\alpha
\cdot u + \beta \cdot v) = \partial(u) = \partial(v) \, $  (with  $ \,
\partial(x) := \partial \left( x' \right) \, $  for any  $ \, x' \in \calUMg{M}
\, $  such that  $ \, x'{\big\vert}_{q=1} = x \, $).
\endproclaim

\demo{Proof}  Just mimick the proof of Theorem 6.11 above.   $ \square $
\enddemo

 \eject
%
%
%

   $ \underline{\hbox{\sl Remark}} $:  \, If we extend  $ R $  by
adding a  $ \Delta_\infty $--th  root of unity  ($ \Delta_\infty $
as in \S 1.1) we can perform the previous construction for any pair
of lattices  $ M_1 $,  $ M_2 $  such that  $ \, Q_\infty \leq M_1 ,
M_2 \leq P_\infty \, $,  thus getting corresponding pairings like
the one of Theorem 6.13.

\vskip1,9truecm

\Refs
\endRefs

\vskip7pt

\smallrm

[Be1] \  J.~Beck,  {\smallit Braid group action and quantum affine
algebras\/},  Comm.~Math.~Phys.~{\smallbf 165} (1994), 555--568.

\vskip3pt

[Be2] \  J.~Beck,  {\smallit Convex bases of PBW type for quantum affine
algebras\/},  Comm.~Math.~Phys.~{\smallbf 165} (1994), 193--199.

\vskip3pt

[BK] \  J.~Beck, V.~G.~Kac,  {\smallit Finite dimensional representations of
quantum affine algebras at roots of 1\/},  J.~Amer.~Math.~Soc.~{\smallbf 9}
(1996), 391--423.

\vskip3pt

[Bo] \  N.~Bourbaki,  {\smallit Groupes et alg\`ebres de Lie\/},  Chapitres
4--6, Hermann, Paris, 1968.

\vskip3pt
[CP] \  V.~Chari, A.~Pressley,  {\smallit A guide to Quantum Groups\/},
Cambridge University Press, Cambridge, 1994.

\vskip3pt

[Da] \  I.~Damiani,  {\smallit La  R--matrice  pour les alg\`e{}bres
quantiques de type affine non tordu\/},  Ann.~Scient.~\'Ec. Norm.~Sup.,
4$^e$  s\'erie, {\smallbf 31} (1998), 493--523.

\vskip3pt

[Di] \  J.~Dieudonn\'e,  {\smallit Introduction to the theory of formal
groups\/},  Pure Appl.~Math.~{\smallbf 20}, Marcel Dekker, Inc.~New York, 1973.

\vskip3pt

[DL] \  C.~De Concini, V.~Lyubashenko,  {\smallit Quantum function algebra at
roots of 1\/},  Adv.~Math.~{\smallbf 108} (1994), 205--262.

\vskip3pt

[Dr] \  V.~G.~Drinfel'd,  {\smallit Quantum groups\/},  Proc.~ICM Berkeley 1
(1986), 789--820.

\vskip3pt

[Ga1] \  F.~Gavarini,  {\smallit Quantization of Poisson groups\/},
Pac.~J.~Math.~{\smallbf 186} (1998), 217--266.

\vskip3pt

[Ga2] \  F.~Gavarini,  {\smallit A PBW basis for Lusztig's form of
untwisted affine quantum groups\/},  Commun.~Algebra  {\smallbf 27}
(1999), no.~2, 903--918.

\vskip3pt

[Ga3] \  F.~Gavarini,  {\smallit The quantum duality principle\/},
preprint math.QA/9909071.

\vskip3pt

[Lu1] \  G.~Lusztig,  {\smallit Quantum groups at roots of 1\/},
Geom.~Dedicata~{\smallbf 35} (1990), 89--113.

\vskip3pt

[Lu2] \  G.~Lusztig,  {\smallit Introduction to quantum groups\/},
Progr.~Math.~{\smallbf 110}, Birkh\"auser, Boston, 1993.

\vskip3pt

[Ta] \  T.~Tanisaki,  {\smallit Killing forms, Harish-Chandra Isomorphisms, and
Universal R-Matrices for Quantum Algebras\/},  Internat.~J.~Modern Phys.~A
{\smallbf 7}, Suppl.~1B (1992), 941--961.

\vskip1truecm

{}

\enddocument